%% file: main_revised.tex
\documentclass[3p]{elsarticle}

\usepackage[utf8]{inputenc}

\usepackage{graphicx} 
\usepackage{amssymb}
\usepackage{amsmath}
\usepackage{amsfonts}
\usepackage{amsthm}

\usepackage{subcaption}

\usepackage{aligned-overset}

\newtheorem{theorem}{Theorem}[section]
\theoremstyle{definition}

\newtheorem{lemma}[theorem]{Lemma}
\theoremstyle{remark}
\newtheorem{remark}[theorem]{Remark}

\input{macro_CR.tex}

\makeatletter
\def\ps@pprintTitle{
 \let\@oddhead\@empty
 \let\@evenhead\@empty
 \def\@oddfoot{}
 \let\@evenfoot\@oddfoot}
\makeatother

\newcommand{\revision}[1]{{\color{black}{#1}}}

\newcommand{\revisionT}[1]{{\color{black}{#1}}}

\newcommand{\revisionU}[1]{{\color{black}{#1}}}

\title{Thermodynamically consistent phase-field modeling and numerical simulation for \revision{reactive} two-phase fluid-solid dynamics}

\author[inst]{Cedric Riethm\"uller\corref{cor1}}
\ead{cedric.riethmueller@ians.uni-stuttgart.de}
\author[inst]{Lars von Wolff}
\ead{lars.von-wolff@ians.uni-stuttgart.de}
\author[inst]{Dominik G\"oddeke}
\ead{dominik.goeddeke@ians.uni-stuttgart.de}
\author[inst]{Christian Rohde}
\ead{christian.rohde@ians.uni-stuttgart.de}

\affiliation[inst]{organization={
Institute of Applied Analysis and Numerical Simulation \& Stuttgart Center for Simulation Science (SC SimTech),
University of Stuttgart},
            country={Germany} \\
            (Emails:  $\{\text{cedric.riethmueller, lars.von-wolff, dominik.goeddeke, christian.rohde}\}$@ians.uni-stuttgart.de)
            }
            
\cortext[cor1]{Corresponding author}

\begin{document}

\begin{abstract}
We introduce a coupled Cahn-Hilliard Navier-Stokes model  that governs the \revision{reactive} two-phase dynamics of a  system that consists of a  fluid and a solid phase and prove its thermodynamic consistency. Moreover, we present an associated fully-discrete numerical method that relies on a continuous finite element approach and a semi-implicit time-stepping method. \revisionU{Our main theoretical result establishes that the fully discrete method satisfies an analog of the free energy dissipation inequality, provided that the phase-field variable remains within prescribed bounds.}\\
Numerical experiments confirm the theoretical findings and show the applicability of the method for realistic settings. In this context, we provide a preprocessing strategy that enables computing fluid flow in complex geometries given a sharp-interface formulation of the initial phase distribution. Moreover, we briefly introduce different solution strategies for the novel discretization based on the monolithic and partitioned solution paradigms and assess these in a comparative study.
\end{abstract}

\begin{keyword}
Phase-field model \sep Cahn-Hilliard \sep Navier-Stokes \sep Thermodynamic consistency \sep Solution strategies \sep Precipitation/dissolution 
\end{keyword}

\maketitle

\section{Introduction}

We introduce a phase-field model for \revision{reactive} fluid-solid single-phase flow comprising the Navier-Stokes equations for evaluating the  flow field and the Cahn-Hilliard evolution for calculating the interfacial dynamics of the fluid and the solid phase. The model's intrinsic complexity lies in its coupled multiphysics nature along with the nonlinearity of the subproblems and its nonconvexity.

\revision{For further phase-field models that take precipitation/dissolution of ions into account, we refer to \cite{Xu2008, Bringedal2020, Rohde2021, vanNoorden2011, Lamperti2025}. We further mention the models in \cite{vonWolff2021, Schollenberger2024} that are derived from \cite{Rohde2021}. As in these studies we consider a liquid and a solid phase with fixed solute concentration in the latter. The work \cite{Lamperti2025} provides a comprehensive overview of modeling differences in these works. In this context, we note that in our model fluid flow and the Gibbs-Thomson effect related to the interface curvature are accounted for. Moreover, in \cite{Lamperti2025} an Allen-Cahn--Cahn-Hilliard ansatz is used for the reactive model. One asset of the resulting model is that it has a variational structure, i.e., the model's governing equations are derived from an underlying free-energy functional. This is not the case for our model.
}

A central focus of this work is on considerations related to thermodynamic consistency. Firstly, we prove the latter for the governing two-phase model. \revisionU{Then, we suggest a thermodynamically consistent fully-discrete method, i.e., we show that discrete solutions  obey an energy inequality if the phase-field variable can be controlled appropriately.} 
The literature for unconditionally gradient-stable schemes for Cahn-Hilliard Navier-Stokes models is quite rich, especially in the context of two-phase incompressible flows \revision{\cite{Feng2006, Gruen2014, Guo2014, Tierra2014, Shen2015, Chen2016, Garcke2016, Gruen2016, ShokrpourRoudbari2018, tenEikelder2024, Brunk2026}}. Other representatives of this class of schemes are established for models governing moving contact line problems \cite{Bonart2019}, variable density two-phase incompressible surface flow \cite{Palzhanov2021}, dynamic boundary conditions \cite{Giorgini2023}, reactive multi-species flow \cite{Wang2025}, and three-phase flow \cite{Minjeaud2013}.

\revision{We further reduce the model to cover the case of non-reactive flow. The resulting model describes a viscous fluid phase interacting with an immobile solid phase, which may alter subject to geometric evolution. This is in contrast to phase-field models for fluid-structure interaction \cite{Mokbel2018} where the solid can alter through, e.g., elastic deformation.}

Furthermore, we provide a strategy to use \revision{the non-reactive} model to compute fluid flow in complex geometries given a sharp-interface formulation of the initial phase distribution. Here, the proposed ansatz works as a preprocessing step and hence does not require changes to the model formulation but rather to the parameter configuration. A pivotal point is that the fluid-solid interface is only moving due to geometric effects and thus almost static.
This is in contrast to Cahn-Hilliard Navier-Stokes models developed in the diffuse-domain context \cite{Aland2010, Guo2021} which on the one hand introduce a second phase-field parameter for the domain boundary but on the other hand intrinsically satisfy the staticity requirement. In this context, we also refer to Chapter 4 of \cite{phdvonWolff}.
There also exist other in some sense contrary approaches that allow for changes in the (solid) structure as for example in the field of shape and topology optimization. Here, the underlying mathematical task is to solve a constrained optimization problem in order to minimize a structure-related objective functional subject to PDE constraints. Besides convolution-thresholding methods \cite{Ruuth2001, Chen2024}, phase-field methods are employed to solve this task, e.g.~for fluid flow \cite{Garcke2015}.

Moreover, we briefly assess different solution strategies for the novel discretization based on the monolithic and partitioned solution paradigms comprising a feasible preconditioning strategy for the latter case. Efficient solution approaches for Cahn-Hilliard Navier-Stokes models with a focus on preconditioned iterative solvers are investigated in \cite{Aland2014, Bosch2018} for the monolithic case and in \cite{Kay2007, Ciloglu2025} for the partitioned case.

The outline of this work is as follows. In Section~\ref{sec:model} we start by introducing the governing two-phase model \revision{providing its sharp- and diffuse-interface formulations, along with the asymptotic analysis of the sharp-interface limit}.
We prove its thermodynamic consistency in Section~\ref{sec:thermo}, followed by a derivation of the (fully) discretized system in Section~\ref{sec:disc}. As the main theoretical result, we provide with Theorem \ref{theo:main} a discrete free energy dissipation inequality yielding the thermodynamic consistency of the proposed discrete system. \revision{In Section~\ref{sec:comput} we introduce two computational features for the non-reactive model.} Firstly, in Section~\ref{sec:solution} we establish numerical solution strategies based on the monolithic and partitioned paradigms, along with some notes on preconditioning techniques for the latter case. 
Secondly, we present a flexible preprocessing strategy for generating initial conditions for the phase-field variable taking parametrized sharp-interface descriptions as input in Section~\ref{sec:preprocessing}.
Lastly, Section~\ref{sec:numer} is devoted to numerical experiments that assess the modeling capabilities, thermodynamic consistency, and numerical solution strategies introduced before. 
We end with concluding remarks in Section~\ref{sec:concl}.

\section{Formulation of a model for two-phase fluid-solid dynamics} \label{sec:model}

\revision{
In this section we introduce the governing model for the interfacial dynamics of a fluid and a solid phase. We start by presenting the sharp-interface model in Section~\ref{sec:sharp_model} followed by a comprehensive introduction of its diffuse counterpart in Section~\ref{sec:diffuse_model}. Lastly, we provide the sharp-interface limit in Section~\ref{sec:limit}.

\subsection{Sharp-interface formulation} \label{sec:sharp_model}

Let $\Omega \subset \setR^\nu$, $\nu \in \{2,3\}$, be a bounded domain with Lipschitz boundary $\partial \Omega$ and $[0,T]$ a time interval with final time $0 < T < \infty$. In the context of two-phase fluid-solid dynamics, we assume the domain $\Omega$ to be decomposable such that $\Omega = \Omega_f(t) \cup \Omega_s(t) \cup \Gamma(t)$ (with $\Omega_f(t) \cap \Omega_s(t) = \emptyset$) for all $t \in [0,T]$. Here, $\Omega_f(t)$ and $\Omega_s(t)$ are bulk domains occupied by the fluid and the solid phase respectively and $\Gamma(t)$ denotes the fluid-solid interface. The domains (and hence also the interface) are time-dependent as the fluid domain can change by convection and precipitation/dissolution processes alter the solid domain. Moreover, we suppose the normal unit vector $\vecn \in \setR^\nu$ of $\Gamma(t)$ to point into $\Omega_s(t)$. The scalar $\nu_\Gamma$ denotes the normal velocity of the interface $\Gamma(t)$.  

\subsubsection{Bulk equations} \label{sec:bulk}
In $\Omega_f(t)$, $t \in (0,T]$, we consider the incompressible flow of a viscous fluid phase. Introducing as unknowns the  pressure $p \colon (\vecx,t) \mapsto \setR$ and the velocity $\mathbf{v} \colon (\vecx,t) \mapsto \setR^\nu$, the dynamics is then modeled by the incompressible Navier-Stokes equations
\begin{align}
 \nabla \cdot \mathbf{v} &= 0, \label{eq:incomp_con} \\
 \partial_t(\rho \mathbf{v}) + \nabla \cdot (\rho \mathbf{v} \otimes \mathbf{v}) + \nabla p &= \nabla \cdot (2 \gamma \nabla^s \mathbf{v}), \label{eq:momentum}
\end{align}
where the positive constants $\rho$ and $\gamma$ denote the fluid density and viscosity respectively, and $\nabla^s \mathbf{v} = \frac{1}{2} (\nabla \mathbf{v} + (\nabla \mathbf{v})^T)$ is the symmetric Jacobian.

Moreover, in order to model precipitation/dissolution effects, we assume the presence of ions that can dissolve in the fluid phase. Introducing as unknown the ion concentration $c \colon (\vecx,t) \mapsto \setR_0^+$, the dynamics is then governed by the transport-diffusion equation
\begin{align}
 \partial_t c + \nabla \cdot (\mathbf{v} c) = D \Delta c
 \label{eq:ion}
\end{align}
in $\Omega_f(t)$, $t \in (0,T]$, with positive diffusion rate constant $D$ and equilibrium concentration $c_{\text{eq}}$.

We assume the solid phase to be immobile, i.e. $\mathbf{v} = \mathbf{0}$ in $\Omega_s(t)$, $t \in (0,T]$, and to have a constant ion concentration $c^\ast > 0$. 

\subsubsection{Interface conditions} \label{sec:interface}
We consider the interface conditions between the bulk domains $\Omega_f(t)$ and $\Omega_s(t)$. We assume the velocity $\mathbf{v}$ to be continuous across both domains such that 
\begin{align}
 \mathbf{v} = \mathbf{0} \text{ on } \Gamma(t).
\label{eq:int_cond_v}
\end{align}
The velocity of the interface $\nu_\Gamma$ is determined by the precipitation and dissolution rates.  We introduce the rate function $r(c)$ that is assumed to model both, precipitation and dissolution, and to increase monotonically in the ion concentration $c$. We then choose
\begin{align}
 \nu_\Gamma = - r(c) + \alpha \sigma \kappa,
\label{eq:int_cond_nu}
\end{align}
where the last term models curvature effects on the reaction with geometric scaling constant $\alpha \geq 0$. $\sigma$ is the constant surface tension coefficient and $\kappa$ is the sum of the principal curvatures ($\nu - 1$ times the mean curvature).
Since for the solid domain $\Omega_s(t)$ we assume a constant concentration $c^\ast$, we impose the Rankine-Hugoniot-like condition
\begin{align}
 D \nabla c \cdot \mathbf{n} = \nu_\Gamma (c^\ast - c) \text{ on } \Gamma(t)
\label{eq:int_cond_c}
\end{align}
to ensure the conservation of ions.
}

\revisionT{\subsection{The phase-field model} \label{sec:diffuse_model}}

In this section we  introduce the governing diffuse-interface model for the interfacial dynamics of a fluid and a solid phase. It consists of a coupled Cahn-Hilliard Navier-Stokes system and accounts for the  hydromechanics of the fluid phase coupled to the phase-field evolution of the diffuse fluid-solid interface \revision{and in addition to the reactive transport of ions and the associated precipitation/dissolution  processes at the fluid-solid interface}.   
It can be understood as an adaption of the ternary Cahn-Hilliard Navier-Stokes model from 
\cite{Rohde2021} that describes the motion of two fluid phases. 

\revisionT{\subsubsection{Diffuse-interface formulation}}

As the phase-field variable $\phi$ we consider the volume fraction of the fluid phase, that is, $\phi = 1$ in the fluid phase and $\phi = 0$ in the solid phase. In a thin layer around the fluid-solid interface, called the diffuse transition zone, $\phi$ runs smoothly between 0 and 1. 
The two-phase model reads
\revisionT{
\begin{equation}
\label{eq:diff_model}
\begin{array}{rcl}
 \nabla \cdot (\phi \mathbf{v}) &=& 0,
 \\
 \partial_t (\rho_f \mathbf{v}) + \nabla \cdot ((\rho_f \mathbf{v} + \mathbf{J}_f) \otimes \mathbf{v}) &=& - \phi \nabla p + \nabla \cdot (2 \gamma \nabla^s \mathbf{v}) - \rho \frac{1}{\varepsilon} d(\phi)\mathbf{v} + \frac{1}{2} \rho \mathbf{v} R,
 \\
 \partial_t (\phi_c (c - c_{\text{eq}})) + \nabla \cdot ((\phi_c \mathbf{v} + \mathbf{J}_c) (c - c_{\text{eq}})) &=& D \nabla \cdot (\phi_c \nabla c) + R_c,
 \\
 \partial_t \phi + \nabla \cdot (\phi \mathbf{v} + \mathbf{J}) &=& R,
 \\
 \mu &=& \frac{W'_{\text{dw}}(\phi)}{\varepsilon} -  \varepsilon \Delta \phi
\end{array}
\end{equation}
in $\Omega \times (0,T]$.}
In \eqref{eq:diff_model}, the unknowns are the  pressure $p \colon (\vecx,t) \mapsto \setR$, the velocity $\mathbf{v} \colon (\vecx,t) \mapsto \setR^\nu$, the concentration $c \colon (\vecx,t) \mapsto [0,c^\ast]$, the phase-field variable $\phi \colon (\vecx,t) \mapsto  [0,1]$ and the chemical potential $\mu \colon (\vecx,t) \mapsto \setR$.\\ 
We denote  the density of the fluid phase by  $\rho$.
Furthermore, 
\begin{align}
 \phi_c = \phi \quad \textrm{and} \quad \rho_f = \rho \phi.
\label{eq:defrhof}
\end{align}
We introduce 
\begin{align}
 \mathbf{J} = - M \varepsilon \nabla \mu, \quad \mathbf{J}_f = \rho \mathbf{J} \quad \textrm{and} \quad \mathbf{J}_c = \mathbf{J}
\label{eq:defJf}
\end{align}
as the Cahn-Hilliard fluxes with constant state-independent mobility $M$. Moreover, $R$ and $R_c$ denote \revision{reaction terms that model ion precipitation and dissolution localized to the fluid-solid interface} satisfying
\begin{align}
 R_c = (c^\ast - c_{\text{eq}}) R.
\label{eq:Rc}
\end{align}
The explicit form of $R$ is discussed in equation~\eqref{eq:react}. 

As is usual for phase-field models we introduce a parameter $\varepsilon > 0$ such that the width of the diffuse transition zone is of order $\varepsilon$. The function $W_\text{dw}(\phi)$ is a double-well potential taking minimal values in the pure phases ($\phi = 0$ and $\phi = 1$). 

In the momentum equations~$\eqref{eq:diff_model}_2$, the constant  $\gamma>0$ denotes the state-independent viscosity and $d(\phi)$ accounts for velocity dissipation.
\\
We augment \eqref{eq:diff_model} with initial and boundary conditions given by 
\begin{equation}
\label{icbc}
\begin{array}{c}
\vecv(\cdot, 0) = \vecv_0, c(\cdot, 0) = c_0, \phi(\cdot, 0) = \phi_0 \text{ in } \Omega,\qquad \vecv = \mathbf{0}, \nabla c \cdot \vecn_\Omega = 0, \nabla\phi \cdot \vecn_\Omega = 0, \nabla 
\mu \cdot \vecn_\Omega = 0 \text{ on } \partial \Omega.\\
\end{array}
\end{equation}
Here, $\vecv_0: \Omega \to \setR^\nu$, $c_0: \Omega \to [0,c^\ast]$ and $\phi_0: \Omega \to [0,1]$ denote the given initial velocity, concentration and phase field, and $\vecn_\Omega \in \setR^\nu$ denotes the outer normal unit vector on $\partial \Omega$.\\
In the following remark we discuss the 
key modeling aspects related to \eqref{eq:diff_model}.
\begin{remark}[Modeling aspects] \label{rem: model}
\hfill
\begin{itemize}
\item[(i)] We follow the seminal modeling approach  for the  incompressible flow of two fluid phases as introduced in \cite{Abels2012}. 
In particular, for the     incompressibility equation~$\eqref{eq:diff_model}_1$, we follow the idea of volume averaging presented in \cite{Abels2012} which in the given setting leads to $\phi \mathbf{v}$ being divergence free, replacing the usual incompressibility constraint on $\mathbf{v}$ alone. \revision{Note that here $\mathbf{v}$ represents the velocity of the fluid not the mixture.

This is in contrast to the classical approach of mass averaging by \cite{Lowengrub1998}. We note, however, that in \cite{tenEikelder2023} a unifying framework for the two modeling approaches and for Cahn-Hilliard Navier-Stokes models with non-matching densities in general is presented.} 

\item[(ii)] As smooth phase-field parameter $\phi$ we consider the volume fraction of the fluid phase. Hence, it would be desirable to ensure that for all $\mathbf{x} \in \Omega$ and $t \in [0,T]$ the range restriction $\phi(\mathbf{x},t) \in [0,1]$ is satisfied. However, the Cahn-Hilliard evolution that models the phase-field dynamics is known to not fulfill a priori such a maximum principle. As a means of remediation, we employ an unbounded potential function of the general form
\begin{align}
 W_\text{dw}(\phi) = \revision{450 \phi^4 (1 - \phi)^4} + \ell(\phi) + \ell(1 - \phi)
\label{eq:defDW}
\end{align}
with the \revision{double-well term of polynomial degree eight (scaled by a constant factor for consistency with the sharp-interface limit)} and a limiter function $\ell(\phi)$ that diverges at $\phi = -\delta$ for some small $\delta > 0$. Thus, $W_\text{dw}(\phi)$ diverges at $\phi = -\delta$ and $\phi = 1 + \delta$, and by establishing an energy estimate \revisionT{in Section~\ref{sec:thermo}} the relaxed constraint $\phi(\mathbf{x},t) \in (-\delta, 1 + \delta)$ is satisfied.

Note that in the limit $\delta \rightarrow 0$ the double-well function $W_\text{dw}(\phi)$ converges point-wise to a potential of double-obstacle type. There exists extensive analysis for the Cahn-Hilliard equation with double-obstacle potential, cf.~e.g.~\cite{Blowey1991}. For numerical reasons, however, we refrain from following this ansatz as it results in a model comprising variational inequalities.

\item[(iii)] The solid phase is assumed to be immobile. By \revision{splitting} the pressure-like term in the momentum equations~$\eqref{eq:diff_model}_2$ \revision{in a conservative and a non-conservative contribution}
as 
\begin{align*}
 - \phi \nabla p = - \nabla (\phi p) + p \nabla \phi
\end{align*}
it becomes apparent that the latter can be a source (or sink) for momentum. \revision{In fact, it acts as a normal force between the two phases as $\nabla \phi$ is perpendicular to the fluid-solid interface.}

In order to achieve thermodynamic consistency below, we follow the approach of adding a flux term in the momentum equations, cf.~\cite{Abels2012}, and generalize it to the case of a solid phase present in the model, yielding the contribution $\nabla \cdot \left(\mathbf{J}_f \otimes \mathbf{v}\right)$.

\revision{Concerning the drag term $\rho \frac{1}{\varepsilon} d(\phi) \mathbf{v}$, the function $d(\cdot)$ can be any smooth, decreasing function satisfying $d(1) = 0$, $d'(1) = 0$ and $d(0) = d_0$ with positive constant $d_0$ independent of $\varepsilon$. This term ensures that in the sharp-interface limit we recover $\mathbf{v} = 0$ for the solid phase, along with a no-slip condition between the fluid and the solid phase. Note that the choice $d(\phi) = d_0 \big(1 - \phi\big)^2$ from \cite{Rohde2021} fulfills the criteria.}

\item[(iv)] The equation for the ion concentration $c$ $\eqref{eq:diff_model}_{3}$ is of convection-diffusion-reaction type. The transport due to the Cahn-Hilliard equation is taken into account through the term $\nabla \cdot(\mathbf{J}_c c)$. \revision{Following \cite{Rohde2021}, for the reaction kinetics we choose
\begin{align}
 r(c) = g(c) + g'(c)(c^\ast - c), \qquad R = - \frac{q(\phi)}{\varepsilon}\left(r(c) + \tilde \alpha \sigma \mu\right)
\label{eq:react}
\end{align}
with convex function $g(c)$ to obtain thermodynamic consistency. The reaction term $R$ depends on the reaction rate $r(c)$ and the chemical potential $\mu$, where the latter's contribution accounts for curvature effects on the reaction, see also \eqref{eq:int_cond_nu}. We define the modified geometric scaling constant
\begin{align}
 \tilde \alpha 
 = 
 \begin{cases}
  \alpha \quad &\text{for } \alpha > 0 \\
  \revisionT{\delta_\alpha} \quad &\text{for } \alpha = 0
 \end{cases}
\label{eq:reg_alpha}
\end{align}
\revisionT{with some small parameter $\delta_\alpha > 0$ yet to be chosen. This definition prevents non-physical behavior in the special case $\alpha = 0$.} Moreover, we set
\begin{align}
 q(\phi) = 
 \begin{cases} 
  \sqrt{2 W_{\text{dw}}(\phi)} = 30 \phi^2 (1 - \phi)^2 \quad &\text{for } 0 \leq \phi \leq 1 \\
  0 \quad &\text{else}
 \end{cases}
\label{eq:defQ}
\end{align}
and hence the reaction is localized to the fluid-solid interface. Note that for the sharp-interface limit we require $q'(1) = 0$ rendering the choice of a quartic double-well potential infeasible.
}
\item[(v)] The classical formulation of  the Cahn-Hilliard equation  
leads to  fourth order  operators. To enable a proper notion of a weak solution we  split  the equation of fourth order in space into the  two equations $\eqref{eq:diff_model}_{4,5}$ of second order introducing the chemical potential $\mu$. The coupling to the Navier-Stokes equations~$\eqref{eq:diff_model}_{1,2}$ is due to the advection of the fluid phase only. Note again that the solid phase is not advected. 
\item[(vi)] The boundary conditions in \eqref{icbc}  provide a simple setting. For the theoretical part we make no attempt to account for inflow/outflow conditions or more advanced concepts like the enforcement of dynamic contact angles.   
 \end{itemize}
\end{remark}

\revisionT{
\subsubsection{Regularized diffuse-interface formulation} \label{sec:reg_diffuse_model}

In order to achieve thermodynamical consistency we establish the following regularization of model \eqref{eq:diff_model}. We need to ensure that $\phi$, $\phi_c$ and $\rho_f$ do not attain negative values as this would lead to a degeneration of the model. Therefore, following Remark~\ref{rem: model}(ii) we redefine these quantities using the parameter $\delta$ from the definition of the double-well potential. We define the regularized fluid fraction $\tilde \phi_f$ by 
\begin{equation}\label{eq:defphif}
 \tilde \phi_f = 2 \delta + (1-2\delta) \phi.
\end{equation}
Moreover, we have
\begin{align}
 \tilde \phi_c = \phi + \delta \quad \textrm{and} \quad \tilde \rho_f = \rho_f + \rho \delta.
\label{eq:defrhof_2}
\end{align}
We directly obtain that for $\phi > -\delta$ the quantities $\tilde \phi_f$, $\tilde \phi_c$ and $\tilde \rho_f$ are positive. 
We formulate the regularized model by
\begin{equation}
\label{eq:model}
\begin{array}{rcl}
 \nabla \cdot (\tilde \phi_f \mathbf{v}) &=& 0, 
 \\
 \partial_t (\tilde \rho_f \mathbf{v}) + \nabla \cdot ((\rho_f \mathbf{v} + \mathbf{J}_f) \otimes \mathbf{v}) &=& -\tilde \phi_f \nabla p + \nabla \cdot (2 \gamma \nabla^s \mathbf{v}) - \rho \frac{1}{\varepsilon} d(\tilde \phi_f)\mathbf{v} + \mathbf{\tilde S} + \frac{1}{2} \rho \mathbf{v} R, 
 \\
 \partial_t (\tilde \phi_c (c - c_{\text{eq}})) + \nabla \cdot ((\phi_c \mathbf{v} + \mathbf{J}_c) (c - c_{\text{eq}})) &=& D \nabla \cdot (\tilde \phi_c \nabla c) + R_c,
 \\
 \partial_t \phi + \nabla \cdot (\phi \mathbf{v} + \mathbf{J}) &=& R, 
 \\
 \mu &=& \frac{W'_{\text{dw}}(\phi)}{\varepsilon} -  \varepsilon \Delta \phi  
\end{array}
\end{equation}
in $\Omega \times (0,T]$.
To keep consistency with the regularized $\tilde \phi_f$ we obtain the (non-physical) effective surface tension term $\mathbf{\tilde S}$ in the momentum equations. It is given as 
\begin{align}
 \mathbf{\tilde S} = 2\delta \sigma (1 - \phi) \nabla \mu
\label{eq:defS}
\end{align}
with the (constant) surface tension coefficient $\sigma > 0$ between the two phases. This term is indispensable for the (proof of) thermodynamic consistency, but has little influence on the global solution as it is scaled with $\delta$. In addition, 
in \eqref{eq:reg_alpha} we employ the scaling $\delta_\alpha = \delta \alpha_{\textrm{ref}}$ with a characteristic scale parameter $\alpha_{\textrm{ref}}$ independent of $\delta$. For the remainder of this work we set $\alpha_{\textrm{ref}} = 1$.
}

\revision{
For classical solutions to model~\eqref{eq:model} satisfying the boundary conditions from \eqref{icbc} using the divergence theorem we obtain conservative properties for the phase-field variable $\phi$ and the ion concentration $c$. Firstly, following from $\eqref{eq:model}_4$ the phase-field variable $\phi$ is balanced by the reaction term $R$ as 
\begin{align*}
 \frac{d}{dt} \int_\Omega \phi \;d\vecx = \int_\Omega R \;d\vecx.
\end{align*}
Secondly, the total amount of ions is given by
\begin{align*}
 \tilde \phi_c c + (1 - \tilde \phi_c) c^\ast = \tilde \phi_c (c - c_{\text{eq}}) - \tilde{\phi_c} (c^\ast - c_{\text{eq}}) + c^\ast
\end{align*}
and hence we have
\begin{align*}
 \partial_t (\tilde \phi_c c + (1 - \tilde \phi_c) c^\ast) = \partial_t (\tilde \phi_c (c - c_{\text{eq}})) - \partial_t \phi (c^\ast - c_{\text{eq}})
\end{align*}
employing the definition of $\tilde \phi_c$ \eqref{eq:defrhof_2}.
Therefore, following from $\eqref{eq:model}_3$ and $\eqref{eq:model}_4$ along with \eqref{eq:Rc} we get
\begin{align*}
 \frac{d}{dt} \int_\Omega \tilde \phi_c c + (1 - \tilde \phi_c) c^\ast \;d\vecx = \int_\Omega R_c - R (c^\ast - c_{\text{eq}}) \;d\vecx = 0,
\end{align*}
i.e., the conservation of ions.
}

\revision{
\begin{remark}[Non-reactive flow model]
We can reduce model~\eqref{eq:model} to a non-reactive flow model by setting $R = 0$ for the reaction term and hence not solving for the ion concentration $c$. The resulting model reads
\begin{equation}
\label{eq:model_non_react}
\begin{array}{rcl}
 \nabla \cdot (\tilde \phi_f \mathbf{v}) &=& 0, 
 \\
 \partial_t (\tilde \rho_f \mathbf{v}) + \nabla \cdot ((\rho_f \mathbf{v} + \mathbf{J}_f) \otimes \mathbf{v}) &=& -\tilde \phi_f \nabla p + \nabla \cdot (2 \gamma \nabla^s \mathbf{v}) - \rho \frac{1}{\varepsilon} d(\tilde \phi_f) \mathbf{v} + \mathbf{\tilde S}, 
 \\ 
 \partial_t \phi + \nabla \cdot (\phi \mathbf{v} + \mathbf{J}) &=& {0}, 
 \\
 \mu &=&\ds  \frac{W'_{\text{dw}}(\phi)}{\varepsilon} -  \varepsilon \Delta \phi  
\end{array}  \qquad \text{in }  \Omega \times (0,T]
\end{equation}
and describes a viscous fluid phase interacting with an immobile solid phase, which may alter subject to a geometric evolution. The theoretical considerations in the following sections also hold true for the non-reactive model.
We note that we do not (need to) use the scaled 8th order potential term $450 \phi^4 (1 - \phi)^4$ here as the classical quartic double-well potential $\phi^2 (1 - \phi)^2$ is sufficient to comply with the corresponding sharp-interface limit, cf.~Remark~\ref{rem: model}(ii) and (iv), and Section~\ref{sec:limit} below.
\end{remark}
}

\revision{
\subsection{Sharp-interface limit} \label{sec:limit}
We use the method of matched asymptotic expansions, see e.g.~\cite{Caginalp1988}, to show that the sharp-interface formulation presented in Section~\ref{sec:sharp_model} is the limit of the \revisionT{regularized} phase-field model presented in Section~\ref{sec:reg_diffuse_model} for $\varepsilon \rightarrow 0$.

We employ the scaling $\delta = \frac{\varepsilon}{L_{\textrm{ref}}}$ with a characteristic length scale $L_{\textrm{ref}}$ independent of $\varepsilon$. Hence, the $\delta$-regularizations vanish for $\varepsilon \rightarrow 0$. For the sake of simplicity we set $L_{\textrm{ref}} = 1$. Moreover, we assume that there are classical solutions to the diffuse-interface model \eqref{eq:model}. Together with the energy estimate from Theorem~\ref{theo:analysis}, we then obtain the bound $\phi \in (-\varepsilon, 1 + \varepsilon)$. Furthermore, we assume the evolution to be on a $\mathcal{O}(1)$-timescale.

\subsubsection{Outer expansions}
We assume that solutions away from the interface can be written in terms of outer expansions of the unknowns. Exemplarily for $\phi$ we thus have
\begin{align*}
 \phi^{\text{out}}(\mathbf{x}, t) = \phi_0^{\text{out}}(\mathbf{x}, t) + \varepsilon \phi_1^{\text{out}}(\mathbf{x}, t) + \varepsilon^2 \phi_2^{\text{out}}(\mathbf{x}, t) + \ldots
\end{align*}
with $\phi_i^{\text{out}}$, $i \in \mathbb{N}_0$, independent of $\varepsilon$. This notation is also applied for non-primary variables, such as
\begin{align*}
 \tilde \phi_f^{\text{out}}(\mathbf{x}, t) = \tilde \phi_{f,0}^{\text{out}}(\mathbf{x}, t) + \varepsilon \tilde \phi_{f,1}^{\text{out}}(\mathbf{x}, t) + \ldots = \phi_0^{\text{out}}(\mathbf{x}, t) + \varepsilon (2 - 2 \phi_0^{\text{out}}(\mathbf{x}, t) + \phi_1^{\text{out}}(\mathbf{x}, t)) + \ldots\,.
\end{align*}
Nonlinearities are expanded through Taylor series, i.e., for a generic function $h$ and a variable $u = u_0 + \varepsilon u_1 + \ldots$ we have
\begin{align*}
 h(u) = h(u_0) + \varepsilon h'(u_0) u_1 + \mathcal{O}(\varepsilon^2)
\end{align*}
given that the respective derivatives exist. The terms are then grouped by powers of $\varepsilon$.  

In order to recover the equations in the bulk domains, see Section~\ref{sec:bulk}, we consider the outer expansions of the unknowns $p$, $\mathbf{v}$, $\phi$, $\mu$ and $c$.

\textit{Expansion of $\eqref{eq:model}_5$, $\mathcal{O}(\varepsilon^{-1})$:}
We obtain $W'_{\text{dw}}(\phi_0^{\text{out}}) = 0$ which for our choice of the double-well potential has the solutions $0, \frac{1}{2}, 1$. Following \cite{vanNoorden2011}, $\phi_0^{\text{out}} = 0$ and $\phi_0^{\text{out}} = 1$ are the stable solutions, whereas $\phi_0^{\text{out}} = \frac{1}{2}$ is unstable due to $W''_{\text{dw}}(\phi_0^{\text{out}} = \frac{1}{2}) < 0$. Now, let $\Omega_f(t)$ and $\Omega_s(t)$ be the corresponding subdomains of $\Omega$ where $\phi_0^{\text{out}}$ is 1 and 0 respectively.

\textit{Expansion of $\eqref{eq:model}_1$, $\mathcal{O}(1)$:} 
For $\phi_0^{\text{out}} = 1$ we have $\tilde \phi_f^{\text{out}} = 1 + \mathcal{O}(\varepsilon)$ and hence get 
\begin{align*} 
 \nabla \cdot \mathbf{v}_0^{\text{out}} = 0,
\end{align*} 
i.e., equation~\eqref{eq:incomp_con}.

\textit{Expansion of $\eqref{eq:model}_3$, $\mathcal{O}(1)$:}
For $\phi_0^{\text{out}} = 1$ we have $\tilde \phi_c^{\text{out}} = 1 + \mathcal{O}(\varepsilon)$. Moreover, following from $q(\phi_0^{\text{out}}) = 0$ and $q'(\phi_0^{\text{out}}) = 0$, we find $q(\phi^{\text{out}}) = \mathcal{O}(\varepsilon^2)$. Therefore, we obtain $R_c = - (c^\ast - c_{\text{eq}}) \frac{q(\phi^{\text{out}})}{\varepsilon} (r(c^{\text{out}}) + \tilde \alpha \sigma \mu^{\text{out}}) = \mathcal{O}(\varepsilon)$ and using that $\nabla \cdot \mathbf{v}_0^{\text{out}} = 0$ we hence get
\begin{align*}
 \partial_t c_0^{\text{out}} + \nabla \cdot (\mathbf{v}_0^{\text{out}} c_0^{\text{out}}) = D \Delta c_0^{\text{out}},
\end{align*}
i.e., equation~\eqref{eq:ion}.

\textit{Expansion of $\eqref{eq:model}_2$, $\mathcal{O}(1)$:}
For $\phi_0^{\text{out}} = 1$ we have $\tilde \phi_f^{\text{out}} = 1 + \mathcal{O}(\varepsilon)$. Following from $d(\phi_0^{\text{out}}) = 0$ and $d'(\phi_0^{\text{out}}) = 0$, we find $d(\phi^{\text{out}}) = \mathcal{O}(\varepsilon^2)$ and thus $\rho \frac{1}{\varepsilon} d(\tilde \phi_f^{\text{out}}) \mathbf{v}^{\text{out}} = \mathcal{O}(\varepsilon)$. From their respective definitions we moreover have $\mathbf{\tilde S} = \mathcal{O}(\varepsilon)$ and $R_f = \mathcal{O}(\varepsilon)$ (analogous to $R_c$ above). We hence get
\begin{align*}
 \partial_t(\rho \mathbf{v}_0^{\text{out}}) + \nabla \cdot (\rho \mathbf{v}_0^{\text{out}} \otimes \mathbf{v}_0^{\text{out}}) + \nabla p_0^{\text{out}} = \nabla \cdot (2 \gamma \nabla^s \mathbf{v}_0^{\text{out}}),
\end{align*}
i.e., equation~\eqref{eq:momentum}.

\textit{Expansion of $\eqref{eq:model}_2$, $\mathcal{O}(\varepsilon^{-1})$:}
For $\phi_0^{\text{out}} = 0$ we have $\tilde \phi_f^{\text{out}} = \mathcal{O}(\varepsilon)$. As $d(0) = d_0$ we find $d(\tilde \phi_f^{\text{out}}) = \mathcal{O}(1)$. Thus, we get $\rho d(\tilde \phi_f^{\text{out}}) \mathbf{v}_0^{\text{out}} = \mathbf{0}$ which is equivalent to $\mathbf{v}_0^{\text{out}} = \mathbf{0}$, i.e., the immobility of the solid phase.

\subsubsection{Inner expansions}
For the interface we set $\Gamma(t) = \left\{\mathbf{x} \in \Omega|\, \phi(\mathbf{x},t) = \frac{1}{2}\right\}$. Let $\mathbf{s} \in \mathbb{R}^{\nu - 1}$ denote a local parametrization along $\Gamma(t)$ such that $\mathbf{x}(\mathbf{s}, t) \in \Gamma(t)$, and $\mathbf{n} \in \mathbb{R}^\nu$ the unit normal vector at $\Gamma(t)$ directed towards the solid phase ($\phi = 0$). Then, for the normal velocity of the interface we can write
\begin{align*}
 \nu_\Gamma(\mathbf{s}, t) = \partial_t \mathbf{x}(\mathbf{s}, t) \cdot \mathbf{n}(\mathbf{s}, t).
\end{align*}
Introducing 
\begin{align*}
 \mathbf{x}(\mathbf{s}, t, \zeta) = \mathbf{x}(\mathbf{s}, t) + \zeta \mathbf{n}(\mathbf{s}, t),
\end{align*}
where $\zeta$ is the signed distance from the interface to the point $\mathbf{x}(\mathbf{s}, t, \zeta)$, we define local curvilinear coordinates $(\mathbf{s}, \zeta)$ near the interface.

Given that the width of the diffuse interface is in the order of $\varepsilon$, we introduce $z = \frac{\zeta}{\varepsilon}$ as a rescaled signed distance to $\Gamma(t)$. 
For a generic scalar variable $u$ and a generic vectorial variable $\mathbf{U}$, the transformation rules
\begin{align*}
 \partial_t u &= - \frac{1}{\varepsilon} \nu_\Gamma \partial_z u + \mathcal{O}(1), \\
 \nabla u &= \frac{1}{\varepsilon} \partial_z u \mathbf{n} + \nabla_{\Gamma^{\text{out}}} u  + \mathcal{O}(\varepsilon), \\
 \nabla \cdot \mathbf{U} &= \frac{1}{\varepsilon} \partial_z \mathbf{U} \cdot \mathbf{n} + \nabla_{\Gamma^{\text{out}}} \cdot \mathbf{U} + \mathcal{O}(\varepsilon), \\
 \nabla \mathbf{U} &= \frac{1}{\varepsilon} \partial_z \mathbf{U} \otimes \mathbf{n} + \nabla_{\Gamma^{\text{out}}} \mathbf{U} + \mathcal{O}(\varepsilon), \\
 \Delta u &= \frac{1}{\varepsilon^2} \partial_{zz} u - \frac{1}{\varepsilon} \kappa \partial_z u + \mathcal{O}(1)
\end{align*}
hold, where $\kappa$ denotes the sum of principal curvatures of $\Gamma^{\text{out}}$ and $\nabla_{\Gamma^{\text{out}}}$ the surface gradient of $\Gamma^{\text{out}}$.

We assume that solutions close to the interface can be written in terms of inner expansions of the unknowns. Exemplarily for $\phi$ we thus have
\begin{align*}
 \phi^{\text{in}}(\mathbf{s}, t, z) = \phi_0^{\text{in}}(\mathbf{s}, t, z) + \varepsilon \phi_1^{\text{in}}(\mathbf{s}, t, z) + \varepsilon^2 \phi_2^{\text{in}}(\mathbf{s}, t, z) + \ldots
\end{align*}
with the local, rescaled coordinates $\mathbf{s}$ and $z$.

For a generic variable $u$ with outer expansion $u^{\text{out}}$ and inner expansion $u^{\text{in}}$, we have the identity $u^{\text{out}}(\mathbf{x}(\mathbf{s}, t, \varepsilon z), t) = u^{\text{in}}(\mathbf{s}, t, z)$. For fixed $t$ and $\mathbf{s}$ we introduce the one-sided limits 
\begin{align*}  
 \lim_{\zeta \rightarrow 0^{\pm}} \mathbf{x}(\mathbf{s}, t, \zeta) = \mathbf{x}_{\pm}
\end{align*} for the outer expansions. Following \cite{Caginalp1988}, we get the following matching conditions in the limit $z \rightarrow \pm \infty$:
\begin{align}
 u_0^{\text{in}}(\mathbf{s}, t, \pm \infty) &= u_0^{\text{out}}(\mathbf{x}_{\pm}, t), \label{eq:match_cond_1} \\
 \partial_z u_0^{\text{in}}(\mathbf{s}, t, \pm \infty) &= 0, \label{eq:match_cond_2} \\
 \partial_z u_1^{\text{in}}(\mathbf{s}, t, \pm \infty) &= \nabla u_0^{\text{out}}(\mathbf{x}_\pm, t) \cdot \mathbf{n}. \label{eq:match_cond_3}
\end{align}

In order to recover the interface conditions, see Section~\ref{sec:interface}, we consider the inner expansions of the unknowns $p$, $\mathbf{v}$, $\phi$, $\mu$ and $c$. For the subsequent analysis we need to investigate leading order and first order terms of these expansions.

\paragraph{Leading order}
We start by investigating the leading order expansions.

\textit{Expansion of $\eqref{eq:model}_5$, $\mathcal{O}(\varepsilon^{-1})$}:
We obtain 
\begin{align}
 0 = \partial_{zz} \phi_0^{\text{in}} - W'_{\text{dw}}(\phi_0^{\text{in}})
\label{eq:ODE_phi}
\end{align}
and following from \eqref{eq:match_cond_1} we have $\phi_0^{\text{in}}(-\infty) = 1$ and $\phi_0^{\text{in}}(\infty) = 0$. Moreover, the definition of $\Gamma(t)$ yields $\phi_0^{\text{in}}(0) = \frac{1}{2}$. The solution to \eqref{eq:ODE_phi} is then implicitly given by
\begin{align*}
 z = \frac{1}{30} \left(\frac{1}{\phi_0^{\text{in}}(z)} - \frac{1}{1 - \phi_0^{\text{in}}(z)} + 2 \ln\!\left(\frac{1 - \phi_0^{\text{in}}(z)}{\phi_0^{\text{in}}(z)}\right)\right).
\end{align*}
We have 
\begin{align*}   
 &\int_{-\infty}^z \partial_{\xi} \phi_0^{\text{in}} \partial_{\xi\xi} \phi_0^{\text{in}} d\xi = \int_{-\infty}^z \partial_{\xi}(\frac{1}{2} (\partial_{\xi} \phi_0^{\text{in}})^2) d\xi = \left[\frac{1}{2} (\partial_{\xi} \phi_0^{\text{in}})^2\right]_{-\infty}^z = \frac{1}{2} (\partial_z \phi_0^{\text{in}})^2 \quad \text{and} \quad \\
 &\int_{-\infty}^z \partial_{\xi} \phi_0^{\text{in}} W'_{\text{dw}}(\phi_0^{\text{in}}) d\xi = \int_1^{\phi_0^{\text{in}}} W'_{\text{dw}}(\phi) d\phi = \left[W_{\text{dw}}(\phi)\right]_1^{\phi_0^{\text{in}}} = W_{\text{dw}}(\phi_0^{\text{in}}),
\end{align*}
such that by multiplying \eqref{eq:ODE_phi} with $\partial_z \phi_0^{\text{in}}$ and integrating, we obtain the equipartition of energy 
\begin{align}   
 \frac{1}{2} (\partial_z \phi_0^{\text{in}})^2 = W_{\text{dw}} (\phi_0^{\text{in}}).
\label{eq:equi}
\end{align}
As $\phi_0^{\text{in}}$ is monotonously decreasing, this yields $\partial_z \phi_0^{\text{in}} = - \sqrt{2 W_{\text{dw}} (\phi_0^{\text{in}})}$.

\textit{Expansion of $\eqref{eq:model}_1$, $\mathcal{O}(\varepsilon^{-1})$}:
We have $\tilde \phi_{f,0}^{\text{in}} = \phi_0^{\text{in}} + \mathcal{O}(\varepsilon)$ and hence obtain 
\begin{align*}
 \partial_z (\phi_0^{\text{in}} \mathbf{v}_0^{\text{in}}) \cdot \mathbf{n} = 0.
\end{align*} 
By integrating and using \eqref{eq:match_cond_1} along with $\phi_0^{\text{in}}(\infty) = 0$ (and $\mathbf{v}_0^{\text{in}}(\infty) = \mathbf{0}$), we get 
\begin{align}
 \mathbf{v}_0^{\text{in}}(z) \cdot \mathbf{n} = \mathbf{v}_0^{\text{out}}(\mathbf{x}_{_{-}},t) \cdot \mathbf{n} = 0 \quad \text{for all } z \in (-\infty, \infty).
\label{eq:v_normal}
\end{align}

\textit{Expansion of $\eqref{eq:model}_3$, $\mathcal{O}(\varepsilon^{-2})$}:
We have $\tilde \phi_c^{\text{in}} = \phi_0^{\text{in}} + \mathcal{O}(\varepsilon)$ and $q(\phi^{\text{in}}) = q(\phi_0^{\text{in}}) + \mathcal{O}(\varepsilon) = \mathcal{O}(1)$. Hence, for the reaction term we obtain $R_c = \mathcal{O}(\varepsilon^{-1})$, yielding in leading order
\begin{align*}
 \partial_z (\phi_0^{\text{in}} \partial_z c_0^{\text{in}}) = 0.
\end{align*}
From the matching conditions \eqref{eq:match_cond_1} and \eqref{eq:match_cond_2} at $z = -\infty$ we then get
\begin{align}
 c_0^{\text{in}}(z) = c_0^{\text{out}}(\mathbf{x}_{_{-}},t) \quad \text{for all } z \in (-\infty, \infty).
\label{eq:c_cont}
\end{align}

\textit{Expansion of $\eqref{eq:model}_2$, $\mathcal{O}(\varepsilon^{-2})$}:
For the reaction term we have $R_f = \mathcal{O}(\varepsilon^{-1})$ (analogous to $R_c$ above) and $\mathbf{\tilde S} = \mathcal{O}(1)$ for the effective surface tension term. In leading order we thus have
\begin{align*}
 0 = \partial_z (\gamma (\partial_z \mathbf{v}_0^{\text{in}} \otimes \mathbf{n} + \mathbf{n} \otimes \partial_z \mathbf{v}_0^{\text{in}})) \mathbf{n} = \partial_z (\gamma ((\mathbf{n} \cdot \mathbf{n}) \partial_z \mathbf{v}_0^{\text{in}} + (\partial_z \mathbf{v}_0^{\text{in}} \cdot \mathbf{n}) \mathbf{n})) = \partial_z (\gamma \partial_z \mathbf{v}_0^{\text{in}}),
\end{align*}
where we employ that $\partial_z \mathbf{v}_0^{\text{in}} \cdot \mathbf{n} = 0$ following from \eqref{eq:v_normal}. We integrate and use \eqref{eq:match_cond_2} along with $\gamma > 0$ to obtain 
\begin{align*}
 0 = \partial_z \mathbf{v}_0^{\text{in}}.
\end{align*}
Using \eqref{eq:match_cond_1} we find 
\begin{align*}
 \mathbf{v}_0^{\text{in}}(z) = \mathbf{v}_0^{\text{out}}(\mathbf{x}_{_{-}},t) = \mathbf{v}_0^{\text{out}}(\mathbf{x}_{_{+}},t) = 0 \quad \text{for all } z \in (-\infty, \infty), 
\end{align*}
i.e., interface condition~\eqref{eq:int_cond_v}.

\textit{Expansion of $\eqref{eq:model}_4$, $\mathcal{O}(\varepsilon^{-1})$}:
In leading order we have 
\begin{align*}
 - \nu_\Gamma \partial_z \phi_0^{\text{in}} + \partial_z (\phi_0^{\text{in}} \mathbf{v}_0^{\text{in}}) \cdot \mathbf{n} - M \partial_{zz} \mu_0^{\text{in}} = - q(\phi_0^{\text{in}}) (r(c_0^{\text{in}}) + \tilde \alpha \sigma \mu_0^{\text{in}}).
\end{align*}
We have $\mathbf{v}_0^{\text{in}} = \mathbf{0}$ and $\tilde \alpha = \alpha + \mathcal{O}(\varepsilon)$. The definition of $q$ \eqref{eq:defQ} and the equipartition of energy \eqref{eq:equi} give the identity $q(\phi_0^{\text{in}}) = - \partial_z \phi_0^{\text{in}}$. Hence, we get 
\begin{align*}
 \nu_\Gamma \partial_z \phi_0^{\text{in}} + M \partial_{zz} \mu_0^{\text{in}} = - \partial_z \phi_0^{\text{in}} (r(c_0^{\text{in}}) + \alpha \sigma \mu_0^{\text{in}}) 
\end{align*}
which can equivalently be written as an ODE for $\mu_0^{\text{in}}$ reading 
\begin{align} 
 (M \partial_{zz} + \partial_z \phi_0^{\text{in}} \alpha \sigma) \mu_0^{\text{in}} = - \partial_z \phi_0^{\text{in}}(r(c_0^{\text{in}}) + \nu_\Gamma).
\label{eq:ODE_mu}
\end{align}
In addition, following from \eqref{eq:match_cond_2} we obtain the asymptotics $\partial_z \mu_0^{\text{in}}(-\infty) = 0 = \partial_z \mu_0^{\text{in}}(\infty)$.

In the case $\alpha = 0$ by integrating we find 
\begin{align*}   
 0 = r(c_0^{\text{in}}) + \nu_\Gamma,
\end{align*}
where $r(c_0^{\text{in}})$ is constant due to \eqref{eq:c_cont}. When this condition is satisfied any $\mu_0^{\text{in}} = \text{const.}$ is a solution to \eqref{eq:ODE_mu}.

For $\alpha > 0$ the homogeneous part has the solution $\mu_0^{\text{in}} = 0$. In total, we obtain the unique, constant solution 
\begin{align}
 \mu_0^{\text{in}} = - \frac{1}{\alpha \sigma}(\nu_\Gamma + r(c_0^{\text{in}}))
\label{eq:mu_const}
\end{align}
leading to
\begin{align}
 \nu_\Gamma = - \alpha \sigma \mu_0^{\text{in}} - r(c_0^{\text{in}}),
\label{eq:nu_Gamma}
\end{align}
which holds true also in the case $\alpha = 0$.

\paragraph{First order}
To complete the sharp-interface limit we consider first order expansions.

\textit{Expansion of $\eqref{eq:model}_3$, $\mathcal{O}(\varepsilon^{-1})$}:
Using $c^{\text{in}} = c_0^{\text{in}} + \varepsilon c_1^{\text{in}} + \mathcal{O}(\varepsilon^2)$ and $\tilde \alpha = \alpha + \mathcal{O}(\varepsilon)$ we obtain
\begin{align*}
 - \nu_\Gamma \partial_z (\phi_0^{\text{in}} c_0^{\text{in}}) + \partial_z (\phi_0^{\text{in}} c_0^{\text{in}} \mathbf{v}_0^{\text{in}}) \cdot \mathbf{n} - \partial_z (c_0^{\text{in}} M \partial_z \mu_0^{\text{in}}) = D \partial_z (\phi_0^{\text{in}} \partial_z c_1^{\text{in}}) - c^\ast q(\phi_0^{\text{in}}) (r(c_0^{\text{in}}) + \alpha \sigma \mu_0^{\text{in}}).
\end{align*}
As before we employ the identities $\mathbf{v}_0^{\text{in}} = \mathbf{0}$ and $q(\phi_0^{\text{in}}) = - \partial_z \phi_0^{\text{in}}$, and integrating yields
\begin{align*}
 \nu_\Gamma c_0^{\text{in}} = - D \partial_z c_1^{\text{in}}(-\infty) - c^\ast (r(c_0^{\text{in}}) + \alpha \sigma \mu_0^{\text{in}}).
\end{align*}
Following from \eqref{eq:nu_Gamma} and \eqref{eq:match_cond_3} we get
\begin{align*}
 \nu_\Gamma (c^\ast - c_0^{\text{in}}) = D \nabla c_0^{\text{out}}(\mathbf{x}_{_{-}},t) \cdot \mathbf{n} \quad \text{for all } z \in (-\infty, \infty), 
\end{align*}
i.e., interface condition~\eqref{eq:int_cond_c}.

\textit{Expansion of $\eqref{eq:model}_5$, $\mathcal{O}(1)$}:
Inserting $\phi^{\text{in}} = \phi_0^{\text{in}} + \varepsilon \phi_1^{\text{in}} + \mathcal{O}(\varepsilon^2)$ we find
\begin{align}
 \mu_0^{\text{in}} = W''_{\text{dw}}(\phi_0^{\text{in}}) \phi_1^{\text{in}} - (\partial_{zz} \phi_1^{\text{in}} - \kappa \partial_z \phi_0^{\text{in}}) = (W''_{\text{dw}}(\phi_0^{\text{in}}) - \partial_{zz}) \phi_1^{\text{in}} + \kappa \partial_z \phi_0^{\text{in}}.
\label{eq:mu_det}
\end{align}
We have
\begin{align*}
 \int_{-\infty}^{\infty} (\partial_z \phi_0^{\text{in}}) (W''_{\text{dw}}(\phi_0^{\text{in}}) - \partial_{zz}) \phi_1^{\text{in}}\;dz = \int_{-\infty}^{\infty} \partial_z(W'_{\text{dw}}(\phi_0^{\text{in}}) - \partial_{zz} \phi_0^{\text{in}}) \phi_1^{\text{in}} \;dz = 0
\end{align*}
from partial integration, matching conditions \eqref{eq:match_cond_2} and \eqref{eq:match_cond_1}, and \eqref{eq:ODE_phi}. Multiplying \eqref{eq:mu_det} by $\partial_z \phi_0^{\text{in}}$ and integrating we thus obtain
\begin{align*}
 \int_{-\infty}^{\infty} \mu_0^{\text{in}} \partial_z \phi_0^{\text{in}} \;dz = \int_{-\infty}^{\infty} \kappa (\partial_z \phi_0^{\text{in}})^2 \;dz.
\end{align*}
$\mu_0^{\text{in}}$ is independent of $z$ due to \eqref{eq:mu_const}. Hence, by \eqref{eq:match_cond_1} we obtain $\int_{-\infty}^{\infty} \mu_0^{\text{in}} \partial_z \phi_0^{\text{in}} \;dz = - \mu_0^{\text{in}}$. Moreover, we have
\begin{align*}
 \int_{-\infty}^{\infty} (\partial_z \phi_0^{\text{in}})^2 \;dz = \int_{-\infty}^{\infty} 2 W_{\text{dw}}(\phi_0^{\text{in}}) \;dz = \int_1^0 - \sqrt{2 W_{\text{dw}}(\phi_0^{\text{in}})} \;d\phi_0^{\text{in}} = 1
\end{align*}
from integration by substitution and explicitly calculating the integral after inserting \eqref{eq:defDW}.
Therefore, we find $- \mu_0^{\text{in}} = \kappa$ and from \eqref{eq:nu_Gamma} we get 
\begin{align*}
 \nu_\Gamma = \alpha \sigma \kappa - r(c_0^{\text{in}}), 
\end{align*}
i.e., interface condition \eqref{eq:int_cond_nu}.
}

\section{Thermodynamic consistency} \label{sec:thermo}
In this section we show in Theorem \ref{theo:analysis} that smooth solutions of \eqref{eq:model} dissipate an appropriate free energy functional. The results are the basis for  the  corresponding 
discrete result on energy dissipation in Section \ref{sec:disc}.

Consider the free energy functional
\revision{
\begin{align}
 F[\phi(\cdot,t), \nabla\phi(\cdot,t), \mathbf{v}(\cdot,t), c(\cdot,t)] = \int_\Omega \frac{1}{2} \tilde\rho_f(\vecx,t)|\mathbf{v}(\vecx,t)|^2 + f(\phi(\vecx,t), \nabla\phi(\vecx,t)) + \frac{1}{\tilde \alpha} g(c(\vecx,t)) \tilde \phi_c(\vecx,t) \;d\vecx.
\label{eq:free_energy}
\end{align}
}
It comprises the kinetic energy of the fluid phase, the potential energy contribution as the  Ginzburg-Landau energy
\begin{align*}
 f(\varphi, \boldsymbol{\xi}) = \sigma \left(\frac{W_{\text{dw}}(\varphi)}{\varepsilon} + \frac{\varepsilon}{2} |\boldsymbol{\xi}|^2\right)
\end{align*}
\revision{and the free energy associated with the fluid-ion mixture.}
We readily deduce: 
\begin{theorem}\label{theo:analysis}
Let $(p, \vecv, c, \phi, \mu)$  be a  smooth solution  of the 
two-phase model 
\eqref{eq:model} 
which satisfies the boundary conditions from \eqref{icbc} \revisionT{and $F[\phi_0, \nabla\phi_0, \mathbf{v}_0, c_0] < \infty$}.\\
Then,  $(\vecv, c, \phi, \mu)$  fulfills the free energy dissipation inequality
\revision{
\begin{align}
 \frac{d}{dt} &F[\phi(\cdot,t), \nabla\phi(\cdot,t), \mathbf{v}(\cdot,t), c(\cdot,t)] \nonumber \\
 &= \int_\Omega - 2 \gamma \nabla^s \mathbf{v}(\vecx,t)\mathbin{:} \nabla \mathbf{v}(\vecx,t) - \rho \frac{1}{\varepsilon} d(\tilde \phi_f(\vecx,t)) |\mathbf{v}(\vecx,t)|^2 - \sigma M \varepsilon |\nabla\mu(\vecx,t)|^2 \nonumber
 \\
 &\qquad- \frac{1}{\tilde \alpha} D g''(c(\vecx,t)) \tilde \phi_c(\vecx,t) |\nabla c(\vecx,t)|^2 - \frac{q(\phi(\vecx,t))}{\varepsilon \tilde \alpha} (r(c(\vecx,t)) + \tilde \alpha \sigma \mu(\vecx,t))^2 \;d\vecx \nonumber \\
 &\leq 0
\label{eq:TC}
\end{align}
}
for all $t \in [0,T]$.
\end{theorem}
\begin{proof}
By basic vector calculus (cf.\ \eqref{eq:aux_1} in Lemma \ref{auxlemma} in the appendix) and  the boundary condition on $\vecv$ in \eqref{icbc} we deduce after partial integration 
\begin{align*}
 \int_\Omega \mathbf{v} \cdot \Big(\nabla \cdot ((\rho_f \mathbf{v} + \mathbf{J}_f) \otimes \mathbf{v})\Big) \;d\vecx 
 &= \int_\Omega \mathbf{v} \cdot \Big(\nabla \cdot (\rho_f \mathbf{v} + \mathbf{J}_f) \mathbf{v} + ((\rho_f \mathbf{v} + \mathbf{J}_f) \cdot \nabla) \mathbf{v}\Big) \;d\vecx \\
 &= \int_\Omega |\mathbf{v}|^2 \nabla \cdot (\rho_f \mathbf{v} + \mathbf{J}_f) + \frac{1}{2} ((\rho_f \mathbf{v} + \mathbf{J}_f) \cdot \nabla) |\mathbf{v}|^2\;d\vecx \\
 &= \int_\Omega \frac{1}{2} |\mathbf{v}|^2 \nabla \cdot (\rho_f \mathbf{v} + \mathbf{J}_f) \;d\vecx.
\end{align*}
Using the definitions of $\tilde\rho_f$ from \eqref{eq:defrhof_2}, $\rho_f$ from \eqref{eq:defrhof} and $\mathbf{J}_f$ from \eqref{eq:defJf}, we obtain with $\eqref{eq:model}_4$ the evolution equation 
\begin{align*}
 \partial_t \tilde\rho_f = \partial_t (\rho\phi + \rho\delta) = \rho \partial_t \phi = \rho R - \nabla \cdot (\rho\phi\mathbf{v} + \rho\mathbf{J}) = \rho R - \nabla \cdot (\rho_f\mathbf{v} + \mathbf{J}_f).
\end{align*}
Together with $\eqref{eq:model}_1$ and $\eqref{eq:model}_2$ the time derivative of the kinetic energy (of the fluid phase) can then be calculated as
\begin{align*}
 \frac{d}{dt} \int_\Omega \frac{1}{2} \tilde\rho_f |\mathbf{v}|^2 \;d\vecx 
 &= \int_\Omega \mathbf{v} \cdot \partial_t (\tilde\rho_f \mathbf{v}) - \frac{1}{2} |\mathbf{v}|^2 \partial_t \tilde\rho_f \;d\vecx \\
 &= \int_\Omega - \mathbf{v} \cdot \Big(\nabla \cdot ((\rho_f \mathbf{v} + \mathbf{J}_f) \otimes \mathbf{v})\Big) - \tilde \phi_f \mathbf{v} \cdot \nabla p + \mathbf{v} \cdot \big(\nabla \cdot (2 \gamma \nabla^s \mathbf{v})\big) - \rho \frac{1}{\varepsilon} d(\tilde \phi_f)|\mathbf{v}|^2 \\
 &\hspace{1.5cm} + \mathbf{\tilde S} \cdot \mathbf{v} + \mathbf{v} \cdot \frac{1}{2} \rho \mathbf{v} R - \frac{1}{2} |\mathbf{v}|^2 \rho R + \frac{1}{2} |\mathbf{v}|^2 \nabla \cdot (\rho_f \mathbf{v} + \mathbf{J}_f) \;d\vecx \\
 &= \int_\Omega \mathbf{v} \cdot \big(\nabla \cdot (2 \gamma \nabla^s \mathbf{v})\big) - \rho \frac{1}{\varepsilon} d(\tilde \phi_f)|\mathbf{v}|^2 + \mathbf{\tilde S} \cdot \mathbf{v} \;d\vecx \\
 &= \int_\Omega - 2 \gamma \nabla^s \mathbf{v} \mathbin{:} \nabla \mathbf{v} - \rho \frac{1}{\varepsilon} d(\tilde \phi_f)|\mathbf{v}|^2 + \mathbf{\tilde S} \cdot \mathbf{v} \;d\vecx,
\end{align*}
using again \eqref{icbc}. Next, for the time derivative of the Ginzburg-Landau free energy $f$ we obtain
\begin{align*}
 \frac{d}{dt} \int_\Omega f(\phi, \nabla\phi) \;d\vecx &= \int_\Omega \frac{\partial{f}}{\partial\phi}\partial_t\phi + \frac{\partial{f}}{\partial\nabla\phi} \partial_t\nabla\phi \;d\vecx = \int_\Omega \left(\frac{\partial{f}}{\partial\phi} - \nabla \cdot \frac{\partial{f}}{\partial\nabla\phi}\right) \partial_t\phi \;d\vecx \\
 &= \int_\Omega \sigma \left(\frac{W'_{\text{dw}}(\phi)}{\varepsilon} - \varepsilon \Delta\phi\right) \partial_t\phi \;d\vecx \overset{\eqref{eq:model}_5}{=} \int_\Omega \sigma \mu \partial_t\phi \;d\vecx \\
 \overset{\eqref{eq:model}_4}&{=} \int_\Omega \sigma \mu \big(R - \nabla \cdot (\phi\mathbf{v} + \mathbf{J})\big) \;d\vecx \\
 &= \int_\Omega \sigma \mu R \;d\vecx + \int_\Omega - \sigma \mu \nabla \cdot (\phi\mathbf{v}) \;d\vecx + \int_\Omega - \sigma \mu \nabla \cdot \mathbf{J} \;d\vecx.
\end{align*}
By $\eqref{eq:model}_1$ we have
\begin{align*}
 0 = \nabla \cdot (\tilde \phi_f \mathbf{v}) = \nabla \cdot \left((2 \delta + (1-2\delta) \phi) \mathbf{v}\right) = \nabla \cdot (\phi \mathbf{v}) + 2 \delta \nabla \cdot \left((1 - \phi) \mathbf{v}\right),
\end{align*}
yielding
\begin{align}
 \int_\Omega \sigma \mu \nabla \cdot (\phi \mathbf{v}) \;d\vecx &= \int_\Omega - 2\delta \sigma \mu \nabla \cdot \left((1 - \phi) \mathbf{v}\right) \;d\vecx = \int_\Omega 2\delta \sigma (1 - \phi) \nabla \mu \cdot \mathbf{v} \;d\vecx \nonumber \\
 &= \int_\Omega \mathbf{\tilde S} \cdot \mathbf{v} \;d\vecx.
\label{eq:surface_tension}
\end{align}
Hence, we obtain
\begin{align*}
 \frac{d}{dt} \int_\Omega f(\phi, \nabla\phi) \;d\vecx &= \int_\Omega \sigma \mu R \;d\vecx - \int_\Omega \mathbf{\tilde S} \cdot \mathbf{v} \;d\vecx + \int_\Omega - \sigma \mu \nabla \cdot \mathbf{J} \;d\vecx \\
 &= \int_\Omega \sigma \mu R \;d\vecx - \int_\Omega \mathbf{\tilde S} \cdot \mathbf{v} \;d\vecx + \int_\Omega - \sigma M \varepsilon |\nabla\mu|^2 \;d\vecx. 
\end{align*}
\revision{
By the boundary condition on $c$ in \eqref{icbc} we deduce after partial integration
\begin{align*}
 \int_\Omega g'(c) \nabla \cdot ((\phi_c \mathbf{v} + \mathbf{J}_c) (c - c_{\text{eq}})) \;d\vecx &= \int_\Omega \nabla(g(c)) \cdot (\phi_c \mathbf{v} + \mathbf{J}_c) + g'(c) (c - c_{\text{eq}}) \nabla \cdot (\phi_c \mathbf{v} + \mathbf{J}_c) \;d\vecx \\
 &= \int_\Omega - (g(c) - g'(c) (c - c_{\text{eq}})) \nabla \cdot (\phi_c \mathbf{v} + \mathbf{J}_c) \;d\vecx.
\end{align*}
Together with $\eqref{eq:model}_3$ and $\eqref{eq:model}_4$ we then calculate
\begin{align*}
 \int_\Omega \partial_t (g(c) \tilde \phi_c) \;d\vecx &= \int_\Omega g'(c) \partial_t c \tilde \phi_c + g(c) \partial_t \tilde \phi_c \;dx = \int_\Omega g'(c) \partial_t (c - c_{\text{eq}}) \tilde \phi_c + g(c) \partial_t \phi_c \;d\vecx \\
 &= \int_\Omega g'(c) \partial_t (\tilde \phi_c (c - c_{\text{eq}})) + (g(c) - g'(c) (c - c_{\text{eq}})) \partial_t \phi_c \;d\vecx \\
 &= \int_\Omega g'(c) \left(D \nabla \cdot (\tilde \phi_c \nabla c) + R_c - \nabla \cdot ((\phi_c \mathbf{v} + \mathbf{J}_c) (c - c_{\text{eq}}))\right) \\
 &\qquad+ (g(c) - g'(c) (c - c_{\text{eq}})) \left(R - \nabla \cdot (\phi_c \mathbf{v} + \mathbf{J}_c)\right) \;d\vecx \\
 &= \int_\Omega g'(c) D \nabla \cdot (\tilde \phi_c \nabla c) + \left(g(c) - g'(c) (c - c_{\text{eq}}) + g'(c) (c^\ast - c_{\text{eq}})\right) R \;d\vecx \\
 &= \int_\Omega D g'(c) \nabla \cdot (\tilde \phi_c \nabla c) + r(c) R \;d\vecx
 = \int_\Omega - D g''(c) \tilde \phi_c |\nabla c|^2 + r(c) R \;d\vecx.
\end{align*}
For the time derivative of the fluid-ion mixture energy we hence obtain
\begin{align*}
 \frac{d}{dt} \int_\Omega \frac{1}{\tilde \alpha} g(c) \tilde \phi_c \;d\vecx = \int_\Omega - \frac{1}{\tilde \alpha} D g''(c) \tilde \phi_c |\nabla c|^2 + \frac{1}{\tilde \alpha} r(c) R \;d\vecx.
\end{align*}
}
Thus, for the time derivative of the total free energy functional we have
\revision{
\begin{align*}
 \frac{d}{dt} F[\phi, \nabla\phi, \mathbf{v}, c] = \int_\Omega - 2 \gamma \nabla^s \mathbf{v} \mathbin{:} \nabla \mathbf{v} - \rho \frac{1}{\varepsilon} d(\tilde \phi_f) |\mathbf{v}|^2 - \sigma M \varepsilon |\nabla\mu|^2 - \frac{1}{\tilde \alpha} D g''(c) \tilde \phi_c |\nabla c|^2 + \left(\frac{1}{\tilde \alpha} r(c) + \sigma \mu\right)\!R \;d\vecx. 
\end{align*}
Inserting the definition of $R$ \eqref{eq:react} and using that $\nabla^s \mathbf{v} \mathbin{:} \nabla \mathbf{v} \geq 0$, this concludes the proof.
}
\end{proof}
\section{An energy-consistent discretization for the fluid-solid system} \label{sec:disc}
In this section we present a fully-discrete numerical method to approximate the initial-boundary value problem for  the fluid-solid dynamics in \eqref{eq:model} with initial and boundary conditions from \eqref{icbc}.
The method relies on a finite element discretization in space.  We employ a specific time-stepping method to ensure the thermodynamic consistency of the fully-discrete approach. This main result 
on the discretization is given in Theorem \ref{theo:main} below.
\subsection{Weak formulation and fully-discrete approach}
The discretization relies on a slight reformulation  of  the 
momentum and ion concentration equations 
in \eqref{eq:model}. In fact, using $\rho= \rm{const.}$, we  can rewrite   the Cahn-Hilliard evolution  in $\eqref{eq:model}_4$ in terms of the fluid density, i.e., 
\[
 \partial_t \phi + \nabla \cdot (\phi \mathbf{v} + \mathbf{J}) = R 
 \,
 \Longleftrightarrow
 \, \partial_t \left(\rho \phi + \rho\delta\right) + \nabla \cdot (\rho \phi \mathbf{v} + \rho \mathbf{J}) =  \rho R 
 \,
 \Longleftrightarrow \,  \partial_t \tilde \rho_f + \nabla \cdot ( \rho_f \mathbf{v} + \mathbf{J}_f) = \rho R. 
\]
Thus, we have   for the  inertial operator 
\begin{align*}
 \partial_t (\tilde \rho_f \mathbf{v}) &+ \nabla \cdot ((\rho_f \mathbf{v} + \mathbf{J}_f) \otimes \mathbf{v}) \\
 &= \tilde\rho_f \partial_t \mathbf{v} + \partial_t \tilde\rho_f \mathbf{v} + \nabla \cdot (\rho_f \mathbf{v} + \mathbf{J}_f) \mathbf{v} + ((\rho_f \mathbf{v} + \mathbf{J}_f) \cdot \nabla) \mathbf{v} \\
 &= \tilde\rho_f \partial_t \mathbf{v} + ((\rho_f \mathbf{v} + \mathbf{J}_f) \cdot \nabla) \mathbf{v} + 
 {\left[\partial_t \tilde\rho_f + \nabla \cdot (\rho_f \mathbf{v} + \mathbf{J}_f)\right]}
 \mathbf{v} \\
 &= \tilde\rho_f \partial_t \mathbf{v} + ((\rho_f \mathbf{v} + \mathbf{J}_f) \cdot \nabla) \mathbf{v} + \rho R \mathbf{v},
\end{align*}
and  obtain for  $\eqref{eq:model}_2$ the equivalent relation
\begin{align} \label{eq:model22}
 \tilde\rho_f \partial_t \mathbf{v} + ((\rho_f \mathbf{v} + \mathbf{J}_f) \cdot \nabla) \mathbf{v} = - \tilde\phi_f \nabla p + \nabla \cdot (2 \gamma \nabla^s \mathbf{v}) - \rho \frac{1}{\varepsilon} d(\tilde\phi_f) \mathbf{v} + \mathbf{\tilde S} - \frac{1}{2} \rho R \mathbf{v}.
\end{align}
Next, from the Cahn-Hilliard evolution in $\eqref{eq:model}_4$ we derive
\begin{align*}
 \partial_t \tilde \phi_c + \nabla \cdot (\phi_c \mathbf{v} + \mathbf{J}_c) = R
\end{align*}
using the definitions of $\phi_c$ \eqref{eq:defrhof} and $\tilde \phi_c$ \eqref{eq:defrhof_2}.
Hence, we have
\begin{align*}
 \partial_t (\tilde \phi_c (c - c_{\text{eq}})) &+ \nabla \cdot ((\phi_c \mathbf{v} + \mathbf{J}_c) (c - c_{\text{eq}})) \\
 &= \tilde \phi_c \partial_t c + \partial_t \tilde \phi_c (c - c_{\text{eq}}) + \nabla \cdot (\phi_c \mathbf{v} + \mathbf{J}_c) (c - c_{\text{eq}}) + (\phi_c \mathbf{v} + \mathbf{J}_c) \cdot \nabla c \\
 &= \tilde \phi_c \partial_t c + (\phi_c \mathbf{v} + \mathbf{J}_c) \cdot \nabla c + [\partial_t \tilde \phi_c + \nabla \cdot (\phi_c \mathbf{v} + \mathbf{J}_c)] (c - c_{\text{eq}}) \\
 &= \tilde \phi_c \partial_t c + (\phi_c \mathbf{v} + \mathbf{J}_c) \cdot \nabla c + R (c - c_{\text{eq}}),
\end{align*}
and obtain for $\eqref{eq:model}_3$ the equivalent relation
\begin{align} \label{eq:model32}
 \tilde \phi_c \partial_t c + (\phi_c \mathbf{v} + \mathbf{J}_c) \cdot \nabla c = D \nabla \cdot (\tilde \phi_c \nabla c) + (c^\ast - c) R.
\end{align}
Based on this form we introduce a weak formulation of the two-phase model.\\
We use the standard notation for Sobolev spaces, see e.g.~\cite{Brezis2011}, and introduce the space of square-integrable functions with zero mean $L_0^2(\Omega) := \{f \in L^2(\Omega)| \frac{1}{|\Omega|} \int_\Omega f \; d\vecx = 0\}$.\\
The quintuple 
\begin{equation} \label{space}
(p,\mathbf{v},c,\phi,\mu) \in L_0^2(\Omega) \times H_0^1(\Omega)^\nu \times H^1(\Omega) \times H^1(\Omega) \times H^1(\Omega)
\end{equation}
is  called a weak solution of the
two-phase model \eqref{eq:model}, with the given initial and boundary data from \eqref{icbc}, if the initial conditions in \eqref{icbc} hold a.e.~in $\Omega$ and if 
\begin{equation}
\label{eq:weak_form}
\begin{array}{rcl}
 0&= &\ds - \int_\Omega \check p \nabla \cdot (\tilde \phi_f \mathbf{v}) \;d\vecx
 ,
\\[2.0ex]
 0&=&\ds 
 \int_\Omega \tilde\rho_f \partial_t \mathbf{v} \cdot \mathbf{\check v} \;d\vecx + \int_\Omega \left(\left(\left(\rho_f \mathbf{v} + \mathbf{J}_f\right) \cdot \nabla \right) \mathbf{v}\right) \cdot \mathbf{\check v} \;d\vecx \\[2.0ex]
 & & \ds - \int_\Omega p \nabla \cdot (\tilde \phi_f \mathbf{\check v}) \;d\vecx 
 + \int_\Omega 2 \gamma \nabla^s \mathbf{v} \mathbin{:} \nabla \mathbf{\check v} \;d\vecx 
 + \int_\Omega \rho \frac{1}{\varepsilon} d(\tilde \phi_f) \mathbf{v} \cdot \mathbf{\check v} \;d\vecx - \int_\Omega \mathbf{\tilde S} \cdot \mathbf{\check v} \;d\vecx + \int_\Omega \frac{1}{2} \rho R \mathbf{v} \cdot \mathbf{\check v} \;d\vecx, 
\\[2.0ex]
 0&=& \ds
 \int_\Omega (\tilde \phi_c \partial_t c) \check c \;d\vecx + \int_\Omega ((\phi_c \mathbf{v} + \mathbf{J}_c) \cdot \nabla c) \check c \;d\vecx + \int_\Omega D (\tilde \phi_c \nabla c) \cdot \nabla \check c \;d\vecx - \int_\Omega (c^\ast - c) R \check c \;d\vecx,
\\[2.0ex]
 0&=&\ds 
 \int_\Omega (\partial_t \phi) \check \phi \;d\vecx
 - \int_\Omega (\phi \mathbf{v} + \mathbf{J}) \cdot \nabla \check \phi \;d\vecx - \int_\Omega R \check \phi \;d\vecx,
\\[2.0ex]
 0&=& \ds 
 \int_\Omega \mu \check \mu \;d\vecx
 - \int_\Omega \frac{W_{\text{dw}}'(\phi)}{\varepsilon} \check \mu \;d\vecx
 - \int_\Omega \varepsilon \nabla \phi \cdot \nabla \check \mu \;d\vecx
\end{array}
\end{equation}
holds for all $\check p \in L_0^2(\Omega)$, $\mathbf{\check v} \in H_0^1(\Omega)^\nu$ and $\check c$, $\check \phi$, $\check \mu \in H^1(\Omega)$, and a.e.~in $(0,T]$.
\medskip 

Note that the boundary conditions from \eqref{icbc} are encoded in the ansatz spaces. In what follows we assume that there exists  a unique weak solution $
(p,\mathbf{v},c,\phi,\mu) \in L_0^2(\Omega) \times H_0^1(\Omega)^\nu \times H^1(\Omega) \times H^1(\Omega) \times H^1(\Omega)$
 of the  
two-phase model \eqref{eq:model}. We conjecture that a rigorous well-posedness analysis for weak solutions can be achieved following e.g.~\cite{Abels2013, Colli2012}.

 Before discretizing \eqref{eq:weak_form}, we make \revision{three} equivalent reformulations that are essential for showing thermodynamic consistency afterwards. Firstly, as proposed in \cite{Gruen2014}, we split the convective term in $\eqref{eq:weak_form}_2$ as
\begin{align*}
 &\int_\Omega \left(\left(\left(\rho_f \mathbf{v} + \mathbf{J}_f\right) \cdot \nabla \right) \mathbf{v}\right) \cdot \mathbf{\check v} \;d\vecx \\
 &\overset{\eqref{eq:aux_3}}{=} \frac{1}{2} \int_\Omega \left(\rho_f \mathbf{v} + \mathbf{J}_f\right) \cdot \nabla \left(\mathbf{v} \cdot \mathbf{\check v}\right) d\vecx + \frac{1}{2} \int_\Omega \left(\rho_f \mathbf{v} + \mathbf{J}_f\right) \cdot \left(\left(\nabla \mathbf{v}\right) \mathbf{\check v}\right) d\vecx - \frac{1}{2} \int_\Omega \left(\rho_f \mathbf{v} + \mathbf{J}_f\right) \cdot \left(\left(\nabla \mathbf{\check v}\right) \mathbf{v}\right) d\vecx.
\end{align*}
\revision{
We proceed similarly for the convective term in $\eqref{eq:weak_form}_3$ and by adding zeros we obtain
\begin{align*}
 &\int_\Omega ((\phi_c \mathbf{v} + \mathbf{J}_c) \cdot \nabla c) \check c \;d\vecx = \int_\Omega ((\phi_c \mathbf{v} + \mathbf{J}_c) \cdot \nabla (c - c_\text{eq})) \check c \;d\vecx \\
 &= \frac{1}{2} \int_\Omega (\phi_c \mathbf{v} + \mathbf{J}_c) \cdot \nabla ((c - c_{\text{eq}}) \check c) \;d\vecx + \frac{1}{2} \int_\Omega ((\phi_c \mathbf{v} + \mathbf{J}_c) \cdot \nabla c) \check c \;d\vecx - \frac{1}{2} \int_\Omega ((\phi_c \mathbf{v} + \mathbf{J}_c) \cdot \nabla \check c) (c - c_{\text{eq}}) \;d\vecx.
\end{align*}
}
\revision{Thirdly}, for the surface tension term in $\eqref{eq:weak_form}_2$ we use \eqref{eq:surface_tension} and partial integration to obtain
\begin{align}
 \int_\Omega \mathbf{\tilde S} \cdot \mathbf{\check v} \;d\vecx = \int_\Omega \sigma \mu \nabla \cdot (\phi \mathbf{\check v}) \;d\vecx = - \int_\Omega \sigma (\phi \mathbf{\check v}) \cdot \nabla \mu \;d\vecx.
\label{eq:S_weak}
\end{align}
Mimicking the continuous case, all three reformulations are needed to eventually cancel contributions from the convective term in $\eqref{eq:weak_form}_4$. \revisionT{In particular, we also establish the weak definition of $\mathbf{\tilde S}$ in our numerical scheme.}

To introduce our numerical method for weak solutions,  we employ a Galerkin finite element method on a conforming triangular mesh $\mathcal{T}_h$ that consists of closed simplices $\{T\}_{T \in \mathcal{T}_h}$ with spatial measures $h_T = \text{diam}(T)$ and $h = \max_{T \in \mathcal{T}_h} h_T$, such that $\overline \Omega = \cup_{T \in \mathcal{T}_h} \overline T$. 
As discrete ansatz spaces we choose Lagrange finite element spaces denoted by $\mathcal{V}^p_h$, $\mathcal{V}^{\mathbf{v}}_h$, $\mathcal{V}^c_h$ and $\mathcal{V}^{\text{ch}}_h$ (for $\phi$ and $\mu$), respectively with
\begin{align} \label{eq:FE_spaces}
 \mathcal{V}^{\text{ch}}_h &= \left\{\eta_h \in C(\overline \Omega) \Big| \text{ } \eta_h|_T \in \mathcal{P}^{i}(T) \text{ for all } T \in \mathcal{T}_h\right\}, \nonumber \\
 \mathcal{V}^c_h &= \left\{\eta_h \in C(\overline \Omega) \Big| \text{ } \eta_h|_T \in \mathcal{P}^{j}(T) \text{ for all } T \in \mathcal{T}_h\right\},
 \nonumber \\
 \mathcal{V}^p_h &= \left\{p_h \in C(\overline \Omega) \Big| \text{ } p_h|_T \in \mathcal{P}^{k}(T) \text{ for all } T \in \mathcal{T}_h\right\} \cap L_0^2(\Omega), \nonumber \\
 \mathcal{V}^{\mathbf{v}}_h &= \left\{\mathbf{v}_h \in C(\overline \Omega)^\nu \Big| \text{ } \mathbf{v}_h|_T \in \left(\mathcal{P}^{l}(T)\right)^\nu \text{ for all } T \in \mathcal{T}_h, \quad \mathbf{v}_h|_{\partial \Omega} = 0\right\}.
\end{align}
Here,  $\mathcal{P}^l(T)$ denotes the
set of polynomials of degree at most $l$ on the simplex $T$. The choice of the Lagrangian finite element spaces is for the sake of definiteness. For our main result in Theorem \ref{theo:main} we only require the ansatz spaces to be embedded into the ansatz space $ L_0^2(\Omega) \times H_0^1(\Omega)^\nu \times H^1(\Omega) \times H^1(\Omega) \times H^1(\Omega) $ from \eqref{space}. Let us also note that \revision{for the generic case where the polynomial exponents $i,j,k,l$ are not related we need to introduce projections to the numerical scheme in order to prove Theorem \ref{theo:main}. This is further detailed in Remark~\ref{rem:proj}.}  For conditions that imply stability for the velocity-pressure pair we refer to \cite{Verfürth2010} and for the specific spaces used in our work to Section \ref{sec:detail} below. \revision{Moreover, let $P_h$ denote the orthogonal $L^2$-projection onto $\mathcal{V}^{\text{ch}}_h$.}

In time we divide the temporal domain $[0,T]$ into subintervals $[t^n,t^{n+1}[$ for $n = 0,\ldots,N_t-1$ with $t^{n+1} = t^n + \tau$ for constant time step $\tau = \frac{T}{N_t} > 0$. Given the values for the primary variables $p^n$, $\mathbf{v}^n$, $c^n$, $\phi^n$ and $\mu^n$ at time $t^n$, we then search for $p^{n+1}$, $\mathbf{v}^{n+1}$, $c^{n+1}$, $\phi^{n+1}$ and $\mu^{n+1}$ in the corresponding finite element spaces. For the time discretization, we employ a tailor-made semi-implicit Euler method in order to account for thermodynamic consistency\revision{, see also Remark~\ref{rem:si_td} for more details}. The extension to higher-order in time methods remains an open challenge. \\
For the double-well potential $W_{\text{dw}}(\phi)$ we employ a convex-concave split of the form 
\begin{align}
 W_{\text{dw}}(\phi) = W_c(\phi) + W_e(\phi)
\label{eq:conv_conc_split}
\end{align} 
with convex part $W_c(\phi)$ and concave part $W_e(\phi)$. The convex part is then evaluated implicitly and the concave part explicitly following \cite{Eyre1998}.

With these preparations the fully-discrete method is given as follows.
Let $(p^0, \vecv^0, c^0, \phi^0, \mu^0) \in \mathcal{V}^p_h \times \mathcal{V}^{\mathbf{v}}_h \times \mathcal{V}^c_h \times \mathcal{V}^{\text{ch}}_h \times \mathcal{V}^{\text{ch}}_h$ denote the discrete initial data that are determined from $L^2$-projection to the finite element spaces. 
\\
We search for \revision{$n = 0,\ldots,N_t-1$ the approximations 
 $ (p^{n+1},\vecv^{n+1}, c^{n+1}, \phi^{n+1},\mu^{n+1}) \in  \mathcal{V}^p_h \times \mathcal{V}^{\mathbf{v}}_h \times \mathcal{V}^c_h\times \mathcal{V}^{\text{ch}}_h \times \mathcal{V}^{\text{ch}}_h  $} such that
\begin{equation}
\label{eq:disc_weak_form}
\begin{array}{rcl}
 0 &= & \ds - \int_\Omega \check p \nabla \cdot (\tilde \phi_f^{n+1} \mathbf{v}^{n+1}) \;d\vecx
 ,
\\[2ex]
 0 &=&\ds  \int_\Omega \frac{\tilde \rho_f^{n} + \tilde \rho_f^{n+1}}{2} \frac{\mathbf{v}^{n+1} - \mathbf{v}^n}{\tau} \cdot \mathbf{\check v} \;d\vecx
 + \frac{1}{2} \int_\Omega \left(\rho_f^{n} \mathbf{v}^{n+1} + \mathbf{J}_f^{n+1}\right) \cdot \nabla \revision{P_h}\!\left(\mathbf{v}^{n} \cdot \mathbf{\check v}\right) d\vecx \\[2ex]
 &&\ds \quad + \frac{1}{2} \int_\Omega \left(\rho_f^{n} \mathbf{v}^{n} + \mathbf{J}_f^{n}\right) \cdot \left(\left(\nabla \mathbf{v}^{n+1}\right) \mathbf{\check v}\right) d\vecx
 - \frac{1}{2} \int_\Omega \left(\rho_f^{n} \mathbf{v}^{n} + \mathbf{J}_f^{n}\right) \cdot \left(\left(\nabla \mathbf{\check v}\right) \mathbf{v}^{n+1}\right) d\vecx \\[2ex]
 &&\ds \quad 
 - \int_\Omega p^{n+1} \nabla \cdot (\tilde \phi_f^{n+1} \mathbf{\check v}) \;d\vecx 
 + \int_\Omega 2 \gamma \nabla^s \mathbf{v}^{n+1} \mathbin{:} \nabla \mathbf{\check v} \;d\vecx + \int_\Omega \rho \frac{1}{\varepsilon} d(\tilde \phi_f^{n}) \mathbf{v}^{n+1} \cdot \mathbf{\check v} \;d\vecx \\[2ex]
 &&\ds \quad - \int_\Omega \mathbf{\tilde S}^{n+1} \cdot \mathbf{\check v} \;d\vecx + \int_\Omega \frac{1}{2} \rho R^{n+1} \revision{P_h}(\mathbf{v}^{n} \cdot \mathbf{\check v}) \;d\vecx,
\\[2ex]
0 &=&\ds  \int_\Omega \frac{\tilde \phi_c^n + \tilde \phi_c^{n+1}}{2} \frac{c^{n+1} - c^n}{\tau} \check c \;d\vecx + \frac{1}{2} \int_\Omega \left(\phi_c^{n} \mathbf{v}^{n+1} + \mathbf{J}_c^{n+1}\right) \cdot \nabla \revision{P_h}\!\left((c^{n} - c_{\text{eq}}) \check c\right) d\vecx \\[2ex]
 &&\ds \quad + \frac{1}{2} \int_\Omega \left((\phi_c^n \mathbf{v}^n + \mathbf{J}_c^n) \cdot \nabla c^{n+1}\right) \check c \;d\vecx - \frac{1}{2} \int_\Omega \left((\phi_c^n \mathbf{v}^n + \mathbf{J}_c^n) \cdot \nabla \check c\right) (c^{n+1} - c_{\text{eq}}) \;d\vecx \\[2ex]
 &&\ds \quad + \int_\Omega D (\tilde \phi_c^n \nabla c^{n+1}) \cdot \nabla \check c \;d\vecx - \int_\Omega (c^\ast - c^n) R^{n+1} \check c \;d\vecx \revision{- \int_\Omega \frac{1}{2} R^{n+1} (I - P_h) ((c^n - c_{\text{eq}}) \check c) \;d\vecx},
\\[2ex]
 0 &=&\ds 
 \int_\Omega \frac{\phi^{n+1}-\phi^n}{\tau} \check \phi \;d\vecx
 + \int_\Omega M \varepsilon \nabla \mu^{n+1} \cdot \nabla \check \phi \;d\vecx - \int_\Omega \left(\phi^{n} \mathbf{v}^{n+1}\right) \cdot \nabla \check \phi \;d\vecx - \int_\Omega R^{n+1} \check \phi \;d\vecx, 
\\[2ex]
 0 &=& \ds 
 \int_\Omega \mu^{n+1} \check \mu \;d\vecx
 - \int_\Omega \frac{W_c'(\phi^{n+1}) + W_e'(\phi^{n})}{\varepsilon} \check \mu \;d\vecx
 - \int_\Omega \varepsilon \nabla \phi^{n+1} \cdot \nabla \check \mu \;d\vecx
\end{array}
\end{equation}
holds for all $\check p \in \mathcal{V}^p_h$, $\mathbf{\check v} \in \mathcal{V}^\vecv_h $, $\check c \in \mathcal{V}^c_h$ and $\check \phi$, $\check \mu \in \mathcal{V}^{\text{ch}}_h$ with
\begin{align*}
 \int_\Omega \mathbf{\tilde S}^{n+1} \cdot \mathbf{\check v} \;d\vecx = - \int_\Omega \sigma (\phi^n \mathbf{\check v}) \cdot \nabla \mu^{n+1} \;d\vecx
\end{align*}
and
\begin{align*}
 \mathbf{J}_f^{n+1} = \rho \mathbf{J}^{n+1} = - \rho M \varepsilon \nabla \mu^{n+1}.
\end{align*}
\revision{The explicit form of $R^{n+1}$ discretizing \eqref{eq:react} is given in \eqref{eq:R_disc} (and \eqref{eq:react_disc}).}
The terms $\tilde \phi_f^{n+1}$, $\tilde \phi_c^{n+1}$ and $\tilde \rho_f^{n+1}$ are to be understood as in their original definition, with all variables replaced by the finite element functions at time $t_{n+1}$, i.e., $\tilde \phi_f^{n+1} = 2 \delta + (1-2\delta) \phi^{n+1}$, $\tilde \phi_c^{n+1} = \phi^{n+1} + \delta$ and $\tilde \rho_f^{n+1} = \rho \big(\phi^{n+1} + \delta\big)$. \revisionT{The weak definition of $\mathbf{\tilde S}^{n+1}$ is motivated by the computation in \eqref{eq:S_weak}.}

\revision{
\begin{remark}[Semi-implicit time discretization]
 The time discretization is tailored such that we can mimic the continuous case in order to prove discrete thermodynamic consistency. In principle, choices such as the time averaging of $\tilde \rho_f$ in $\eqref{eq:disc_weak_form}_2$ and $\tilde \phi_c$ in $\eqref{eq:disc_weak_form}_3$, and the splitting of the convective terms in the respective equations are essential for the proof. There are, however, also terms such as for the state-dependent dissipation $d(\tilde \phi_f)$ that can be discretized explicitly or implicitly without effect on the result. From a solver perspective, for these we favor the explicit form due to the (partial) linearization implied by this choice.
\label{rem:si_td}
\end{remark}

\begin{remark}[Projection]
 The need for the projection $P_h$ in the respective terms becomes apparent in the proof of thermodynamic consistency in Section~\ref{sec:disc_thermo}. It results from testing $\eqref{eq:disc_weak_form}_4$ with products of $\mathbf{v}$ and $c$. Using $i = 2 \max\{l, j\}$ in \eqref{eq:FE_spaces} one would obtain $P_h = I$. Dependent on the polynomial exponents, this, however, may lead to a remarkable increase in the number of unknowns and hence the computational need. Hence, we refrain from doing so. Conversely, in Section~\ref{sec:proj_free} we show that we can rewrite the discrete scheme \eqref{eq:disc_weak_form} in a projection-free form.
\label{rem:proj}
\end{remark}
}

We compose  the discrete solution $(p_h,\vecv_h, \phi_h,\mu_h): \Omega\times [0,T] \to \setR \times \setR^\nu \times \setR \times \setR  $
through 
\[(p_h,\vecv_h, c_h, \phi_h,\mu_h)(\vecx,t) =  
  (p^n,\vecv^n, c^n, \phi^n,\mu^n) (\vecx,t) \qquad \left(t \in [t^n,t^{n+1}[\right).
\]
As for the weak solution we assume in what follows that the discrete solution $(p_h,\vecv_h, c_h, \phi_h,\mu_h)$ is uniquely determined and refer to 
\cite{Garcke2016}
for analytical background. 
\subsection{Thermodynamic consistency of the fully-discrete solution} \label{sec:disc_thermo}

We now establish the discrete thermodynamic consistency by showing that the fully-discrete method satisfies an analog of the free energy dissipation inequality \eqref{eq:TC}. \revision{Due to the occurrence of $g''(c)$ in the energy estimate \eqref{eq:TC} we refrain from establishing the proof for a generic, convex function $g(c)$. Conversely, we make a special choice for $g(c)$ (and hence by \eqref{eq:react} also for $r(c)$). We choose the convex function
\begin{align*}
 g(c) = \frac{1}{2} \frac{k_c}{c_{\text{eq}}^2} (c - c_{\text{eq}})^2
\end{align*}
leading to a reaction rate of
\begin{align*}
 r(c) = \frac{1}{2} \frac{k_c}{c_{\text{eq}}^2} \left[(c - c_{\text{eq}})^2 + 2 (c - c_{\text{eq}}) (c^\ast -c)\right]
\end{align*}
following the definition in \eqref{eq:react}. Here, $k_c$ denotes the reaction rate constant.
The free-energy contribution from the fluid-ion mixture thus takes the form
\begin{align*}
 \int_\Omega \frac{1}{\tilde \alpha} g(c) \tilde \phi_c \;dx = \int_\Omega \frac{1}{2 \tilde \alpha} \frac{k_c}{c_{\text{eq}}^2} \tilde \phi_c (c - c_{\text{eq}})^2 \;dx.
\end{align*}
We depict the specific choices for $g(c)$ and $r(c)$ in Figure~\ref{fig:react}.

\begin{figure}[h]
 \centering
 \includegraphics[width=0.5\linewidth]{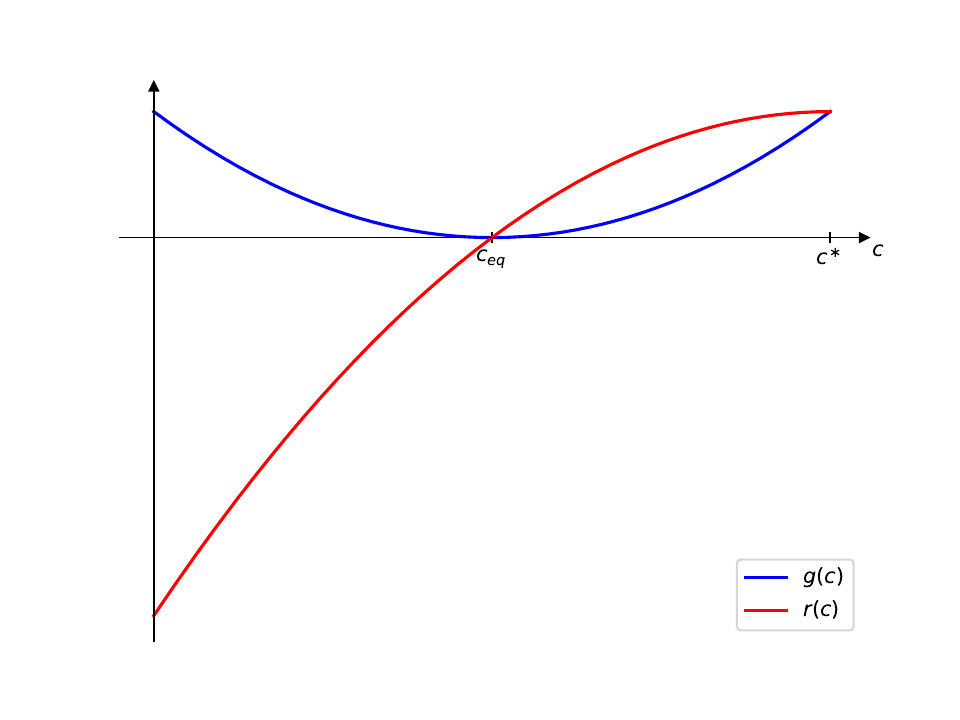}
 \caption{Schematic plot of the reaction kinetics. All axes are scaled linearly.}
 \label{fig:react}
\end{figure}

For the discretization of the reaction term we employ
\begin{align}
 R^{n+1} = - \frac{q(\phi^n)}{\varepsilon}\left(r(c^{n,n+1}) + \tilde \alpha \sigma \mu^{n+1}\right)
\label{eq:R_disc}
\end{align}
with semi-implicit reaction rate
\begin{align}
 r(c^{n,n+1}) = \frac{1}{2} \frac{k_c}{c_{\text{eq}}^2} \left[(c^{n} - c_{\text{eq}}) (c^{n+1} - c_{\text{eq}}) + 2 (c^{n+1} - c_{\text{eq}}) (c^\ast - c^{n})\right].
\label{eq:react_disc}
\end{align}

\begin{remark}[Reaction rate]
 The reaction rate $r(c) = \frac{1}{2} \frac{k_c}{c_{\text{eq}}^2} \left[(c - c_{\text{eq}})^2 + 2 (c - c_{\text{eq}}) (c^\ast - c)\right]$ (obtained by the special choice of $g(c)$) is identical to the one presented in \cite{Lamperti2025} for $k_c = 2k$. While there $r(c)$ follows from (consistency with) the sharp-interface limit, we establish this choice in order to obtain a specific form of the fluid-ion mixture energy that allows us to prove discrete thermodynamic consistency for the numerical scheme \eqref{eq:disc_weak_form}.
\end{remark}
}

\begin{theorem} \label{theo:main}
Assume that  \eqref{eq:disc_weak_form}  can be solved by a function  $(p_h,\vecv_h, c_h, \phi_h,\mu_h)(\cdot,t)  \in  \mathcal{V}^p_h \times \mathcal{V}^{\mathbf{v}}_h \times \mathcal{V}^c_h \times \mathcal{V}^{\text{ch}}_h \times \mathcal{V}^{\text{ch}}_h$ for $t\in [0,T]$ \revisionT{with $\phi_h(\cdot,t) \in (- \delta, 1 + \delta)$}.\\
Then, the discrete free energy dissipation inequality
\revision{
\begin{align}
 &\frac{1}{\tau} \left(F[\phi^{n+1}, \nabla\phi^{n+1}, \mathbf{v}^{n+1}, c^{n+1}] -  F[\phi^{n}, \nabla\phi^{n}, \mathbf{v}^{n}, c^n]\right) \nonumber \\
 &\qquad\leq \int_\Omega - 2 \gamma \nabla^s \mathbf{v}^{n+1} \mathbin{:} \nabla \mathbf{v}^{n+1} - \rho \frac{1}{\varepsilon} d(\tilde \phi_f^{n}) \left|\mathbf{v}^{n+1}\right|^2 - \sigma M \varepsilon \left|\nabla \mu^{n+1}\right|^2 \nonumber \\
 &\qquad \qquad- \frac{1}{\tilde \alpha} D \frac{k_c}{c_{\text{eq}}^2} \tilde \phi_c^n \left|\nabla c^{n+1}\right|^2 - \frac{q(\phi^n)}{\varepsilon \tilde \alpha} \left(r(c^{n,n+1}) + \tilde \alpha \sigma \mu^{n+1}\right)^2 \;d\vecx \nonumber \\
 &\qquad \leq 0 
\label{eq:TC_disc}
\end{align}
}
is satisfied for  all  
$n \in \left\{0,\ldots,N_t-1\right\}$, and hence
\begin{align*}
 F[\phi^{n+1}, \nabla\phi^{n+1}, \mathbf{v}^{n+1}, c^{n+1}]\leq F[\phi^n, \nabla\phi^n, \mathbf{v}^n, c^n].
\end{align*}
\end{theorem}
\begin{proof}
We test $\eqref{eq:disc_weak_form}_4$ with $\check \phi = \sigma \mu^{n+1}$, yielding
\begin{align}
 0 =
 \int_\Omega \sigma \frac{\phi^{n+1}-\phi^n}{\tau} \mu^{n+1} \;d\vecx
 + \int_\Omega \sigma M \varepsilon \left|\nabla \mu^{n+1}\right|^2 \;d\vecx - \int_\Omega \sigma \left(\phi^{n} \mathbf{v}^{n+1}\right) \cdot \nabla \mu^{n+1} \;d\vecx - \int_\Omega \sigma \mu^{n+1} R^{n+1} \;d\vecx,
\label{eq:femmodel_4_tested_mu}
\end{align}
and for $\check \mu = - \sigma ({\phi^{n+1}-\phi^n})/{\tau}$ in $\eqref{eq:disc_weak_form}_5$ we obtain
\[
\begin{array}{c}
 \ds 0 = 
 - \int_\Omega \mu^{n+1} \sigma \frac{\phi^{n+1}-\phi^n}{\tau} \;d\vecx
 + \int_\Omega \frac{W_c'(\phi^{n+1}) + W_e'(\phi^{n})}{\varepsilon} \sigma \frac{\phi^{n+1}-\phi^n}{\tau} \;d\vecx\\[2ex]
\ds \hspace*{2cm} + \int_\Omega \varepsilon \nabla \phi^{n+1} \cdot \nabla \left(\sigma \frac{\phi^{n+1} - \phi^n}{\tau}\right) d\vecx.
\end{array}
\]
By the convexity of $W_c(\phi)$ and the concavity of $W_e(\phi)$ from 
\eqref{eq:conv_conc_split} we have
\begin{align*}
 W_c'(\phi^{n+1}) (\phi^{n+1} - \phi^n) \geq W_c(\phi^{n+1}) - W_c(\phi^n) \qquad \text{and} \qquad W_e'(\phi^{n}) (\phi^{n+1} - \phi^n) \geq W_e(\phi^{n+1}) - W_e(\phi^n)
\end{align*}
respectively, implying
\begin{align*}
 \left(W_c'(\phi^{n+1}) + W_e'(\phi^{n})\right) (\phi^{n+1} - \phi^n) \geq W_\text{dw}(\phi^{n+1}) - W_\text{dw}(\phi^n).
\end{align*}
Hence, by the monotonicity of the integral we get
\begin{align*}
 \int_\Omega \frac{W_c'(\phi^{n+1}) + W_e'(\phi^{n})}{\varepsilon} \frac{\phi^{n+1}-\phi^n}{\tau} \;d\vecx \geq \int_\Omega \frac{1}{\tau} \frac{W_\text{dw}(\phi^{n+1}) - W_\text{dw}(\phi^n)}{\varepsilon} \;d\vecx.
\end{align*}
Moreover, utilizing \eqref{eq:aux_2} we obtain
\begin{align*}
 \int_\Omega \varepsilon \nabla \phi^{n+1} \cdot \frac{\nabla \phi^{n+1}-\nabla \phi^n}{\tau} \;d\vecx = \int_\Omega \frac{1}{\tau} \frac{\varepsilon}{2} |\nabla \phi^{n+1}|^2 \;d\vecx - \int_\Omega \frac{1}{\tau} \frac{\varepsilon}{2} |\nabla \phi^{n}|^2 \;d\vecx + \int_\Omega \frac{1}{\tau} \frac{\varepsilon}{2} |\nabla \phi^{n+1} - \nabla \phi^n|^2 \;d\vecx
\end{align*}
and thus
\begin{align}
\begin{split}
 0 &\geq - \sigma \int_\Omega \frac{\phi^{n+1}-\phi^n}{\tau} \mu^{n+1} \;d\vecx + \sigma \left(\int_\Omega \frac{1}{\tau} \frac{W_\text{dw}(\phi^{n+1}) - W_\text{dw}(\phi^n)}{\varepsilon} \;d\vecx \right.\\
 &\qquad\left.+ \int_\Omega \frac{1}{\tau} \frac{\varepsilon}{2} |\nabla \phi^{n+1}|^2 \;d\vecx - \int_\Omega \frac{1}{\tau} \frac{\varepsilon}{2} |\nabla \phi^{n}|^2 \;d\vecx + \int_\Omega \frac{1}{\tau} \frac{\varepsilon}{2} |\nabla \phi^{n+1} - \nabla \phi^n|^2 \;d\vecx\right).
\label{eq:femmodel_5_tested}
\end{split}
\end{align}
Next, testing ${\eqref{eq:disc_weak_form}}_2$ with $\mathbf{\check v} = \mathbf{v}^{n+1}$ yields
\begin{align*}
 0 &=
 \int_\Omega \frac{\tilde \rho_f^{n} + \tilde \rho_f^{n+1}}{2} \frac{\mathbf{v}^{n+1} - \mathbf{v}^n}{\tau} \cdot \mathbf{v}^{n+1} \;d\vecx + \frac{1}{2} \int_\Omega \left(\rho_f^{n} \mathbf{v}^{n+1} + \mathbf{J}_f^{n+1}\right) \cdot \nabla \revision{P_h}\!\left(\mathbf{v}^{n} \cdot \mathbf{v}^{n+1}\right) d\vecx \\
 &\qquad 
 - \int_\Omega p^{n+1} \nabla \cdot (\tilde \phi_f^{n+1} \mathbf{v}^{n+1}) \;d\vecx 
 + \int_\Omega 2 \gamma \nabla^s \mathbf{v}^{n+1} \mathbin{:} \nabla \mathbf{v}^{n+1} \;d\vecx + \int_\Omega \rho \frac{1}{\varepsilon} d(\tilde \phi_f^{n}) \left|\mathbf{v}^{n+1}\right|^2 \;d\vecx \\ 
 &\qquad
 - \int_\Omega \mathbf{\tilde S}^{n+1} \cdot \mathbf{v}^{n+1} \;d\vecx + \int_\Omega \frac{1}{2} \rho R^{n+1} \revision{P_h}(\mathbf{v}^{n} \cdot \mathbf{v}^{n+1}) \;d\vecx.
\end{align*}
We have the algebraic relation
\begin{align*}
 (\tilde \rho_f^{n} + \tilde \rho_f^{n+1})(\mathbf{v}^{n+1} - \mathbf{v}^n) \cdot \mathbf{v}^{n+1} &= \tilde \rho_f^{n+1} |\mathbf{v}^{n+1}|^2 - \tilde \rho_f^{n+1} \mathbf{v}^{n} \cdot \mathbf{v}^{n+1} + \tilde \rho_f^{n} |\mathbf{v}^{n+1} - \mathbf{v}^n|^2 + \tilde \rho_f^{n} \mathbf{v}^{n+1} \cdot \mathbf{v}^{n} - \tilde \rho_f^{n} |\mathbf{v}^{n}|^2 \\[1.5ex]
 &= \tilde \rho_f^{n+1} |\mathbf{v}^{n+1}|^2 - \tilde \rho_f^{n} |\mathbf{v}^{n}|^2 + \tilde \rho_f^{n} |\mathbf{v}^{n+1} - \mathbf{v}^n|^2 - (\tilde \rho_f^{n+1} - \tilde \rho_f^{n}) \mathbf{v}^{n} \cdot \mathbf{v}^{n+1}
\end{align*}
and hence 
\begin{align*}
 \int_\Omega \frac{\tilde \rho_f^{n} + \tilde \rho_f^{n+1}}{2} \frac{\mathbf{v}^{n+1} - \mathbf{v}^n}{\tau} \cdot \mathbf{v}^{n+1} \;d\vecx &= \int_\Omega \frac{1}{\tau} \frac{1}{2} \tilde \rho_f^{n+1} |\mathbf{v}^{n+1}|^2 \;d\vecx - \int_\Omega \frac{1}{\tau} \frac{1}{2} \tilde \rho_f^{n} |\mathbf{v}^{n}|^2 \;d\vecx \\
 &\qquad+ \int_\Omega \frac{1}{\tau} \frac{1}{2} \tilde \rho_f^{n} |\mathbf{v}^{n+1} - \mathbf{v}^n|^2 \;d\vecx - \int_\Omega \frac{1}{2} \frac{\tilde \rho_f^{n+1} - \tilde \rho_f^{n}}{\tau} \mathbf{v}^{n} \cdot \mathbf{v}^{n+1} \;d\vecx.
\end{align*}
From testing $\eqref{eq:disc_weak_form}_1$ with $\check p = p^{n+1}$, we thus obtain
\begin{align}
\begin{split}
 0 &= \int_\Omega \frac{1}{\tau} \frac{1}{2} \tilde \rho_f^{n+1} |\mathbf{v}^{n+1}|^2 \;d\vecx - \int_\Omega \frac{1}{\tau} \frac{1}{2} \tilde \rho_f^{n} |\mathbf{v}^{n}|^2 \;d\vecx + \int_\Omega \frac{1}{\tau} \frac{1}{2} \tilde \rho_f^{n} |\mathbf{v}^{n+1} - \mathbf{v}^n|^2 \;d\vecx - \int_\Omega \frac{1}{2} \frac{\tilde \rho_f^{n+1} - \tilde \rho_f^{n}}{\tau} \mathbf{v}^{n} \cdot \mathbf{v}^{n+1} \;d\vecx \\
 &\qquad + \frac{1}{2} \int_\Omega \left(\rho_f^{n} \mathbf{v}^{n+1} + \mathbf{J}_f^{n+1}\right) \cdot \nabla \revision{P_h}\!\left(\mathbf{v}^{n} \cdot \mathbf{v}^{n+1}\right) d\vecx
 + \int_\Omega 2 \gamma \nabla^s \mathbf{v}^{n+1} \mathbin{:} \nabla \mathbf{v}^{n+1} \;d\vecx + \int_\Omega \rho \frac{1}{\varepsilon} d(\tilde \phi_f^{n}) \left|\mathbf{v}^{n+1}\right|^2 \;d\vecx \\
 &\qquad - \int_\Omega \mathbf{\tilde S}^{n+1} \cdot \mathbf{v}^{n+1} \;d\vecx + \int_\Omega \frac{1}{2} \rho R^{n+1} \revision{P_h}(\mathbf{v}^{n} \cdot \mathbf{v}^{n+1}) \;d\vecx.
\label{eq:femmodel_2_tested}
\end{split}
\end{align}
Note that by definition, see \eqref{eq:defrhof} and \eqref{eq:defrhof_2}, we have $\tilde \rho_f^{n+1} - \tilde \rho_f^{n} = \rho_f^{n+1} + \rho\delta - (\rho_f^{n} + \rho\delta) = \rho (\phi^{n+1} - \phi^{n})$ and therefore
\begin{align*}
 \int_\Omega \frac{1}{2} \frac{\tilde \rho_f^{n+1} - \tilde \rho_f^{n}}{\tau} \mathbf{v}^{n} \cdot \mathbf{v}^{n+1} \;d\vecx = \int_\Omega \frac{1}{2} \rho \frac{\phi^{n+1} - \phi^{n}}{\tau} \mathbf{v}^{n} \cdot \mathbf{v}^{n+1} \;d\vecx.
\end{align*}
Using $\check \phi = \frac{1}{2} \rho \revision{P_h} (\mathbf{v}^{n} \cdot \mathbf{v}^{n+1})$ in  $\eqref{eq:disc_weak_form}_4$ \revision{and the orthogonality of $P_h$} yields
\begin{align}
 0 &=
 \int_\Omega \frac{\phi^{n+1}-\phi^n}{\tau} \frac{1}{2} \rho \mathbf{v}^{n} \cdot \mathbf{v}^{n+1} \;d\vecx
 + \int_\Omega M \varepsilon \nabla \mu^{n+1} \cdot \nabla \left(\frac{1}{2} \rho \revision{P_h} (\mathbf{v}^{n} \cdot \mathbf{v}^{n+1})\right) d\vecx \nonumber \\
 &\qquad- \int_\Omega \left(\phi^{n} \mathbf{v}^{n+1}\right) \cdot \nabla \left(\frac{1}{2} \rho \revision{P_h} (\mathbf{v}^{n} \cdot \mathbf{v}^{n+1})\right) d\vecx - \int_\Omega \frac{1}{2} \rho R^{n+1} \revision{P_h}(\mathbf{v}^{n} \cdot \mathbf{v}^{n+1}) \;d\vecx\nonumber \\
 &= \int_\Omega \frac{1}{2} \frac{\tilde \rho_f^{n+1} - \tilde \rho_f^{n}}{\tau} \mathbf{v}^{n} \cdot \mathbf{v}^{n+1} \;d\vecx - \frac{1}{2} \int_\Omega \left(\rho_f^{n} \mathbf{v}^{n+1} + \mathbf{J}_f^{n+1}\right) \cdot \nabla \revision{P_h}\!\left(\mathbf{v}^{n} \cdot \mathbf{v}^{n+1}\right) d\vecx \nonumber \\
 &\qquad- \int_\Omega \frac{1}{2} \rho R^{n+1} \revision{P_h}(\mathbf{v}^{n} \cdot \mathbf{v}^{n+1}) \;d\vecx,
\label{eq:femmodel_4_tested_v}
\end{align}
as $\rho_f^{n} = \rho \phi^n$ and $\mathbf{J}_f^{n+1} = - \rho M \varepsilon \nabla \mu^{n+1}$.
Adding \eqref{eq:femmodel_2_tested} and \eqref{eq:femmodel_4_tested_v} we have
\begin{align}
 \begin{split}
 0 &= \int_\Omega \frac{1}{\tau} \frac{1}{2} \tilde \rho_f^{n+1} |\mathbf{v}^{n+1}|^2 \;d\vecx - \int_\Omega \frac{1}{\tau} \frac{1}{2} \tilde \rho_f^{n} |\mathbf{v}^{n}|^2 \;d\vecx + \int_\Omega \frac{1}{\tau} \frac{1}{2} \tilde \rho_f^{n} |\mathbf{v}^{n+1} - \mathbf{v}^n|^2 \;d\vecx \\
 &\qquad+ \int_\Omega 2 \gamma \nabla^s \mathbf{v}^{n+1} \mathbin{:} \nabla \mathbf{v}^{n+1} \;d\vecx + \int_\Omega \rho \frac{1}{\varepsilon} d(\tilde \phi_f^{n}) \left|\mathbf{v}^{n+1}\right|^2 \;d\vecx - \int_\Omega \mathbf{\tilde S}^{n+1} \cdot \mathbf{v}^{n+1} \;d\vecx.
\label{eq:femmodel_2_tested_added}
\end{split}
\end{align}

\revision{
We test $\eqref{eq:disc_weak_form}_3$ with $\check c = \frac{k_c}{c_{\text{eq}}^2} (c^{n+1} - c_{\text{eq}})$, yielding
\begin{align*}
 0 &= \int_\Omega \frac{k_c}{c_{\text{eq}}^2} \frac{\tilde \phi_c^n + \tilde \phi_c^{n+1}}{2} \frac{c^{n+1} - c^n}{\tau} (c^{n+1} - c_{\text{eq}}) \;d\vecx \\
 &\qquad+ \frac{1}{2} \int_\Omega \frac{k_c}{c_{\text{eq}}^2} \left(\phi_c^{n} \mathbf{v}^{n+1} + \mathbf{J}_c^{n+1}\right) \cdot \nabla \revision{P_h}\!\left((c^{n} - c_{\text{eq}}) (c^{n+1} - c_{\text{eq}})\right) d\vecx + \int_\Omega D \frac{k_c}{c_{\text{eq}}^2} (\tilde \phi_c^n \nabla c^{n+1}) \cdot \nabla c^{n+1} \;d\vecx \\
 &\qquad- \int_\Omega \frac{k_c}{c_{\text{eq}}^2} (c^\ast - c^n) R^{n+1} (c^{n+1} - c_{\text{eq}}) \;d\vecx \revision{- \int_\Omega \frac{1}{2} \frac{k_c}{c_{\text{eq}}^2} R^{n+1} (I - P_h) ((c^n - c_{\text{eq}}) (c^{n+1} - c_{\text{eq}})) \;d\vecx}.
\end{align*}
We add a zero to the time evolution and similarly to before obtain the algebraic relation
\begin{align*}
 (\tilde \phi_c^n + \tilde \phi_c^{n+1}) &((c^{n+1} - c_{\text{eq}}) - (c^n - c_{\text{eq}})) (c^{n+1} - c_{\text{eq}}) \\
 &= \tilde \phi_c^{n+1} (c^{n+1} - c_{\text{eq}})^2 - \tilde \phi_c^{n} (c^n - c_{\text{eq}})^2 + \tilde \phi_c^{n} (c^{n+1} - c^{n})^2 - (\tilde \phi_c^{n+1} - \tilde \phi_c^{n}) (c^{n} - c_{\text{eq}}) (c^{n+1} - c_{\text{eq}})
\end{align*}
and hence
\begin{align*}
 &\int_\Omega \frac{k_c}{c_{\text{eq}}^2} \frac{\tilde \phi_c^n + \tilde \phi_c^{n+1}}{2} \frac{(c^{n+1} - c_{\text{eq}}) - (c^n - c_{\text{eq}})}{\tau} (c^{n+1} - c_{\text{eq}}) \;d\vecx \\
 &\quad= \int_\Omega \frac{1}{\tau} \frac{1}{2} \frac{k_c}{c_{\text{eq}}^2} \tilde \phi_c^{n+1} (c^{n+1} - c_{\text{eq}})^2 \;d\vecx - \int_\Omega \frac{1}{\tau} \frac{1}{2} \frac{k_c}{c_{\text{eq}}^2} \tilde \phi_c^n (c^{n} - c_{\text{eq}})^2 \;d\vecx + \int_\Omega \frac{1}{\tau} \frac{1}{2} \frac{k_c}{c_{\text{eq}}^2} \tilde \phi_c^{n} (c^{n+1} - c^n)^2 \;d\vecx \\
 &\qquad- \int_\Omega \frac{1}{2} \frac{k_c}{c_{\text{eq}}^2} \frac{\tilde \phi_c^{n+1} - \tilde \phi_c^n}{\tau} (c^{n} - c_{\text{eq}}) (c^{n+1} - c_{\text{eq}}) \;d\vecx.
\end{align*}
We thus obtain
\begin{align}
\begin{split}
 0 &= \int_\Omega \frac{1}{\tau} \frac{1}{2} \frac{k_c}{c_{\text{eq}}^2} \tilde \phi_c^{n+1} (c^{n+1} - c_{\text{eq}})^2 \;d\vecx - \int_\Omega \frac{1}{\tau} \frac{1}{2} \frac{k_c}{c_{\text{eq}}^2} \tilde \phi_c^n (c^{n} - c_{\text{eq}})^2 \;d\vecx + \int_\Omega \frac{1}{\tau} \frac{1}{2} \frac{k_c}{c_{\text{eq}}^2} \tilde \phi_c^{n} (c^{n+1} - c^n)^2 \;d\vecx \\
 &\qquad- \int_\Omega \frac{1}{2} \frac{k_c}{c_{\text{eq}}^2} \frac{\tilde \phi_c^{n+1} - \tilde \phi_c^n}{\tau} (c^{n} - c_{\text{eq}}) (c^{n+1} - c_{\text{eq}}) \;d\vecx \\
 &\qquad+ \frac{1}{2} \int_\Omega \frac{k_c}{c_{\text{eq}}^2} \left(\phi_c^{n} \mathbf{v}^{n+1} + \mathbf{J}_c^{n+1}\right) \cdot \nabla \revision{P_h}\!\left((c^{n} - c_{\text{eq}}) (c^{n+1} - c_{\text{eq}})\right) d\vecx + \int_\Omega D \frac{k_c}{c_{\text{eq}}^2} (\tilde \phi_c^n \nabla c^{n+1}) \cdot \nabla c^{n+1} \;d\vecx \\
 &\qquad- \int_\Omega \frac{k_c}{c_{\text{eq}}^2} (c^\ast - c^n) R^{n+1} (c^{n+1} - c_{\text{eq}}) \;d\vecx \revision{- \int_\Omega \frac{1}{2} \frac{k_c}{c_{\text{eq}}^2} R^{n+1} (I - P_h) ((c^n - c_{\text{eq}}) (c^{n+1} - c_{\text{eq}})) \;d\vecx}.
\end{split}
\label{eq:femmodel_3_tested}
\end{align}
Note that by definition, see \eqref{eq:defrhof_2}, we have $\tilde \phi_c^{n+1} - \tilde \phi_c^{n} = \phi^{n+1} + \delta - (\phi^{n} + \delta) = \phi^{n+1} - \phi^{n}$ and therefore
\begin{align*}
 \int_\Omega \frac{1}{2} \frac{k_c}{c_{\text{eq}}^2} \frac{\tilde \phi_c^{n+1} - \tilde \phi_c^n}{\tau} (c^{n} - c_{\text{eq}}) (c^{n+1} - c_{\text{eq}}) \;d\vecx = \int_\Omega \frac{1}{2} \frac{k_c}{c_{\text{eq}}^2} \frac{\phi^{n+1} - \phi^n}{\tau} (c^{n} - c_{\text{eq}}) (c^{n+1} - c_{\text{eq}}) \;d\vecx.
\end{align*}
Using $\check \phi = \frac{1}{2} \frac{k_c}{c_{\text{eq}}^2} \revision{P_h} ((c^{n} - c_{\text{eq}}) (c^{n+1} - c_{\text{eq}}))$ in  $\eqref{eq:disc_weak_form}_4$ \revision{and the orthogonality of $P_h$} yields
\begin{align}
 0 &=
 \int_\Omega \frac{\phi^{n+1}-\phi^n}{\tau} \frac{1}{2} \frac{k_c}{c_{\text{eq}}^2} ((c^{n} - c_{\text{eq}}) (c^{n+1} - c_{\text{eq}})) \,d\vecx
 + \int_\Omega M \varepsilon \nabla \mu^{n+1} \cdot \nabla\!\left(\frac{1}{2} \frac{k_c}{c_{\text{eq}}^2} \revision{P_h} ((c^{n} - c_{\text{eq}}) (c^{n+1} - c_{\text{eq}}))\!\right) d\vecx \nonumber \\
 &\quad- \int_\Omega \left(\phi^{n} \mathbf{v}^{n+1}\right) \cdot \nabla\!\left(\frac{1}{2} \frac{k_c}{c_{\text{eq}}^2} \revision{P_h} ((c^{n} - c_{\text{eq}}) (c^{n+1} - c_{\text{eq}}))\!\right) d\vecx - \int_\Omega \frac{1}{2} \frac{k_c}{c_{\text{eq}}^2} R^{n+1} \revision{P_h} ((c^{n} - c_{\text{eq}}) (c^{n+1} - c_{\text{eq}})) \;d\vecx \nonumber \\
 &= \int_\Omega \frac{1}{2} \frac{k_c}{c_{\text{eq}}^2} \frac{\tilde \phi_c^{n+1} - \tilde \phi_c^n}{\tau} (c^{n} - c_{\text{eq}}) (c^{n+1} - c_{\text{eq}}) \;d\vecx - \frac{1}{2} \int_\Omega \frac{k_c}{c_{\text{eq}}^2} \left(\phi_c^{n} \mathbf{v}^{n+1} + \mathbf{J}_c^{n+1}\right) \cdot \nabla \revision{P_h}\!\left((c^{n} - c_{\text{eq}}) (c^{n+1} - c_{\text{eq}})\right) d\vecx \nonumber \\
 &\qquad- \int_\Omega \frac{1}{2} \frac{k_c}{c_{\text{eq}}^2} R^{n+1} P_h ((c^n - c_{\text{eq}}) (c^{n+1} - c_{\text{eq}})) \;d\vecx,
\label{eq:femmodel_4_tested_c}
\end{align}
as $\phi_c^n = \phi^n$ and $\mathbf{J}_c^{n+1} = - M \varepsilon \nabla \mu^{n+1}$. Adding \eqref{eq:femmodel_3_tested} and \eqref{eq:femmodel_4_tested_c}, and inserting definition \eqref{eq:react_disc} we have
\begin{align}
\begin{split}
 0 
 &= \int_\Omega \frac{1}{\tau} \frac{1}{2} \frac{k_c}{c_{\text{eq}}^2} \tilde \phi_c^{n+1} (c^{n+1} - c_{\text{eq}})^2 \;d\vecx - \int_\Omega \frac{1}{\tau} \frac{1}{2} \frac{k_c}{c_{\text{eq}}^2} \tilde \phi_c^n (c^{n} - c_{\text{eq}})^2 \;d\vecx + \int_\Omega \frac{1}{\tau} \frac{1}{2} \frac{k_c}{c_{\text{eq}}^2} \tilde \phi_c^{n} (c^{n+1} - c^n)^2 \;d\vecx \\
 &\qquad+ \int_\Omega D \frac{k_c}{c_{\text{eq}}^2} \tilde \phi_c^n \left|\nabla c^{n+1}\right|^2 \;d\vecx - \int_\Omega r(c^{n,n+1}) R^{n+1} \;d\vecx.
\end{split}
\label{eq:femmodel_3_tested_added}
\end{align}
}

Finally, adding \eqref{eq:femmodel_4_tested_mu}, \eqref{eq:femmodel_5_tested}, \eqref{eq:femmodel_2_tested_added} and $\frac{1}{\tilde \alpha}\cdot$\eqref{eq:femmodel_3_tested_added}, and using that for the surface tension term we have
\begin{align*}
 \int_\Omega \mathbf{\tilde S}^{n+1} \cdot \mathbf{v}^{n+1} \;d\vecx = - \int_\Omega \sigma (\phi^n \mathbf{v}^{n+1}) \cdot \nabla \mu^{n+1} \;d\vecx,
\end{align*}
we obtain
\begin{align*}
 0 &\geq \sigma \int_\Omega \frac{1}{\tau} \frac{W_\text{dw}(\phi^{n+1}) - W_\text{dw}(\phi^n)}{\varepsilon} \;d\vecx \\
 &\quad+ \sigma \left(\int_\Omega \frac{1}{\tau} \frac{\varepsilon}{2} |\nabla \phi^{n+1}|^2 \;d\vecx - \int_\Omega \frac{1}{\tau} \frac{\varepsilon}{2} |\nabla \phi^{n}|^2 \;d\vecx + \int_\Omega \frac{1}{\tau} \frac{\varepsilon}{2} |\nabla \phi^{n+1} - \nabla \phi^n|^2 \;d\vecx\right) \\
 &\quad+ \int_\Omega \frac{1}{\tau} \frac{1}{2} \tilde \rho_f^{n+1} |\mathbf{v}^{n+1}|^2 \;d\vecx - \int_\Omega \frac{1}{\tau} \frac{1}{2} \tilde \rho_f^{n} |\mathbf{v}^{n}|^2 \;d\vecx + \int_\Omega \frac{1}{\tau} \frac{1}{2} \tilde \rho_f^{n} |\mathbf{v}^{n+1} - \mathbf{v}^n|^2 \;d\vecx \\
 &\quad+ \int_\Omega \frac{1}{\tau} \frac{1}{2 \tilde \alpha} \frac{k_c}{c_{\text{eq}}^2} \tilde \phi_c^{n+1} (c^{n+1} - c_{\text{eq}})^2 \;d\vecx - \int_\Omega \frac{1}{\tau} \frac{1}{2 \tilde \alpha} \frac{k_c}{c_{\text{eq}}^2} \tilde \phi_c^n (c^{n} - c_{\text{eq}})^2 \;d\vecx + \int_\Omega \frac{1}{\tau} \frac{1}{2 \tilde \alpha} \frac{k_c}{c_{\text{eq}}^2} \tilde \phi_c^{n} (c^{n+1} - c^n)^2 \;d\vecx \\
 &\quad+ \int_\Omega \sigma M \varepsilon \left|\nabla \mu^{n+1}\right|^2 \;d\vecx + \int_\Omega 2 \gamma \nabla^s \mathbf{v}^{n+1} \mathbin{:} \nabla \mathbf{v}^{n+1} \;d\vecx + \int_\Omega \rho \frac{1}{\varepsilon} d(\tilde \phi_f^{n}) \left|\mathbf{v}^{n+1}\right|^2 \;d\vecx \\
 &\quad+ \int_\Omega \frac{1}{\tilde \alpha} D \frac{k_c}{c_{\text{eq}}^2} \tilde \phi_c^n \left|\nabla c^{n+1}\right|^2 \;d\vecx - \int_\Omega \left(\frac{1}{\tilde \alpha} r(c^{n,n+1}) + \sigma \mu^{n+1}\right) R^{n+1} \;d\vecx \\
 &= \frac{1}{\tau} \left(F[\phi^{n+1}, \nabla\phi^{n+1}, \mathbf{v}^{n+1}, c^{n+1}] -  F[\phi^{n}, \nabla\phi^{n}, \mathbf{v}^{n}, c^n]\right) \\
 &\quad+ \int_\Omega \frac{1}{\tau} \sigma \frac{\varepsilon}{2} |\nabla \phi^{n+1} - \nabla \phi^n|^2 \;d\vecx + \int_\Omega \frac{1}{\tau} \frac{1}{2} \tilde \rho_f^{n} |\mathbf{v}^{n+1} - \mathbf{v}^n|^2 \;d\vecx + \int_\Omega \frac{1}{\tau} \frac{1}{2 \tilde \alpha} \frac{k_c}{c_{\text{eq}}^2} \tilde \phi_c^{n} (c^{n+1} - c^n)^2 \;d\vecx \\
 &\quad+ \int_\Omega \sigma M \varepsilon \left|\nabla \mu^{n+1}\right|^2 \;d\vecx + \int_\Omega 2 \gamma \nabla^s \mathbf{v}^{n+1} \mathbin{:} \nabla \mathbf{v}^{n+1} \;d\vecx + \int_\Omega \rho \frac{1}{\varepsilon} d(\tilde \phi_f^{n}) \left|\mathbf{v}^{n+1}\right|^2 \;d\vecx \\
 &\quad+ \int_\Omega \frac{1}{\tilde \alpha} D \frac{k_c}{c_{\text{eq}}^2} \tilde \phi_c^n \left|\nabla c^{n+1}\right|^2 \;d\vecx - \int_\Omega \frac{1}{\tilde \alpha} \left(r(c^{n,n+1}) + \tilde \alpha \sigma \mu^{n+1}\right) R^{n+1} \;d\vecx,
\end{align*}
which by inserting definition \eqref{eq:R_disc} yields the discrete free energy dissipation inequality 
\begin{align*}
 &\frac{1}{\tau} \left(F[\phi^{n+1}, \nabla\phi^{n+1}, \mathbf{v}^{n+1}, c^{n+1}] -  F[\phi^{n}, \nabla\phi^{n}, \mathbf{v}^{n}, c^n]\right) \\
 &\quad\leq - \left(\int_\Omega \sigma M \varepsilon \left|\nabla \mu^{n+1}\right|^2 \;d\vecx + \int_\Omega 2 \gamma \nabla^s \mathbf{v}^{n+1} \mathbin{:} \nabla \mathbf{v}^{n+1} \;d\vecx + \int_\Omega \rho \frac{1}{\varepsilon} d(\tilde \phi_f^{n}) \left|\mathbf{v}^{n+1}\right|^2 \;d\vecx\right. \\
 &\qquad \qquad\left.+ \int_\Omega \frac{1}{\tilde \alpha} D \frac{k_c}{c_{\text{eq}}^2} \tilde \phi_c^n \left|\nabla c^{n+1}\right|^2 \;d\vecx + \int_\Omega \frac{q(\phi^n)}{\varepsilon \tilde \alpha} \left(r(c^{n,n+1}) + \tilde \alpha \sigma \mu^{n+1}\right)^2 \;d\vecx\right) \\
 &\quad \leq 0
\end{align*}
for all given $\tau$ and $n$.
\end{proof}

\revision{
\subsection{Projection-free scheme} \label{sec:proj_free}
Following the ansatz in \cite{Gruen2014} we can eliminate the projection operator $P_h$ from the discrete scheme \eqref{eq:disc_weak_form}. To do so, we insert $\check \phi = \frac{1}{2} \rho P_h (\mathbf{v}^{n} \cdot \mathbf{\check v})$ in $\eqref{eq:disc_weak_form}_4$, yielding
\begin{align*}
 0 &=
 \int_\Omega \frac{\phi^{n+1}-\phi^n}{\tau} \frac{1}{2} \rho \mathbf{v}^{n} \cdot \mathbf{\check v} \;d\vecx
 + \int_\Omega M \varepsilon \nabla \mu^{n+1} \cdot \nabla\!\left(\frac{1}{2} \rho P_h (\mathbf{v}^{n} \cdot \mathbf{\check v})\right) d\vecx \nonumber \\
 &\qquad- \int_\Omega \left(\phi^{n} \mathbf{v}^{n+1}\right) \cdot \nabla\!\left(\frac{1}{2} \rho P_h (\mathbf{v}^{n} \cdot \mathbf{\check v})\right) d\vecx \nonumber - \int_\Omega \frac{1}{2} \rho R^{n+1} P_h(\mathbf{v}^{n} \cdot \mathbf{\check v}) \;d\vecx \\
 &= \int_\Omega \frac{1}{2} \frac{\tilde \rho_f^{n+1} - \tilde \rho_f^{n}}{\tau} \mathbf{v}^{n} \cdot \mathbf{\check v} \;d\vecx - \frac{1}{2} \int_\Omega \left(\rho_f^{n} \mathbf{v}^{n+1} + \mathbf{J}_f^{n+1}\right) \cdot \nabla P_h\!\left(\mathbf{v}^{n} \cdot \mathbf{\check v}\right) d\vecx - \int_\Omega \frac{1}{2} \rho R^{n+1} P_h(\mathbf{v}^{n} \cdot \mathbf{\check v}) \;d\vecx. 
\end{align*}
We can now rewrite $\eqref{eq:disc_weak_form}_2$ in a projection-free form. Using
\begin{align*}
 \int_\Omega \frac{\tilde \rho_f^{n} + \tilde \rho_f^{n+1}}{2} \frac{\mathbf{v}^{n+1} - \mathbf{v}^n}{\tau} \cdot \mathbf{\check v} \;d\vecx + \int_\Omega \frac{1}{2} \frac{\tilde \rho_f^{n+1} - \tilde \rho_f^{n}}{\tau} \mathbf{v}^{n} \cdot \mathbf{\check v} \;d\vecx = \int_\Omega \frac{\tilde \rho_f^{n} + \tilde \rho_f^{n+1}}{2} \frac{\mathbf{v}^{n+1}}{\tau} \cdot \mathbf{\check v} \;d\vecx - \int_\Omega \frac{\tilde \rho_f^{n}}{\tau} \mathbf{v}^{n} \cdot \mathbf{\check v} \;d\vecx
\end{align*}
we obtain
\begin{align}
 &0 = \int_\Omega \frac{\tilde \rho_f^{n} + \tilde \rho_f^{n+1}}{2} \frac{\mathbf{v}^{n+1}}{\tau} \cdot \mathbf{\check v} \;d\vecx - \int_\Omega \frac{\tilde \rho_f^{n}}{\tau} \mathbf{v}^{n} \cdot \mathbf{\check v} \;d\vecx + \frac{1}{2} \int_\Omega \left(\rho_f^{n} \mathbf{v}^{n} + \mathbf{J}_f^{n}\right) \cdot \left(\left(\nabla \mathbf{v}^{n+1}\right) \mathbf{\check v}\right) d\vecx \nonumber \\
 &\qquad - \frac{1}{2} \int_\Omega \left(\rho_f^{n} \mathbf{v}^{n} + \mathbf{J}_f^{n}\right) \cdot \left(\left(\nabla \mathbf{\check v}\right) \mathbf{v}^{n+1}\right) d\vecx 
 - \int_\Omega p^{n+1} \nabla \cdot (\tilde \phi_f^{n+1} \mathbf{\check v}) \;d\vecx 
 + \int_\Omega 2 \gamma \nabla^s \mathbf{v}^{n+1} \mathbin{:} \nabla \mathbf{\check v} \;d\vecx \nonumber \\
 &\qquad + \int_\Omega \rho \frac{1}{\varepsilon} d(\tilde \phi_f^{n}) \mathbf{v}^{n+1} \cdot \mathbf{\check v} \;d\vecx - \int_\Omega \mathbf{\tilde S}^{n+1} \cdot \mathbf{\check v} \;d\vecx.
 \label{eq:disc_weak_form_2_proj_free}
\end{align}

Similarly, we insert $\check \phi = \frac{1}{2} P_h ((c^n - c_{\text{eq}}) \check c)$ in $\eqref{eq:disc_weak_form}_4$, yielding
\begin{align*}
 0 &= \int_\Omega \frac{1}{2} \frac{\tilde \phi_c^{n+1} - \tilde \phi_c^n}{\tau} (c^n - c_{\text{eq}}) \check c \;d\vecx - \frac{1}{2} \int_\Omega \left(\phi_c^n \mathbf{v}^{n+1} + \mathbf{J}_c^{n+1}\right) \cdot \nabla P_h\!\left((c^n - c_{\text{eq}}) \check c\right) d\vecx \\
 &\qquad- \int_\Omega \frac{1}{2} R^{n+1} P_h ((c^n - c_{\text{eq}}) \check c) \;d\vecx.
\end{align*}
We can now rewrite $\eqref{eq:disc_weak_form}_3$ in a projection-free form as
\begin{align*}
 0 &= \int_\Omega \frac{\tilde \phi_c^n + \tilde \phi_c^{n+1}}{2} \frac{c^{n+1} - c^n}{\tau} \check c \;d\vecx + \int_\Omega \frac{1}{2} \frac{\tilde \phi_c^{n+1} - \tilde \phi_c^n}{\tau} (c^n - c_{\text{eq}}) \check c \;d\vecx + \frac{1}{2} \int_\Omega \left((\phi_c^n \mathbf{v}^n + \mathbf{J}_c^n) \cdot \nabla c^{n+1}\right) \check c \;d\vecx \\
 &\qquad- \frac{1}{2} \int_\Omega \left((\phi_c^n \mathbf{v}^n + \mathbf{J}_c^n) \cdot \nabla \check c\right) (c^{n+1} - c_{\text{eq}}) \;d\vecx + \int_\Omega D (\tilde \phi_c^n \nabla c^{n+1}) \cdot \nabla \check c \;d\vecx \\
 &\qquad- \int_\Omega \frac{1}{2} R^{n+1}\!\left(2 (c^\ast - c^n) + (c^n - c_{\text{eq}})\right)\!\check c \;d\vecx.
\end{align*}
Adding a zero in the time evolution we further have
\begin{align*}
 \int_\Omega \frac{\tilde \phi_c^n + \tilde \phi_c^{n+1}}{2} \frac{c^{n+1} - c^n}{\tau} &\check c \;d\vecx + \int_\Omega \frac{1}{2} \frac{\tilde \phi_c^{n+1} - \tilde \phi_c^n}{\tau} (c^n - c_{\text{eq}}) \check c \;d\vecx \\
 &= \int_\Omega \frac{\tilde \phi_c^n + \tilde \phi_c^{n+1}}{2} \frac{c^{n+1} - c_{\text{eq}}}{\tau} \check c \;d\vecx - \int_\Omega \frac{\tilde \phi_c^n}{\tau} (c^n - c_{\text{eq}}) \check c \;d\vecx
\end{align*}
and hence obtain
\begin{align}
 0 &= \int_\Omega \frac{\tilde \phi_c^n + \tilde \phi_c^{n+1}}{2} \frac{c^{n+1} - c_{\text{eq}}}{\tau} \check c \;d\vecx - \int_\Omega \frac{\tilde \phi_c^n}{\tau} (c^n - c_{\text{eq}}) \check c \;d\vecx + \frac{1}{2} \int_\Omega \left((\phi_c^n \mathbf{v}^n + \mathbf{J}_c^n) \cdot \nabla c^{n+1}\right) \check c \;d\vecx \nonumber \\
 &\quad- \frac{1}{2} \int_\Omega \left((\phi_c^n \mathbf{v}^n + \mathbf{J}_c^n) \cdot \nabla \check c\right) (c^{n+1} - c_{\text{eq}}) \;d\vecx + \int_\Omega D (\tilde \phi_c^n \nabla c^{n+1}) \cdot \nabla \check c \;d\vecx - \int_\Omega \frac{1}{2} R^{n+1}\!\left(2 c^\ast - c^n - c_{\text{eq}}\right)\!\check c \;d\vecx.
\label{eq:disc_weak_form_3_proj_free}
\end{align}

For our numerical experiments in Section~\ref{sec:numer} we employ the projection-free discrete scheme and hence replace $\eqref{eq:disc_weak_form}_2$ by \eqref{eq:disc_weak_form_2_proj_free} and $\eqref{eq:disc_weak_form}_3$ by \eqref{eq:disc_weak_form_3_proj_free}.
}

\revision{
\section{Extended computational considerations} \label{sec:comput}
In this section we briefly discuss two computational features that may accelerate the overall simulation time, namely the choice of the solution strategy (monolithic or partitioned) and a preprocessing strategy to obtain an appropriate initial value $\phi_0$ for the phase-field variable. For the sake of simplicity  we consider both methods for the  non-reactive model~\eqref{eq:model_non_react} only.  
}

\subsection{Numerical solution strategies} \label{sec:solution}
In order to solve the \revision{projection-free} discretized system of equations \revision{derived from \eqref{eq:disc_weak_form} for the non-reactive case}, we implement two schemes following the monolithic and partitioned solution paradigms, respectively.

\paragraph{Monolithic case}
For the monolithic case we apply Newton's method in order to solve the coupled equations 
collectively resolving the complex nonlinearities of the problem. Employing Newton's method requires the solution of linear systems with the Jacobian  ${\mathcal A} \in \setR^{N_{\text{m}} \times N_{\text{m}}}$ as system matrix. For the reaction-free setting, it exhibits the general block structure
\begin{align}
\mathcal{A} 
= 
\begin{pmatrix}
A_{\textrm{CH}} & C_T \\
C_I & A_{\textrm{NS}}
\end{pmatrix}
\label{eq:lin_system}
\end{align}
where $A_{\textrm{CH}} \in \setR^{N_{\text{CH}} \times N_{\text{CH}}}$ denotes the Cahn-Hilliard block matrix, $A_{\textrm{NS}}\in \setR^{N_{\text{NS}} \times N_{\text{NS}}}$ the Navier-Stokes block matrix, and $C_T \in \setR^{N_{\text{CH}} \times N_{\text{NS}}}$ and $C_I \in \setR^{N_{\text{NS}} \times N_{\text{CH}}}$ account for the coupling through the transport at the interface and the coupling through the interfacial force respectively, and hence $N_{\text{m}} = N_{\text{NS}} + N_{\text{CH}}$. In our implementation, the unknowns are ordered as $(\phi,\mu,\mathbf{v},p)^T$ and the Cahn-Hilliard block refers to $\phi$ and $\mu$ while the Navier-Stokes block refers to $\mathbf{v}$ and $p$. The matrix dimensions depend on the Lagrange finite element spaces used as discrete ansatz spaces. For the specific choice of $i = k = 1$, $l = 2$ in \eqref{eq:FE_spaces}, that we employ in the numerical experiments in Section~\ref{sec:numer}, we obtain the relation $N_{\textrm{NS}} \approx 4.5 \cdot N_{\textrm{CH}}$.

This approach directly addresses the interrelated dynamics of the coupled equations by taking into account the off-diagonal blocks of the Jacobian in each Newton iteration.

\paragraph{Partitioned case}
For the partitioned case we decouple the system in the sense that we  alternate in solving the Navier-Stokes equations
with operator $A_{\textrm{NS}}$ and the Cahn-Hilliard equations
with operator $A_{\textrm{CH}}$ with updated data from the other in an iterative loop, similar to the Strang splitting \cite{Strang1968}. The resulting operator splitting has the advantage that problem-tailored solvers can be applied for the subsystems, rather than merely problem-tailored preconditioners to the diagonal blocks as in the monolithic case. Furthermore, for the fully coupled system there are also nonlinearities introduced in the coupling terms so that the splitting not only decouples but also partly linearizes the system. However, it is clear that both subproblems remain highly nonlinear so that in order to solve them Newton's method is applied respectively. 

This approach replaces the direct coupling by coupling iterations with smaller, potentially better conditioned Jacobians.

\medskip
In Section~\ref{sec:solution_strategies} we assess both solution strategies in a comparative study for different flow regimes, which shows that both approaches are suitable for the solution of the discretized system.

\subsection{Generation of the initial condition $\phi_0$} \label{sec:preprocessing}

For a non-trivial initial phase distribution it is often more convenient to parametrize the sharp-interface formulation than to directly parametrize the required smooth phase field $\phi_0$. Here, we provide a two-step strategy to obtain $\phi_0$ from arbitrary sharp-interface configurations.   
The formal workflow 
is as follows:
\begin{enumerate}
 \item
 Choose an open set $\Phi \subset \Omega \subset \setR^\nu$ representing the fluid phase. The sharp-interface configuration is then given by the characteristic function $\chi_\Phi$.
 \item Do a small number of iterations $n_{\text{pre}}$ of the Cahn-Hilliard operator $A_{\textrm{CH}}$ alone with $\chi_\Phi$ as initial condition and high mobility $M_{\text{pre}}$ in order to obtain a diffuse-interface configuration $\phi_0$.
\end{enumerate}
This strategy is \revision{straightforwardly applicable} in the partitioned setting as due to the underlying operator splitting the Cahn-Hilliard and Navier-Stokes operators are separately given, implemented and directly accessible.

Employing the strategy as a preprocessing step for the two-phase model enables us to calculate flow fields in complex geometries, see Section~\ref{sec:channel_flow}. To do so, the phase distribution should be almost static. We achieve this in practice by using a small mobility $M$ for the Cahn-Hilliard operator after the preprocessing step. 

Note that the switch of the mobility is key to obtain a valid strategy as the latter determines the time scale of the Cahn-Hilliard evolution. It is favorable to use a high preprocessing mobility $M_{\text{pre}}$ in order for the diffuse profile to develop rapidly within only a small number of iterations $n_{\text{pre}}$ of the Cahn-Hilliard operator. The parameters $n_{\text{pre}}$ and $M_{\text{pre}}$ are chosen due to numerical evidence. Conversely, for the actual flow calculations the time scale of interest is linked to the flow and hence the Cahn-Hilliard evolution should be comparably slow. Otherwise, this may, e.g., result in corners smoothing out over time further altering the complex geometry.

\section{Numerical experiments} \label{sec:numer}

The numerical experiments presented here are divided mainly into three studies\revision{, where the first two focus on assessing the non-reactive flow model \eqref{eq:model_non_react}}. Firstly, we investigate the thermodynamic consistency of the discrete numerical scheme, followed by a comparative study of the monolithic and partitioned strategies (cf.~Section~\ref{sec:solution}). To do so, we introduce a problem setup inspired by the classical fluid-flow problem of the lid-driven cavity extended to our two-phase model by adding circular solid inclusions in the flow. Furthermore, we want to assess the model's capability to capture flow scenarios in complex geometries following the strategy presented in Section~\ref{sec:preprocessing}. To do so, we introduce another problem setup that comprises a classical channel-flow setup extended by adding rectangular obstacles. Lastly, \revision{we assess the reactive model~\eqref{eq:model} that accounts for precipitation/dissolution effects at the liquid-solid interface. We numerically study the discrete thermodynamic consistency for the lid-driven cavity scenario, followed by the investigation of \revisionT{two channel-flow scenarios; one for the precipitation case that leads to clogging and one for the dissolution case.}} 
Before we present the simulations we collect some features 
of the numerical method.
\medskip

\textbf{Implementation details.}\label{sec:detail}
The implementation is carried out in DUNE \cite{Sander2020}, based on a previously published implementation \cite{DARUS2021} \revisionT{and made publicly available \cite{DARUS2026}}. The latter builds on the finite element framework PDELab \cite{Bastian2010}, which provides the usage of DUNE-ALUGrid \cite{Alkämper2016} as a grid interface. All computations are carried out single-threaded on an AMD Ryzen Threadripper 2950X 16-Core Processor. In the following, we highlight some algorithmic features and parameters.

\paragraph{Spatial discretization and grid adaptivity}\label{adap}
We use LBB-stable Taylor-Hood elements of degree $(2,1)$ for the pair $(\mathbf{v},p)$ and piecewise linear Lagrange elements for $c$, $\phi$ and $\mu$. Note that this corresponds to $i = j = k = 1$, $l = 2$ in \eqref{eq:FE_spaces}.

As the diffuse transition zone of width $\mathcal{O}(\varepsilon)$ is small yet needs to be resolved properly, a resolution of the entire grid with finite elements of same size is computationally infeasible. Therefore, if not stated otherwise, we use grid adaptivity \cite{Alkämper2016}, where the simplicial mesh is refined based on the gradient of the phase-field variable $\phi$ \cite{DARUS2021, phdvonWolff} leading to a fine resolution of the evolving transition layer.

\paragraph{Newton's method}
We apply Newton's method from \cite{Bastian2010}, where for the update step we employ the Hackbusch-Reusken line search \cite{Hackbusch1989}. 
Let $m$ denote the iteration index of Newton's method and $\mathbf{r}^m$ the corresponding residual vector. Newton's method is converged when the current residual reaches the stop tolerance $\text{maxTol}$ defined as  
\begin{align*}
 \text{maxTol} = \max\Big\{\left\|\mathbf{r}^0\right\|_2 \text{relTol}, \text{absTol}\Big\},
\end{align*}
with relative tolerance $\text{relTol}$ and absolute tolerance $\text{absTol}$. The linear solver tolerance $\text{linRed}$ is set to
\begin{align*}
 \text{linRed} = \max\left\{\frac{\text{maxTol}}{10 \left\|\mathbf{r}^m\right\|_2}, \min\left\{\text{minLinRed}, \frac{\left\|\mathbf{r}^m\right\|_2^2}{\left\|\mathbf{r}^{m-1}\right\|_2^2}\right\}\right\}, 
\end{align*}
where $\text{minLinRed}$ denotes the initial reduction parameter. We allow for a maximum number of $m_\text{N}$ iterations.
For the Hackbusch-Reusken line search we prescribe the maximum number of line search iterations $m_{\text{LS}}$ and the damping factor $d_{\text{LS}}$. 

Throughout all numerical experiments we set $m_\text{N} = 100$, $\text{relTol} = 10^{-7}$, $\text{absTol} = 10^{-9}$, $\text{minLinRed} = 10^{-3}$, $m_{\text{LS}} = 100$ and $d_{\text{LS}} = 0.5$. In case we solve the linear systems resulting from Newton's method directly, we use the direct solver UMFPack \cite{Davis2004}. 

\begin{figure}
 \centering
 \begin{subfigure}[b]{0.495\textwidth}
  \centering
  \includegraphics[keepaspectratio=true,width=\columnwidth]{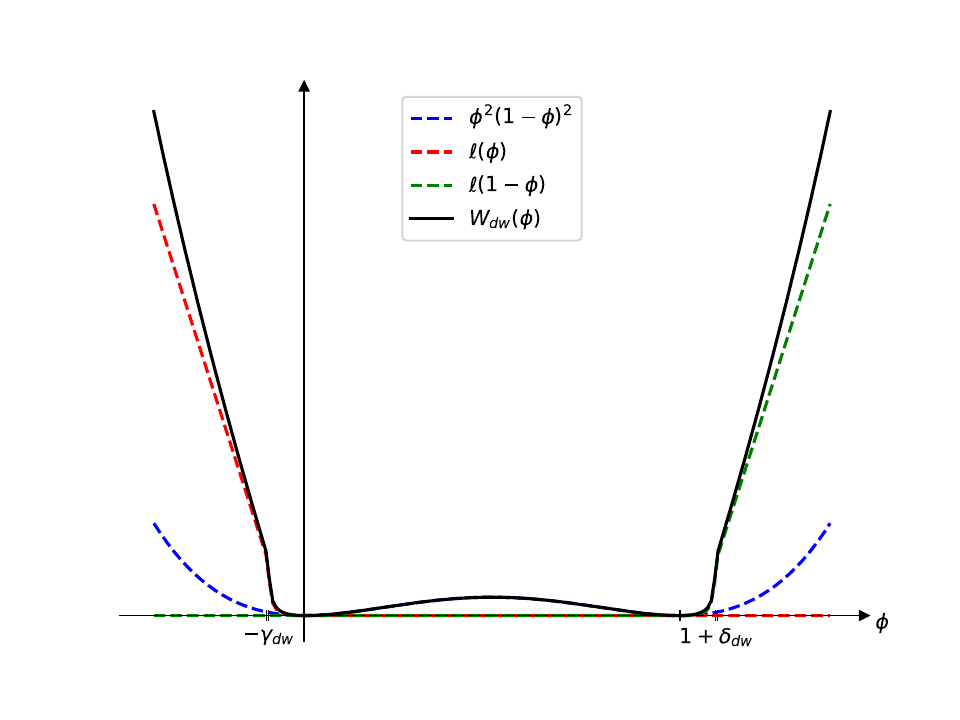}
 \end{subfigure}
 \begin{subfigure}[b]{0.495\textwidth}
  \centering
  \includegraphics[keepaspectratio=true,width=\columnwidth]{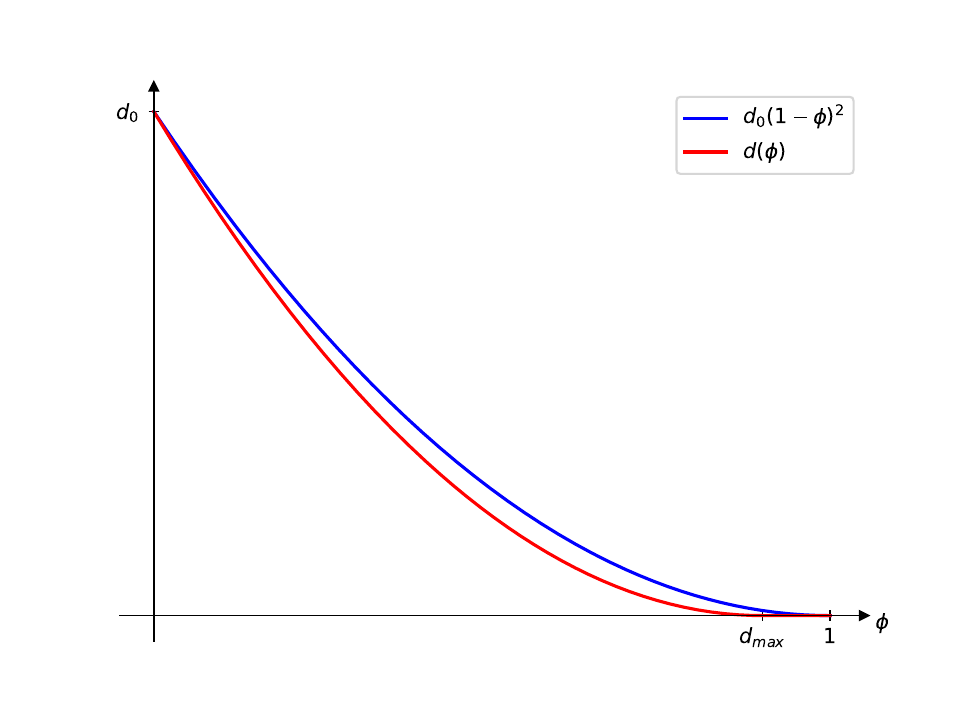}
 \end{subfigure}
 \caption{Schematic plots of the double-well potential (and its constituents) (left) and the drag terms (right). All axes are scaled linearly.}
 \label{fig:dw_dissipation}
\end{figure}

\paragraph{Double-well potential}
The double-well potential is defined by
\revision{
\begin{align*}
 W_{\text{dw}}(\phi) = 
 \left\{\!
 \begin{aligned}
  &450 \phi^4 (1 - \phi)^4 \\
  &\phi^2 (1 - \phi)^2 
 \end{aligned}
 \right\}
 + \ell(\phi) + \ell(1-\phi)
 \quad
 \begin{cases}
  \text{for the full model~\eqref{eq:model}} \\
  \text{for the non-reactive model~\eqref{eq:model_non_react}}
 \end{cases}
\end{align*}
}
with limiter function $\ell$. Here, we use
\begin{align}
 \ell(\phi) =
 \begin{cases}
  \delta_\text{dw} \left[\frac{\gamma_\text{dw}^2}{\delta_\text{dw} - \gamma_\text{dw}} + \left(\phi + \gamma_\text{dw}\right) \frac{-\gamma_\text{dw}\left(-\gamma_\text{dw} + 2\delta_\text{dw}\right)}{\left(\delta_\text{dw} - \gamma_\text{dw}\right)^2}\right] &\text{ for } \phi \leq -\gamma_\text{dw} \\
  \delta_\text{dw} \frac{\phi^2}{\phi + \delta_\text{dw}} &\text{ for } -\gamma_\text{dw} < \phi < 0 \\
  0 &\text{ for } \phi \geq 0
 \end{cases}
\label{eq:limiter}
\end{align}
with positive parameters $\delta_\text{dw} \leq \delta$ and $\gamma_\text{dw}$ satisfying $0 < \delta_\text{dw} - \gamma_\text{dw} \ll 1$, and impose the convex-concave split as in \cite{DARUS2021}. Note that, in contrast to what is stated in Remark~\ref{rem: model}\revision{(ii)}, the given limiter function $\ell(\phi)$ does not diverge at $\phi = -\delta$, and we apply this relaxation for stability reasons related to Newton's method. However, as long as the solution stays within $\phi > -\gamma_\text{dw}$, the solution is the same as for a limiter function that does not have the first case in \eqref{eq:limiter} and would diverge at $\phi = -\delta_\text{dw}$. \revisionT{In Theorem~\ref{theo:main} we assume that $\phi \in (- \delta, 1 + \delta)$ holds. In case one finds that $\phi$ exceeds these regularization induced bounds in a computation, one can a posteriori reset the limiter's parameters and redo the computation. A valid strategy is to set $\delta_{\text{dw}} = \delta$ and increase $\gamma_{\text{dw}}$ towards $\delta$, approximating a limiter that diverges at the bounds $- \delta$ and $1 + \delta$.}
\revision{We depict the double-well potential $W_{\text{dw}}(\phi)$ using the quartic term and its constituents on the left of Figure~\ref{fig:dw_dissipation}.}

\paragraph{Momentum dissipation in the solid phase} 
\revision{Following up Remark~\ref{rem: model}(iii), choosing a large momentum dissipation $d_0$ for the standard drag term $d(\tilde\phi_f) = d_0 \big(1 - \tilde\phi_f\big)^2$ may lead to undesired velocity dissipation in the fluid for $\tilde\phi_f$ close to one.} As remedy we utilize the modified drag term  
\begin{align*}
 d(\tilde\phi_f) = 
 \begin{cases}
  d_0 \frac{\left(d_{\text{max}} - \tilde\phi_f\right)^2}{d_{\text{max}}^2} &\text{ for } \tilde\phi_f \leq d_{\text{max}} \\
  0 &\text{ for } \tilde\phi_f > d_{\text{max}}
 \end{cases}
\end{align*}
with constant dissipation cut-off $0 < d_{\text{max}} \leq 1$.
\revision{We depict the two drag term variants presented on the right of Figure~\ref{fig:dw_dissipation}.}

\subsection{Lid-driven cavity inspired setup} \label{sec:lid-driven}
We modify the classical lid-driven cavity setup known from fluid dynamics by incorporating circular solid inclusions in the fluid flow. Note that using the strategy presented in Section~\ref{sec:preprocessing} is not needed here as a feasible (diffuse) initial profile $\phi_0$ for the phase-field variable can directly be given.
In order to do so, the corresponding parametrization describes a circle in terms of its center $\mathbf{x}_{\text{c}}$ and its radius $r$, resulting in a radially symmetric formula for the spherical initial conditions. Energy considerations (for the Ginzburg-Landau energy) for the simple setting of a planar interface separating two phases result in a closed form expression of the one-dimensional equilibrium profile given by
\begin{align*}
 \phi(x) = \frac{1}{2} \left(1 + \tanh\!\left(\frac{1}{\sqrt{2}\varepsilon}x\right)\right),
\end{align*}
see e.g.~\cite{Liu2003, Donaldson2011}.
Extended to the radially symmetric case, for one circular inclusion the initial phase field then takes the form 
\begin{align*}
 \phi(\mathbf{x},0) = \frac{1}{2} \left(1 + \tanh\Bigg(\frac{1}{\sqrt{2}\varepsilon}\Big(\|\mathbf{x} - \mathbf{x}_{\text{c}}\|_2 - r\Big)\Bigg)\right).
\end{align*}
In general, employing $N_s \in \mathbb{N}$ spheres in the flow field with parametrization
\begin{align*}
 \revision{\phi_i(\mathbf{x})} = \frac{1}{2} \left(1 + \tanh\!\left(\frac{1}{\sqrt{2}\varepsilon}\Big(\|\mathbf{x} - \mathbf{x}_{i,\text{c}}\|_2 - r_i\Big)\right)\right), \quad i \in \{1,\ldots,N_s\},
\end{align*}
the initial phase field takes the form
\begin{align*}
 \phi(\mathbf{x},0) = \prod_{i=1}^{N_s} \revision{\phi_i(\mathbf{x})}.
\end{align*}
For the velocity we utilize lid-driven cavity inspired Dirichlet boundary conditions, i.e., a parabolic velocity profile in horizontal direction on the top boundary of the domain $\Omega = [0,2] \times [0,1]$, 
\begin{align*}
 v_1\!\left((x_1,1)^T,t\right) = \bar{f} \frac{2}{3} \frac{x_1(2 - x_1)}{2}, \quad t \in [0,T],
\end{align*}
with mean flow $\bar{f}$, and value zero on the others. The velocity is initialized with value zero in the entire domain. For the phase field $\phi$, chemical potential $\mu$ and pressure $p$ we employ zero Neumann boundary conditions.

The phase-field model parameters are given in Table~\ref{tab:parameters_ldc}. Moreover, for the mean flow we set $\bar{f} = 0.1$.

\begin{table}
\centering
\caption{Phase-field model parameters for the lid-driven cavity inspired setup. The given dimensions of the parameters are consistent with three spatial dimensions, where $\mathrm{T}$, $\mathrm{L}$ and $\mathrm{M}$ denote the dimension symbols of time, length and mass respectively.}
\begin{tabular}{llll}
 \hline
 Parameter & Symbol & Value & Dimension \\ 
 \hline
 Fluid density & $\rho$ & $1.0$ & $\mathrm L^{-3} \mathrm M$ \\ 
 Fluid viscosity & $\gamma$ & $0.01$ & $\mathrm T^{-1} \mathrm L^{-1} \mathrm M$ \\ 
 Momentum dissipation in solid phase & $d_0$ & $1000$ & $\mathrm T^{-1}$ \\ 
 Dissipation cut-off & $d_{\text{max}}$ & $0.9$ & -- \\
 Phase-field mobility & $M$ & $1.0$ & $\mathrm T^{-1} \mathrm L^2$ \\ 
 Surface tension coefficient & $\sigma$ & $1.0$ & $\mathrm T^{-2} \mathrm M$ \\ 
 Diffuse-interface width & $\varepsilon$ & $0.03$ & $\mathrm L$ \\
 Phase-field regularization & $\delta$  & $0.03$ & -- \\
 Double-well potential modification & $\delta_{\text{dw}}$ & $0.02$ & -- \\
 Double-well potential modification & $\gamma_{\text{dw}}$ & $0.015$ & -- \\
 \hline
\end{tabular}
\label{tab:parameters_ldc}
\end{table}

\subsubsection{Thermodynamic consistency} \label{sec:num_thm}
Due to the potential energy increase caused by coarsening of the computational grid \cite{Garcke2016}, we do not enable grid adaptivity here during the simulation but only for the initial refinement step.

\begin{figure}
 \centering
 \begin{subfigure}[b]{0.495\textwidth}
  \centering
  \includegraphics[keepaspectratio=true,width=\columnwidth]{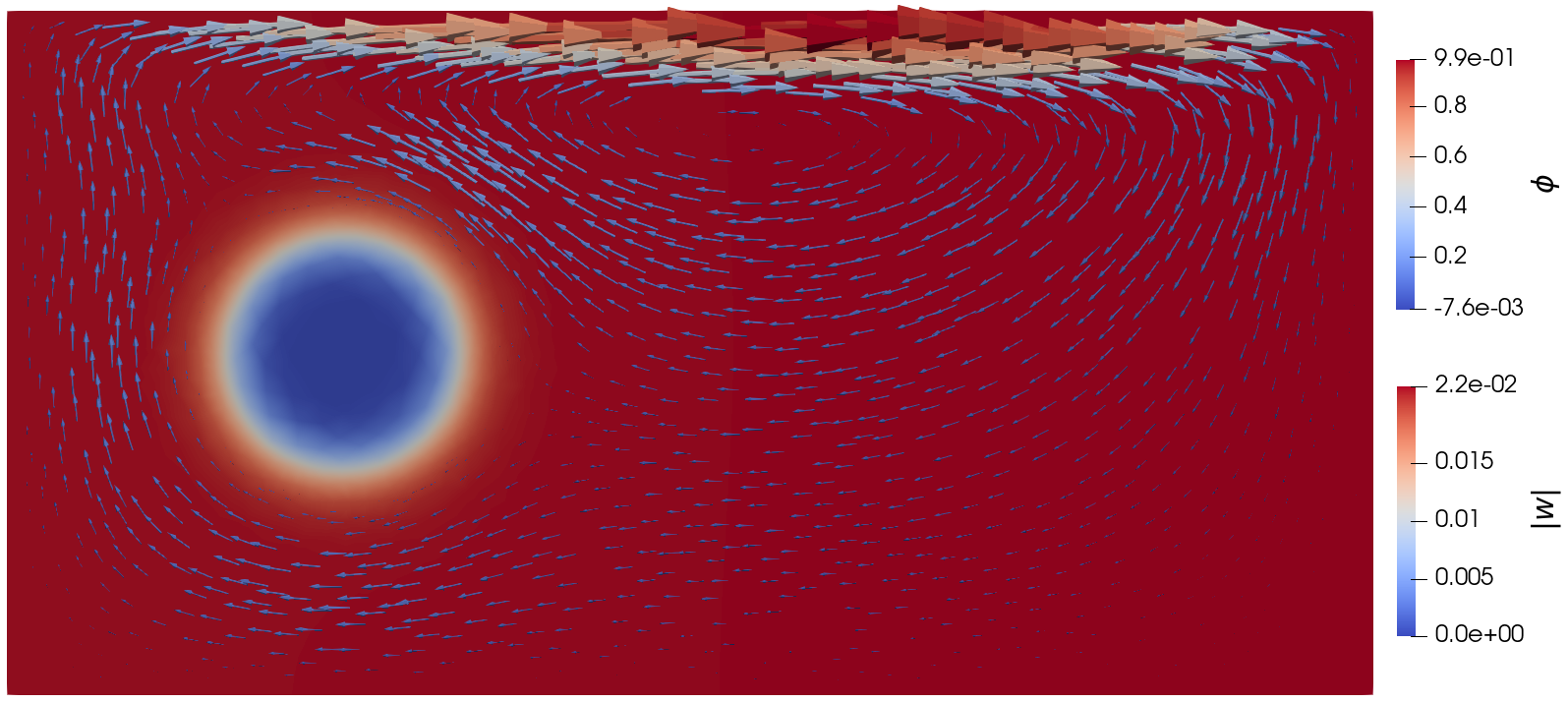}
 \end{subfigure}
 \begin{subfigure}[b]{0.495\textwidth}
  \centering
  \includegraphics[keepaspectratio=true,width=\columnwidth]{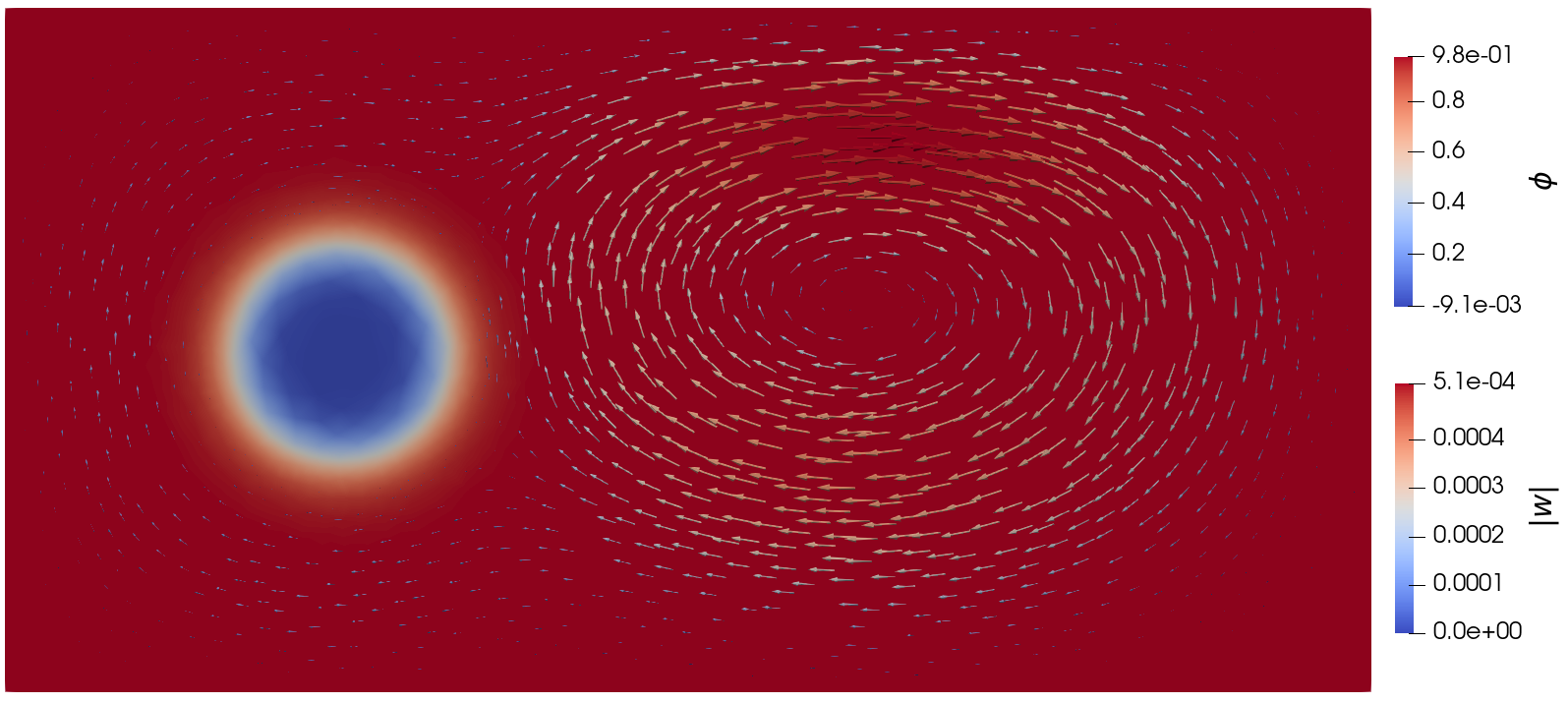}
 \end{subfigure}
 \caption{Lid-driven cavity inspired setup with one solid inclusion. On the left, we depict the result of the first time step after the driving flow is stopped. On the right, we show the result of the final time step. The fluid phase is colored red and the solid phase blue. The phase transition is illustrated by the transition in color. The arrows are scaled according to the magnitude of the conserved quantity $\mathbf{w} := \tilde \phi_f \mathbf{v}$. The color scaling is due to the data ranges presented on the right of the respective subfigure. Note that the scaling is two orders of magnitude smaller in the right subfigure.}
 \label{fig:model_problem_one_sphere}
\end{figure}

For the proof of thermodynamic consistency we have relied on the classical set of boundary conditions that comprises homogeneous Dirichlet conditions for the velocity $\mathbf{v}$. It is obvious that the latter is not satisfied for the lid-driven cavity inspired setup, where the flow boundary condition can be a source for an energy increase. Hence, in order to investigate the thermodynamic consistency for a non-trivial flow field, we use the following modification. In a preparatory step, we compute 30 time steps with the lid-driven cavity conditions. We then stop the driving flow to obtain a theoretically consistent initial flow field (and an equilibrated initial phase field) for the subsequent iterations of the two-phase model. The solid inclusion has center $\mathbf{x}_{c} = (0.5,0.5)^T$ and radius $r = 0.2$. The resulting setup is depicted on the left of Figure~\ref{fig:model_problem_one_sphere}. We calculate a total of $200$ time steps (comprising the preparatory ones) with time step size $\tau = 0.02$ resulting in final time $T = 4$ and the corresponding result is shown on the right of Figure~\ref{fig:model_problem_one_sphere}. \revisionU{We observe that, while the phase-field variable $\phi$ exceeds its physical bounds $[0,1]$, $\phi \in (- \delta, 1 + \delta)$ holds true and hence the $\delta$-regularized quantities $\tilde \phi_f$ \eqref{eq:defphif} and $\tilde \rho_f$ \eqref{eq:defrhof_2} do not attain negative values, preventing degeneration of the model and ensuring the validity of Theorem~\ref{theo:main}. In fact, the latter also holds for all numerical experiments that follow as we respectively obtain $\phi \in (- \delta, 1 + \delta)$.}

For the parameters given, we find that the circular solid inclusion shrinks to some extent, while an inclusion with initial radius $r_\text{v} = 0.175$ vanishes completely. These observations comply with findings of \cite{Yue2007} investigating spontaneous shrinkage of drops subject to Cahn-Hilliard diffusion. We note that the Cahn number for our setup is $Cn = \frac{\varepsilon}{r} = 0.15$. The volume of the domain $\Omega$ is $V = |\Omega| = 2$. The critical radius derived in \cite{Yue2007} then evaluates to
\begin{align*}
 r_\text{c} = \left(\frac{\sqrt{6}}{8 \pi} V \varepsilon\right)^{\!\frac{1}{3}} = \left(\frac{3\sqrt{6}}{400 \pi}\right)^{\!\frac{1}{3}} \approx 0.180,
\end{align*}
i.e., for an initial radius $r < r_\text{c}$ the drop is predicted to vanish. Lastly, the formula for the shrinking of the drop radius presented in \cite{Yue2007} evaluates to
\begin{align*}
 \tilde r = - \frac{\sqrt{2} V}{24 \pi} \frac{\varepsilon}{r^2} = - \frac{\sqrt{2}}{16 \pi} \approx - 0.028
\end{align*}
and hence one obtains the reduced radius $\revision{r + \tilde r} \approx 0.172$.

From these formulae, it is clear that there are measures to mitigate the shrinking process, such as decreasing $\varepsilon$, which may be taken given the implications on mass conservation discussed in \cite{Yue2007}. Other possible adjustments comprise employing an adapted double-well potential \cite{Donaldson2011} and reducing the mobility $M$ to increase the shrinkage time (compared to the time scale of interest) \cite{Yue2007}. We make use of the latter in our preprocessing strategy in Section~\ref{sec:preprocessing}.

\begin{figure}
 \centering
 \includegraphics[width=0.5\linewidth]{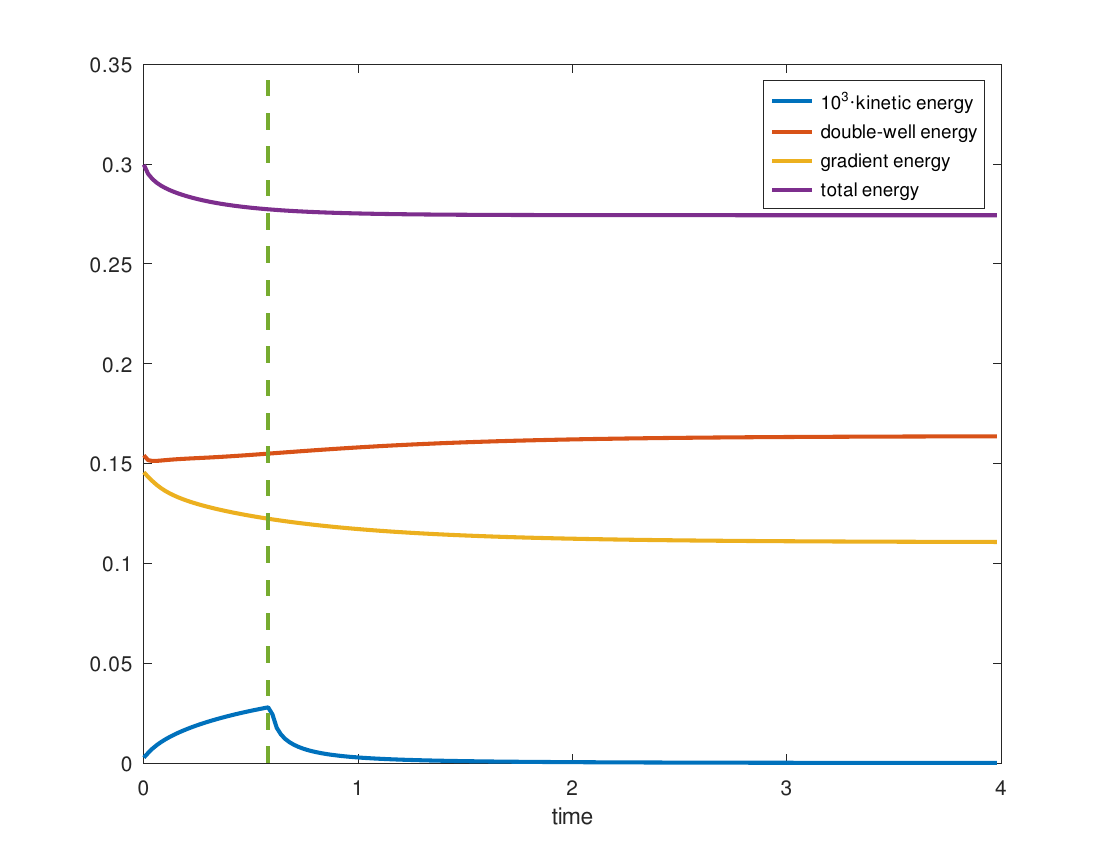}
 \caption{Energy plot for the lid-driven cavity inspired setup with one solid inclusion and stopped driving flow (dashed line).}
 \label{fig:model_problem_one_sphere_energy}
\end{figure}

In Figure~\ref{fig:model_problem_one_sphere_energy} we depict the time evolution of the \revision{non-reactive free energy $F[\phi(\cdot,t), \nabla\phi(\cdot,t), \mathbf{v}(\cdot,t)] $ derived from} \eqref{eq:free_energy} and the subenergies it is composed of, where we further split the Ginzburg-Landau energy into its two constituents. For the given setup, we observe that the Ginzburg-Landau energy for the Cahn-Hilliard evolution is orders of magnitude larger than the kinetic energy. Therefore, we scale the latter by a factor of $10^3$ to properly show its curve in the plot. The dashed line marks the time step where we stop the driving flow which is also reflected in the kinetic energy curve. The total energy is decreasing in accordance with Theorem~\ref{theo:main}.

\subsubsection{Solution strategies} \label{sec:solution_strategies}
This section is devoted to assessing the monolithic and partitioned solution strategies of Section~\ref{sec:solution}. We briefly discuss the blocks of the respective Jacobian matrices along with some matrix properties. Moreover, we introduce a feasible preconditioning technique for the partitioned case and conduct a comparative study with a focus on efficiency (in terms of runtime) and robustness (in terms of parameter variations). 

For the comparative study, we establish three different solver setups; monolithic with direct solver, partitioned with both subsystems solved directly, and partitioned with tailor-made preconditioned iterative solution. We again highlight that the underlying operator split in the partitioned setting leads to the iterative solution of the respective subsystems with matrices $A_{\textrm{CH}}$ and $A_{\textrm{NS}}$, allowing the use of optimized subsystem solvers. It is thus clear that the last option offers various combinations and optimization potential, however, we merely focus on one specific choice to showcase the feasibility of the approach. As linear solver we employ Flexible Restarted GMRES (FGMRES) \cite{Saad2003} for both subsystems, preconditioned with the ILU($0$) preconditioner available in dune-istl \cite{Blatt2007} for Navier-Stokes and a SIMPLE-type preconditioner \cite{Benzi2005} for Cahn-Hilliard. The former is applicable to saddle-point problems where fill-in is produced mostly within the block of the factors corresponding to the (zero) pressure block. A brief discussion of the use of incomplete factorizations as preconditioners for saddle-point matrices is presented in \cite{Benzi2005}. We use a restart parameter of $200$ for FGMRES with a maximum number of $1000$ iterations for the Navier-Stokes and $10000$ iterations for the Cahn-Hilliard solver.

We calculate $20$ time steps with time step size $\tau = 0.02$ of the lid-driven cavity inspired setup with one solid inclusion. As a spatial measure we state the number of unknowns which intrinsically corresponds to the size of the solution vector of the linear system arising from Newton's method. For the initially refined grid, we denote the latter as $N_{\text{m}}^0$ for the monolithic scheme, and for the partitioned scheme we write $N_{\textrm{NS}}^0$ for the Navier-Stokes subsystem and $N_{\textrm{CH}}^0$ for the Cahn-Hilliard subsystem respectively. The computations show that over time there is a slight increase in the number of unknowns due to the adaptive grid refinement. To show that $N_{\ast} \in \mathcal{O}(N_{\ast}^0)$ still holds, we further state the numbers after the final refinement step $N_{\ast}^{\text{f}}$ with $\ast \in \{\text{m}, \textrm{NS}, \textrm{CH}\}$. 

We investigate the dependence of the monolithic and the partitioned matrix blocks and schemes based on the viscosity $\gamma$. As performance metrics we measure the CPU time. Note that for the partitioned scheme we start the iterative procedure by solving the Navier-Stokes subsystem first followed by the Cahn-Hilliard subsystem. In Table~\ref{tab:cond_gamma} we show 1-norm condition number estimates, calculated by invoking the command \textit{condest} from GNU Octave, for the initial respective matrix blocks. Increasing the viscosity $\gamma$ leads to an upward trend in the estimated condition numbers for both the monolithic matrix $\mathcal{A}$ and the Navier-Stokes block $A_{\textrm{NS}}$, whereas the Cahn-Hilliard block $A_{\textrm{CH}}$ does not depend on $\gamma$.

\begin{table}
\centering
\caption{1-norm condition number estimates $\tilde{\kappa}$ for varying viscosity $\gamma$.}
\begin{tabular}{l|llllll}
 \hline
 $\gamma$ & $10^{-4}$ & $10^{-3}$ & $10^{-2}$ & $10^{-1}$ & $1$ & $10$ \\
 \hline
 $\tilde{\kappa}_\mathcal{A}$ & $8.8 \cdot 10^{19}$ & $4.0 \cdot 10^{20}$ & $3.9 \cdot 10^{20}$ & $4.2 \cdot 10^{21}$ & $5.2 \cdot 10^{20}$ & $3.1 \cdot 10^{23}$ \\
 $\tilde{\kappa}_{A_{\textrm{NS}}}$ & $1.5 \cdot 10^{17}$ & $3.5 \cdot 10^{14}$ & $7.5 \cdot 10^{16}$ & $7.5 \cdot 10^{19}$ & $1.8 \cdot 10^{19}$ & $3.1 \cdot 10^{20}$ \\
 $\tilde{\kappa}_{A_{\textrm{CH}}}$ & $1934.2$ & $1934.2$ & $1934.2$ & $1934.2$ & $1934.2$ & $1934.2$ \\
 \hline
\end{tabular}
\label{tab:cond_gamma}
\end{table}

We show the results in Table~\ref{tab:results_gamma}. We highlight that the results reveal an increase in robustness for the fully iterative partitioned solution compared to the monolithic and partitioned schemes using solely the
direct solver UMFPack. While these fail to converge for values $\gamma \geq 1$, the former converges for all viscosity values considered. Considering the range in which all schemes converge, we observe that the monolithic direct solver is the fastest. For the partitioned schemes, the iterative method is faster than its direct counterpart. The average number of coupling iterations per time step for the partitioned
schemes to converge is (about) four. This indicates that the coupling terms $C_T$ and $C_I$ from \eqref{eq:lin_system} are rather weak which benefits the operator splitting.  

The increase in robustness along with the high optimization potential related to the preconditioning strategies make the iterative partitioned scheme a preferred  approach for the solution of the two-phase model.

\revision{The extension and assessment of the partitioned solution strategy to the reactive model~\eqref{eq:model} is part of future work.}

\begin{table}
\centering
\caption{Results for varying viscosity $\gamma$. For the reported CPU times, the symbol $ \rm t(m)$ refers to the monolithic scheme, $ \rm t(d)$  to the partitioned scheme with direct solution and $ \rm t(i) $ to the partitioned scheme with iterative solution. $\Sigma$Ni denotes the total number of Newton iterations and $\Sigma$Ci the total number of coupling iterations. The cases where Newton's method did not converge are marked with "--", and the number in brackets gives the corresponding time step. For the numbers of unknowns we have $N_{\text{m}}^0 = 22310$, $N_{\text{m}}^{\text{f}} = 23828$, $N_{\text{CH}}^0 = 4112$, $N_{\text{CH}}^{\text{f}} = 4388$, $N_{\text{NS}}^0 = 18198$ and $N_{\text{NS}}^{\text{f}} = 19440$.}
\begin{tabular}{l|ll|llll|llll}
 \hline
 $\gamma$ & $\Sigma$Ni & t (m) 
 & $\Sigma$Ci & $\Sigma$Ni$_{\textrm{NS}}$ & $\Sigma$Ni$_{\textrm{CH}}$ & t (d) 
 & $\Sigma$Ci & $\Sigma$Ni$_{\textrm{NS}}$ & $\Sigma$Ni$_{\textrm{CH}}$ & t (i) 
 \\
 \hline
 $10^{-4}$ & $58$  & $46.56$ & $80$ & $87$  & $78$ & $73.18$  & $80$ & $100$ & $78$ & $62.55$  \\
 $10^{-3}$ & $58$  & $46.84$ & $80$ & $81$  & $78$ & $69.48$  & $80$ & $100$ & $78$ & $61.89$  \\
 $10^{-2}$ & $63$  & $50.29$ & $80$ & $101$ & $78$ & $83.19$  & $80$ & $100$ & $78$ & $61.59$  \\
 $10^{-1}$ & $84$  & $64.35$ & $80$ & $159$ & $78$ & $122.22$ & $80$ & $101$ & $78$ & $71.39$  \\
 $1$       & $(2)$ & --      & --   & $(2)$ & --   & --       & $80$ & $103$ & $81$ & $100.00$ \\
 $10$      & $(2)$ & --      & --   & $(1)$ & --   & --       & $80$ & $107$ & $83$ & $124.71$ \\
 \hline
\end{tabular}
\label{tab:results_gamma}
\end{table}

\subsubsection{Long-term model behavior}
\begin{figure}
 \centering
 \begin{subfigure}[b]{0.495\textwidth}
  \centering
  \includegraphics[keepaspectratio=true,width=\columnwidth]{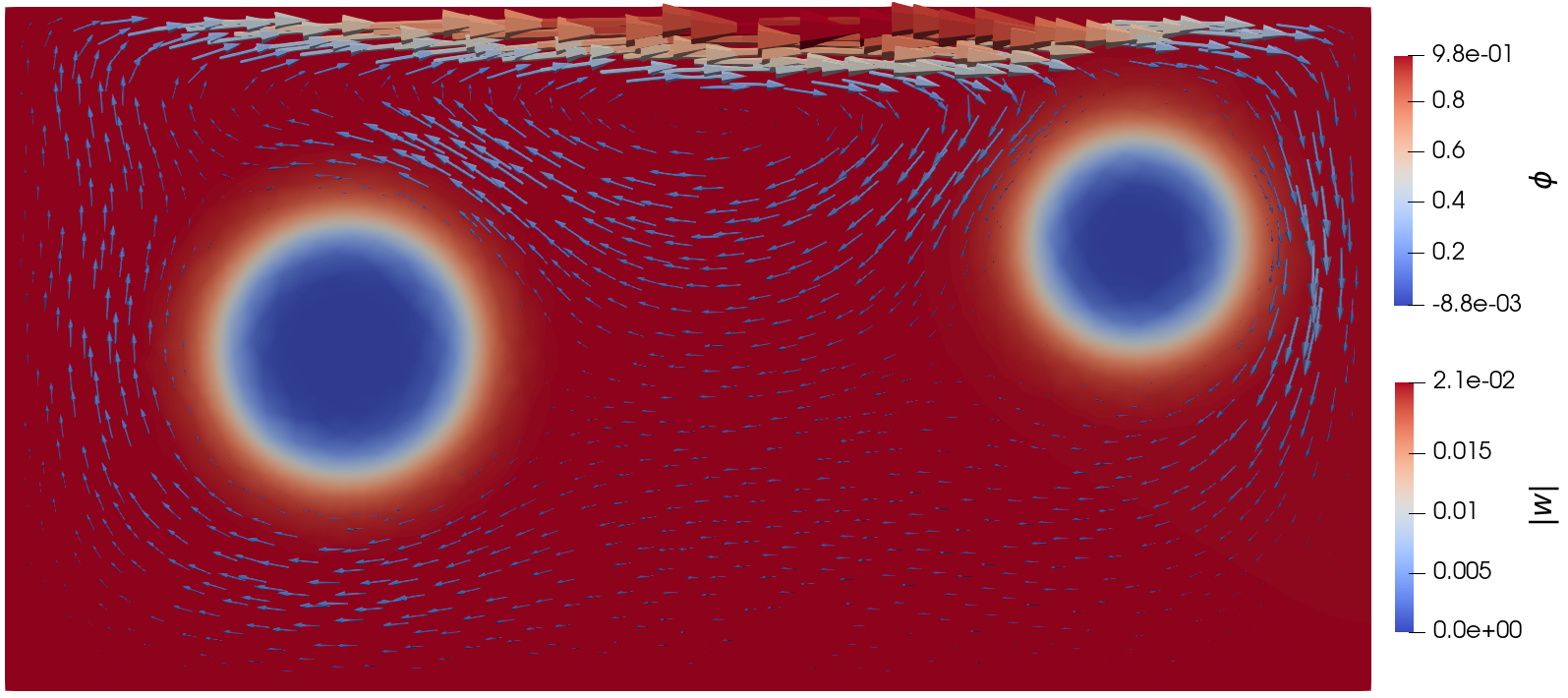}
 \end{subfigure}
 \begin{subfigure}[b]{0.495\textwidth}
  \centering
  \includegraphics[keepaspectratio=true,width=\columnwidth]{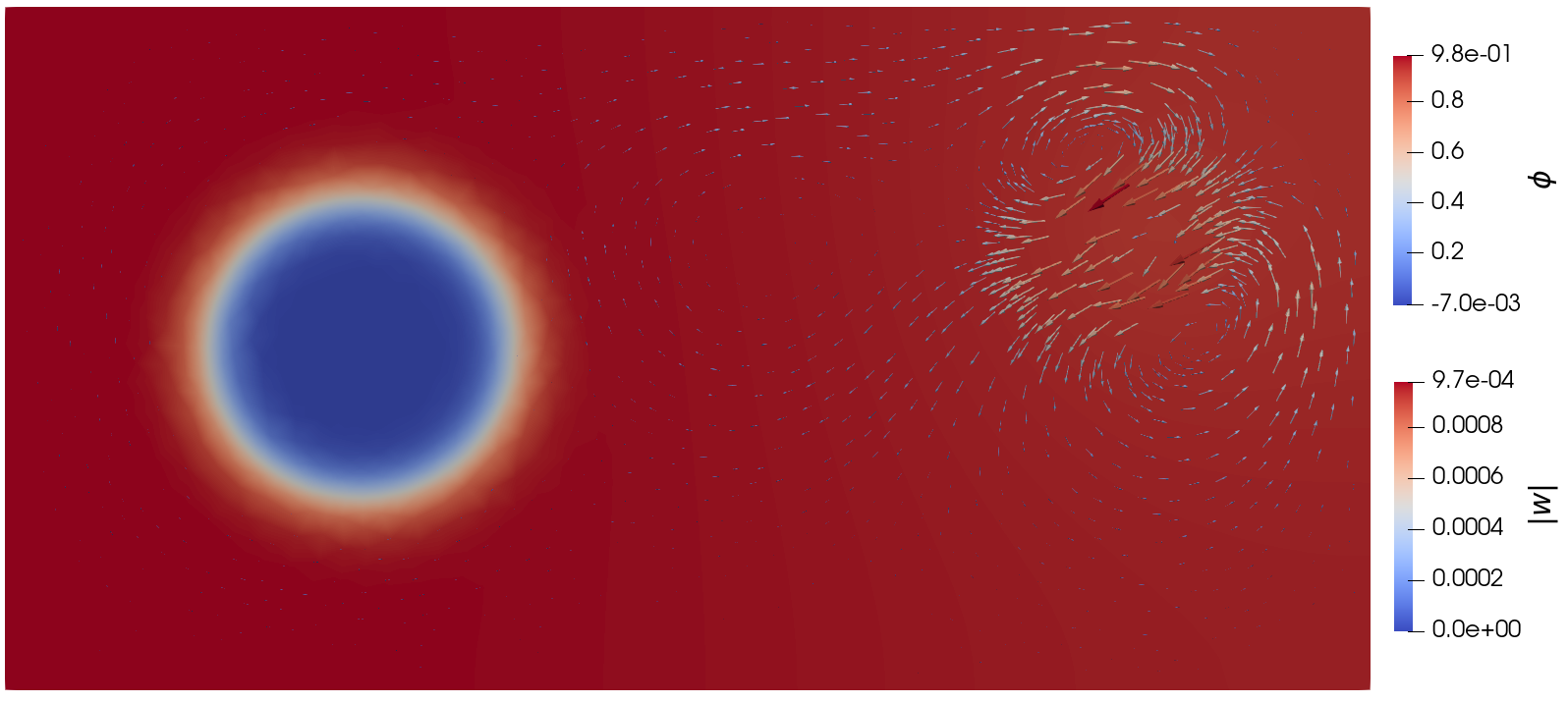}
 \end{subfigure}
 \caption{Lid-driven cavity inspired setup with two solid inclusions. On the left, we depict the result of the first time step after the driving flow is stopped. On the right, we show the result of the time step where the smaller inclusion is completely vanished. The arrows are scaled according to the magnitude of the conserved quantity $\mathbf{w} := \tilde \phi_f \mathbf{v}$. The color scaling is due to the data ranges presented on the right of the respective subfigure, which is an order of magnitude smaller in the right subfigure.}
 \label{fig:model_problem_two_spheres}
\end{figure}

In order to investigate the long-term model behavior (and the related effect of the Cahn-Hilliard evolution), we propose the following experiment. We still use the lid-driven cavity inspired setup, but now with two (non-overlapping) solid inclusions with different radii. We keep the parameters from before including the solid inclusion with center $\mathbf{x}_{1,c} = (0.5,0.5)^T$ and radius $r_1 = 0.2$. The second solid inclusion has its center in $\mathbf{x}_{2,c} = (1.65,0.65)^T$ and radius $r_2 = 0.175 < r_1$. To assess the thermodynamic consistency, we again stop the driving flow after 30 time steps. The resulting setup is depicted on the left of Figure~\ref{fig:model_problem_two_spheres}.

\begin{figure}
 \centering
 \includegraphics[width=0.5\linewidth]{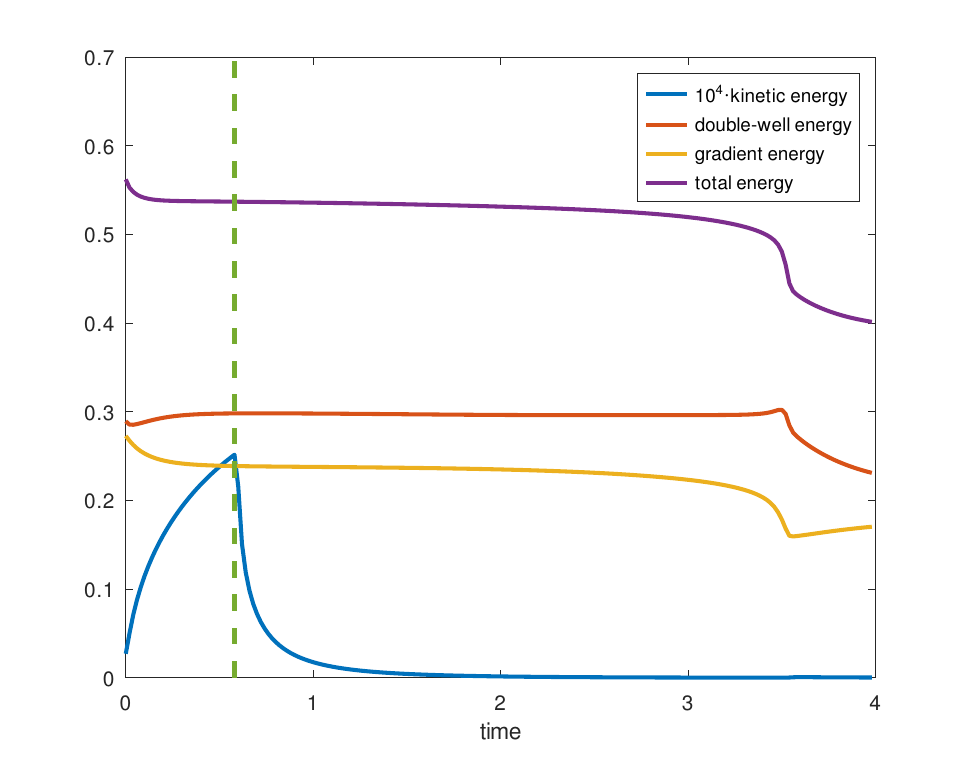}
 \caption{Energy plot for the lid-driven cavity inspired setup with two solid inclusions and stopped driving flow (dashed line).}
 \label{fig:model_problem_two_spheres_energy}
\end{figure}

We observe that the smaller inclusion is assimilated by the larger one as shown on the right of Figure~\ref{fig:model_problem_two_spheres}, which is a well-known phenomenon of the Cahn-Hilliard evolution \cite{Liu2003}. We now check the validity of the energy estimate in this rather involved setting. In Figure~\ref{fig:model_problem_two_spheres_energy} we depict the time evolution of the \revision{reduced} free 
energy  expression
$F[\phi(\cdot,t), \nabla\phi(\cdot,t), \mathbf{v}(\cdot,t)]$ \revision{derived} from  \eqref{eq:free_energy} and its subenergies. For better illustration, we scale the kinetic energy by a factor of $10^4$, whose curve reflects the time step where the driving flow is stopped, marked by the dashed line. The total energy is decreasing with a clearly observable drop at the  time step when the smaller inclusion vanishes. This is due to the expected drop in the gradient energy in the instance of assimilation that is not compensated by the concomitant slight increase in the double-well and the kinetic energy.

\subsection{Channel-flow inspired setup} \label{sec:channel_flow}
\begin{figure}
 \centering
 \includegraphics[width=0.5\linewidth]{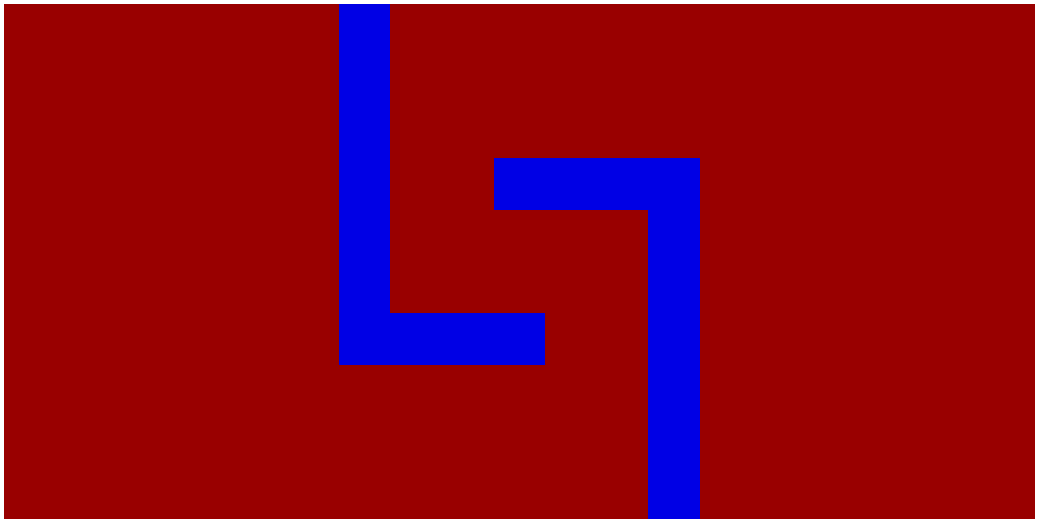}
\caption{Sharp-interface configuration of the channel-flow inspired setup.}
\label{fig:model_problem_presmoothing_illu}
\end{figure}

We modify the classical channel-flow setup known from fluid dynamics by incorporating (overlapping) rectangular obstacles in the fluid flow with the governing scenario depicted in Figure~\ref{fig:model_problem_presmoothing_illu}. This non-trivial initial phase distribution can be parametrized in a sharp-interface formulation in a straightforward manner, which is not the case for the required diffuse-interface formulation. Hence, this model problem serves as a test case to assess the preprocessing strategy proposed in Section~\ref{sec:preprocessing}, where we generate a valid initial phase field $\phi_0$ starting from a sharp-interface configuration of the initial phase distribution.

For the velocity $\mathbf{v}$ we utilize Dirichlet boundary conditions comprising a parabolic velocity profile in horizontal direction on the left boundary (inflow) of the domain $\Omega = [0,2] \times [0,1]$, 
\begin{align*}
 v_1\!\left((0,x_2)^T,t\right) = \bar{f} \frac{2}{3} \frac{x_2(1 - x_2)}{4}, \quad t \in [0,T],
\end{align*}
with mean flow $\bar{f}$, and value zero on the top and bottom boundary. On the right boundary (outlet) we employ zero Neumann boundary conditions for the velocity and prescribe a fixed value $p_b = 0$ for the pressure $p$. Other than that, for the pressure, the phase field $\phi$ and the chemical potential $\mu$ we use zero Neumann boundary conditions. The velocity is initialized with zero value in the entire domain.

\begin{figure}
 \centering
 \begin{subfigure}[b]{0.51\textwidth}
  \centering
  \includegraphics[keepaspectratio=true,width=\columnwidth]{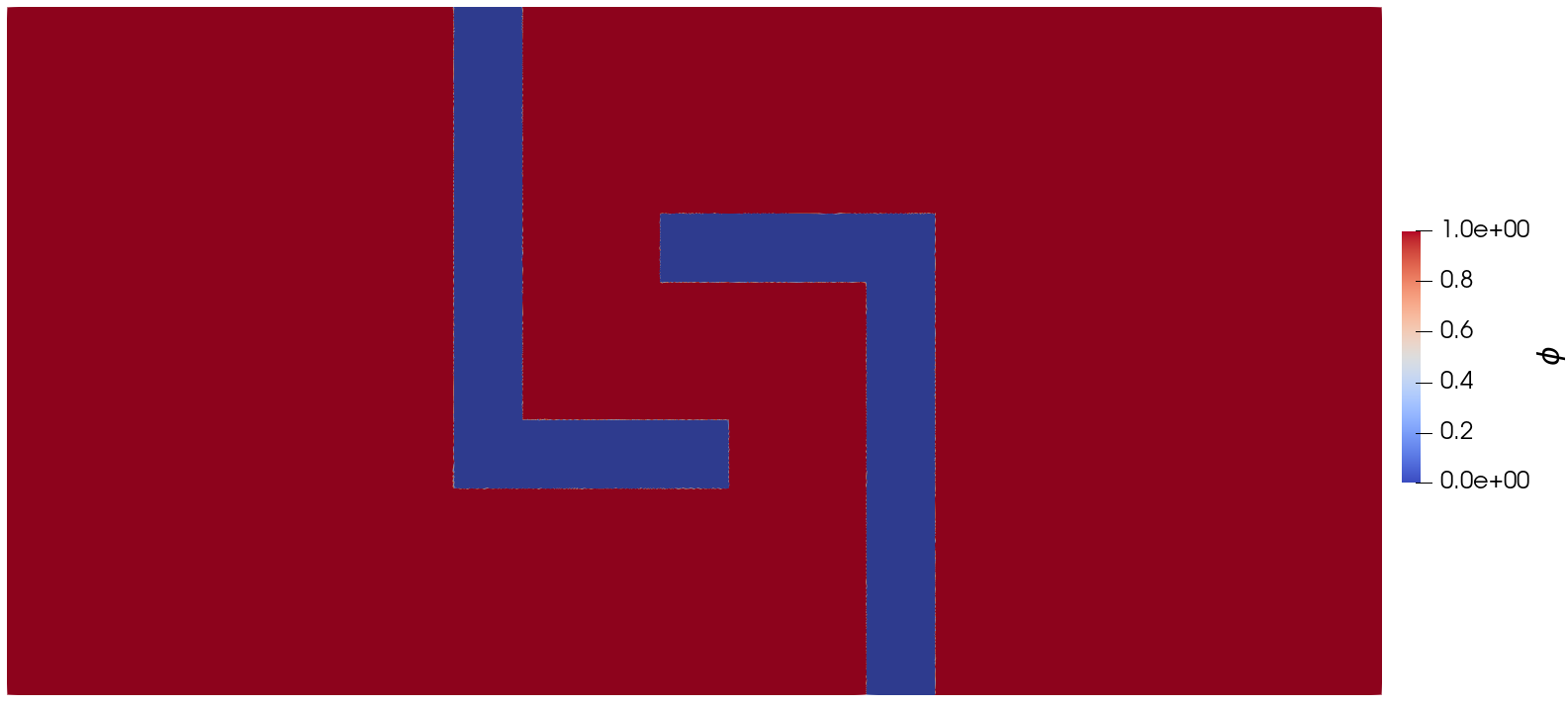}
 \end{subfigure}
 \begin{subfigure}[b]{0.465\textwidth}
  \centering
  \includegraphics[keepaspectratio=true,width=\columnwidth]{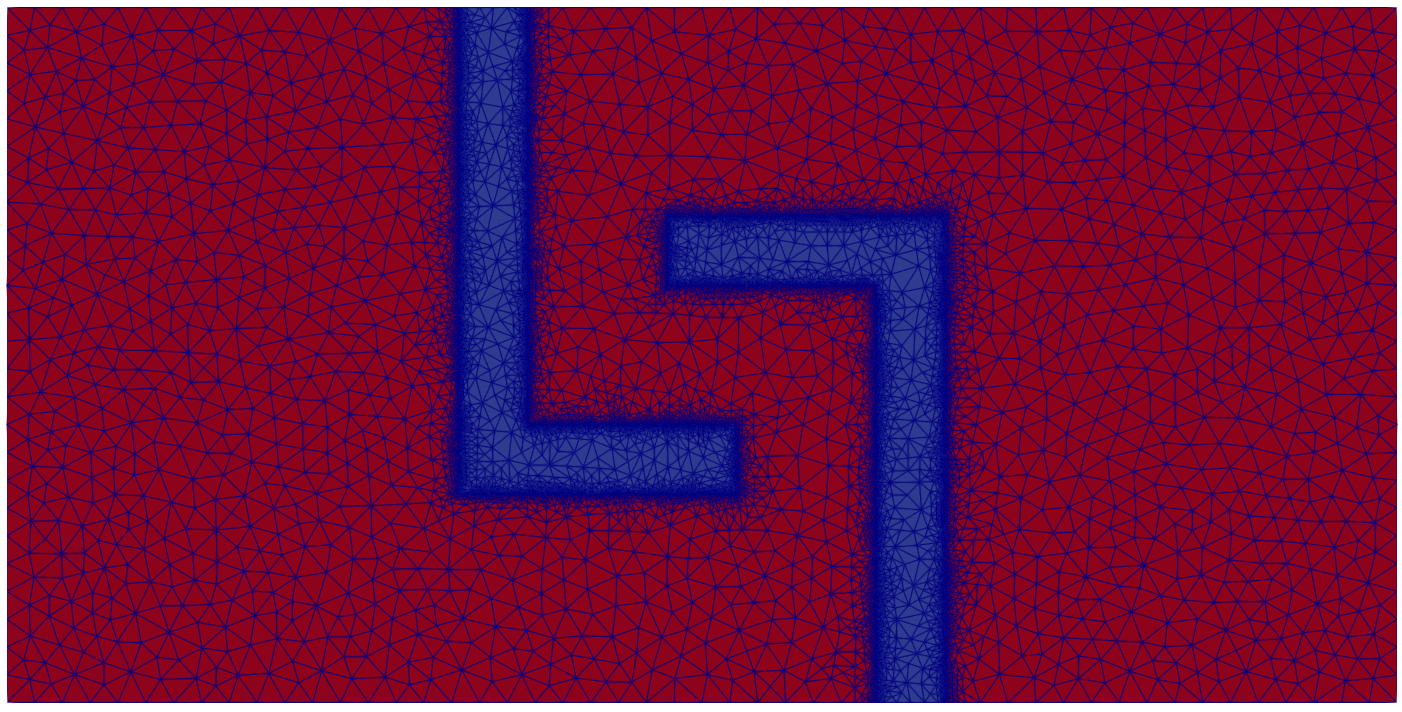}
 \end{subfigure}
 \caption{Sharp-interface configuration of the channel-flow inspired setup (left: implementation, right: refined grid).}
 \label{fig:model_problem_presmoothing_impl}
\end{figure}

The phase-field model parameters are given in Table~\ref{tab:parameters_cf}. Moreover, we set $\bar{f} = 0.1$ for the mean flow. For the second step of the preprocessing strategy introduced in Section~\ref{sec:preprocessing}, we employ $n_{\text{pre}} = 5$ iterations of the Cahn-Hilliard operator and set $M_{\text{pre}} = 10^3$ for the preprocessing mobility. In Figure~\ref{fig:model_problem_presmoothing_impl} we show the initial profile gained from our implementation along with the refined computational grid. The refinement indicator is based on the gradient of the phase-field variable $\phi$, so the grid refinement is applied predominantly in the interfacial region. In order to properly represent the sharp-interface configuration and to prevent spurious grid related effects a sufficiently high resolution of the interface is needed. In practice, we initially employ ten levels of refinement and we further refine the grid throughout the preprocessing strategy. 

\begin{table}
\centering
\caption{Phase-field model parameters for the channel-flow inspired setup. The given dimensions of the parameters are consistent with three spatial dimensions, where $\mathrm{T}$, $\mathrm{L}$ and $\mathrm{M}$ denote the dimension symbols of time, length and mass respectively.}
\begin{tabular}{llll}
 \hline
 Parameter & Symbol & Value & Dimension \\ 
 \hline
 Fluid density & $\rho$ & $10^3$ & $\mathrm L^{-3} \mathrm M$ \\ 
 Fluid viscosity & $\gamma$ & $10^{-3}$ & $\mathrm T^{-1} \mathrm L^{-1} \mathrm M$ \\ 
 Momentum dissipation in solid phase & $d_0$ & $10^3$ & $\mathrm T^{-1}$ \\ 
 Dissipation cut-off & $d_{\text{max}}$ & $0.9$ & -- \\
 Phase-field mobility & $M$ & $10^{-3}$ & $\mathrm T^{-1} \mathrm L^2$ \\ 
 Surface tension coefficient & $\sigma$ & $1.0$ & $\mathrm T^{-2} \mathrm M$ \\ 
 Diffuse-interface width & $\varepsilon$ & $6 \cdot 10^{-3}$ & $\mathrm L$ \\
 Phase-field regularization & $\delta$ & $12 \cdot 10^{-3}$ & -- \\
 Double-well potential modification & $\delta_{\text{dw}}$ & $8 \cdot 10^{-3}$ & -- \\
 Double-well potential modification & $\gamma_{\text{dw}}$ & $6 \cdot 10^{-3}$ & -- \\
 \hline
\end{tabular}
\label{tab:parameters_cf}
\end{table}

\begin{figure}
 \centering
 \begin{subfigure}[b]{0.495\textwidth}
  \centering
  \includegraphics[keepaspectratio=true,width=\columnwidth]{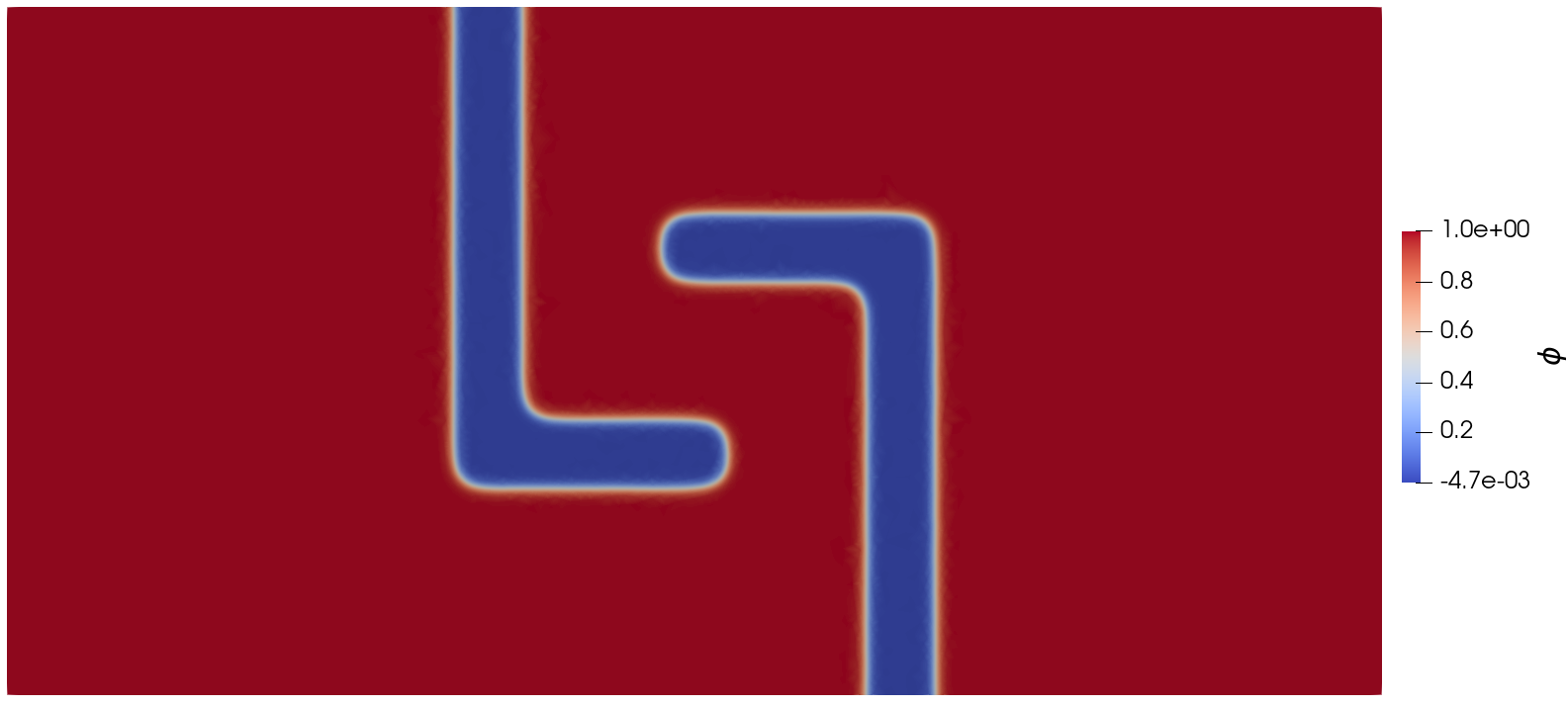}
 \end{subfigure}
 \begin{subfigure}[b]{0.495\textwidth}
  \centering
  \includegraphics[keepaspectratio=true,width=\columnwidth]{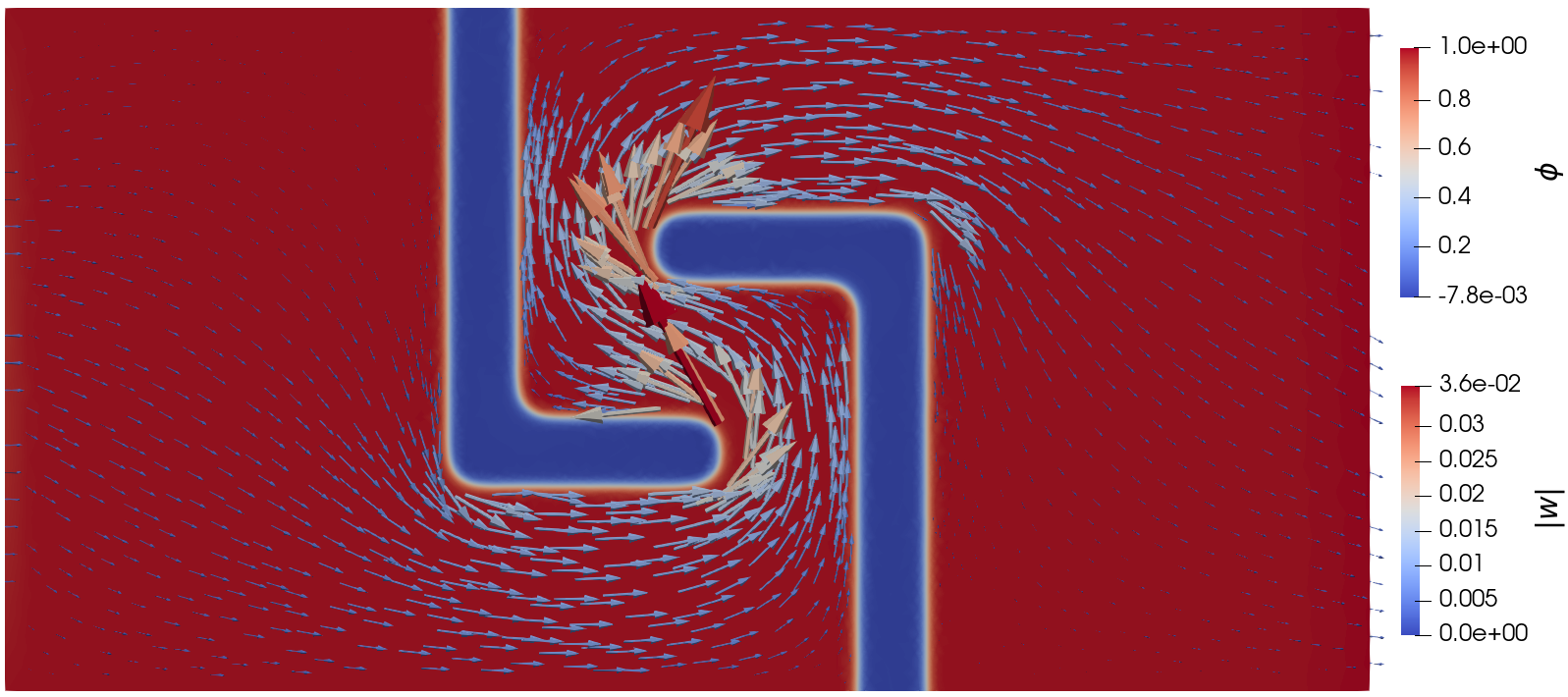}
 \end{subfigure}
 \caption{Diffuse-interface computation of the channel-flow inspired setup. We depict the result after the first time step on the left and the final one on the right. The arrows are scaled according to the magnitude of the conserved quantity $\mathbf{w} := \tilde \phi_f \mathbf{v}$. The color scaling is due to the data ranges presented on the right of the respective subfigure.}
 \label{fig:model_problem_presmoothing_comp}
\end{figure}

On the left of Figure~\ref{fig:model_problem_presmoothing_comp} we depict the diffuse profile obtained after the first time step comprising the preprocessing step. We see that the preprocessing strategy (in combination with an appropriately refined grid) yields a good diffuse-interface representation of the profile in Figure~\ref{fig:model_problem_presmoothing_impl}. We calculate $100$ time steps with time step size $\tau = 0.002$ resulting in final time $T = 0.2$. We depict the phase-field variable $\phi$ and the conserved quantity $\tilde \phi_f \mathbf{v}$ at the final time on the right of Figure~\ref{fig:model_problem_presmoothing_comp}. Firstly, we note that choosing a reasonably low mobility for the two-phase model evaluations we can satisfy the requirement of an almost static phase-field distribution. Secondly, the computed flow field shows the expected behavior of acceleration in the constriction built from the solid phase obstacles and deceleration afterwards.

These findings show the applicability of the two-phase model together with the preprocessing strategy for evaluating fluid flow in complex geometries.

\begin{figure}
 \centering
 \begin{subfigure}[b]{0.495\textwidth}
  \centering
  \includegraphics[keepaspectratio=true,width=\columnwidth]{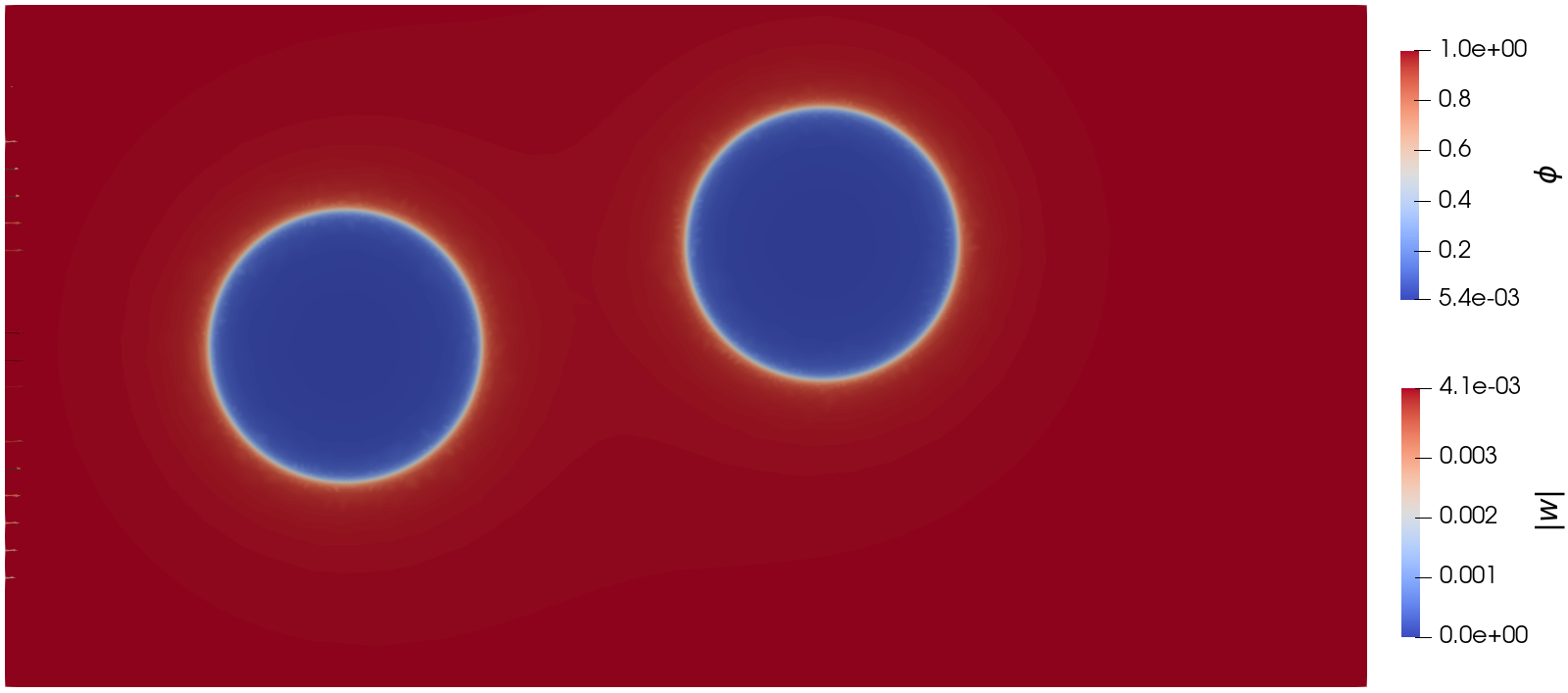}
 \end{subfigure}
 \centering
 \begin{subfigure}[b]{0.495\textwidth}
  \centering
  \includegraphics[keepaspectratio=true,width=\columnwidth]{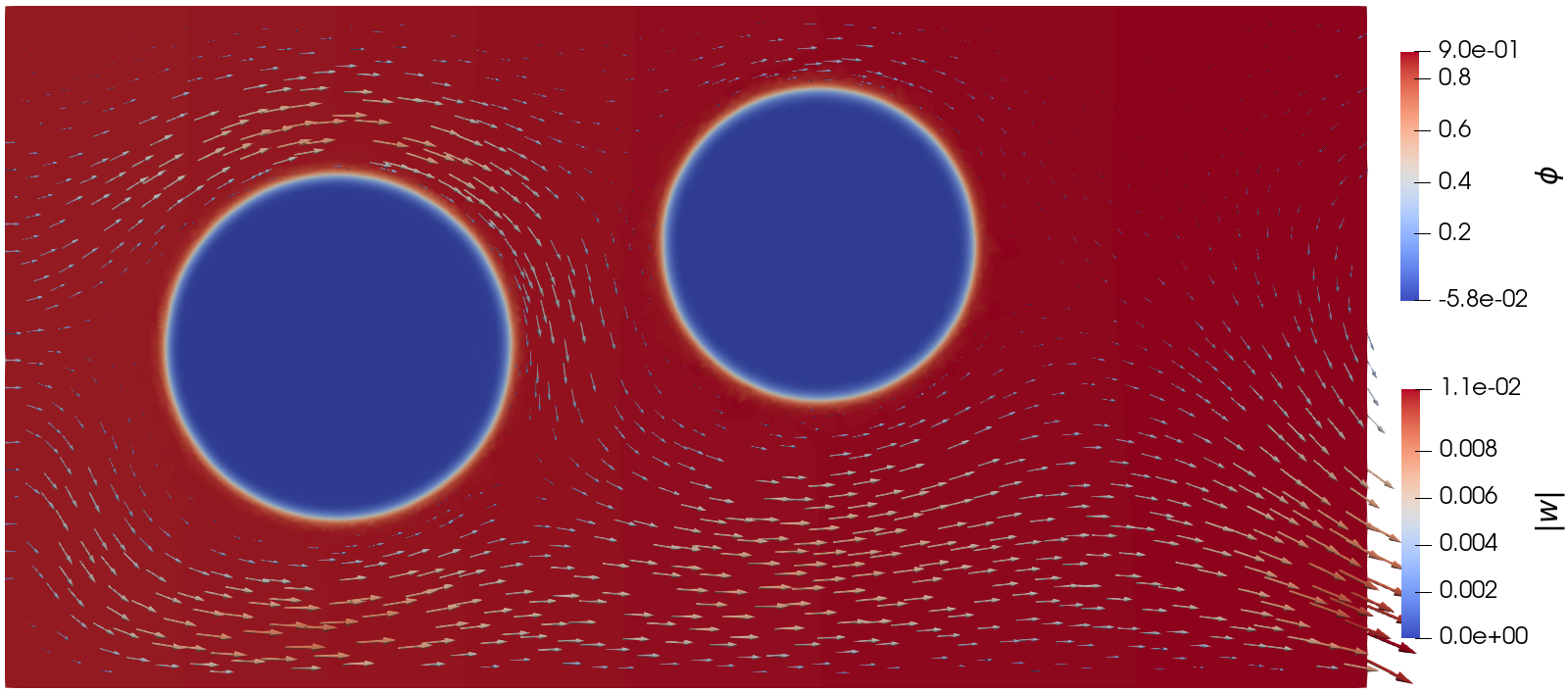}
 \end{subfigure}
 \begin{subfigure}[b]{0.495\textwidth}
  \centering
  \includegraphics[keepaspectratio=true,width=\columnwidth]{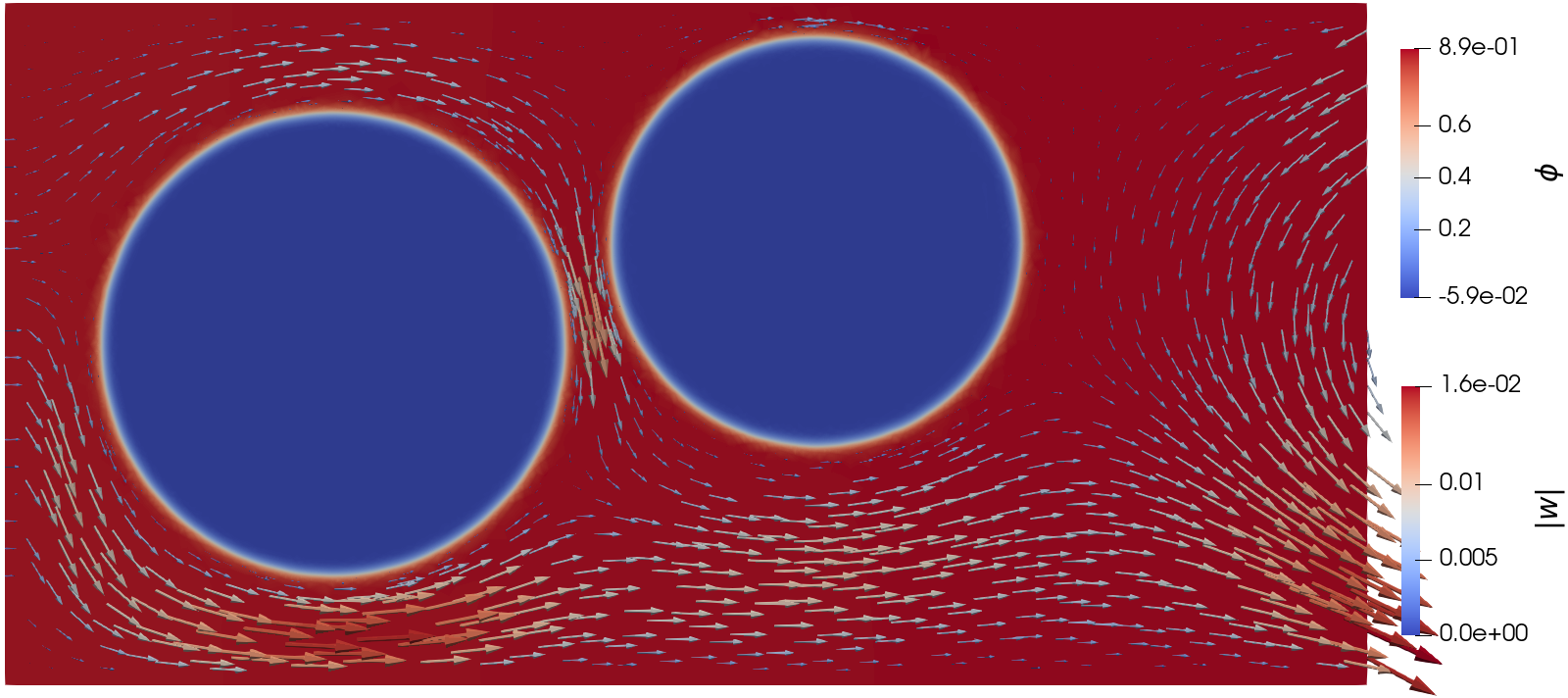}
 \end{subfigure}
 \begin{subfigure}[b]{0.495\textwidth}
 \centering
  \includegraphics[keepaspectratio=true,width=\columnwidth]{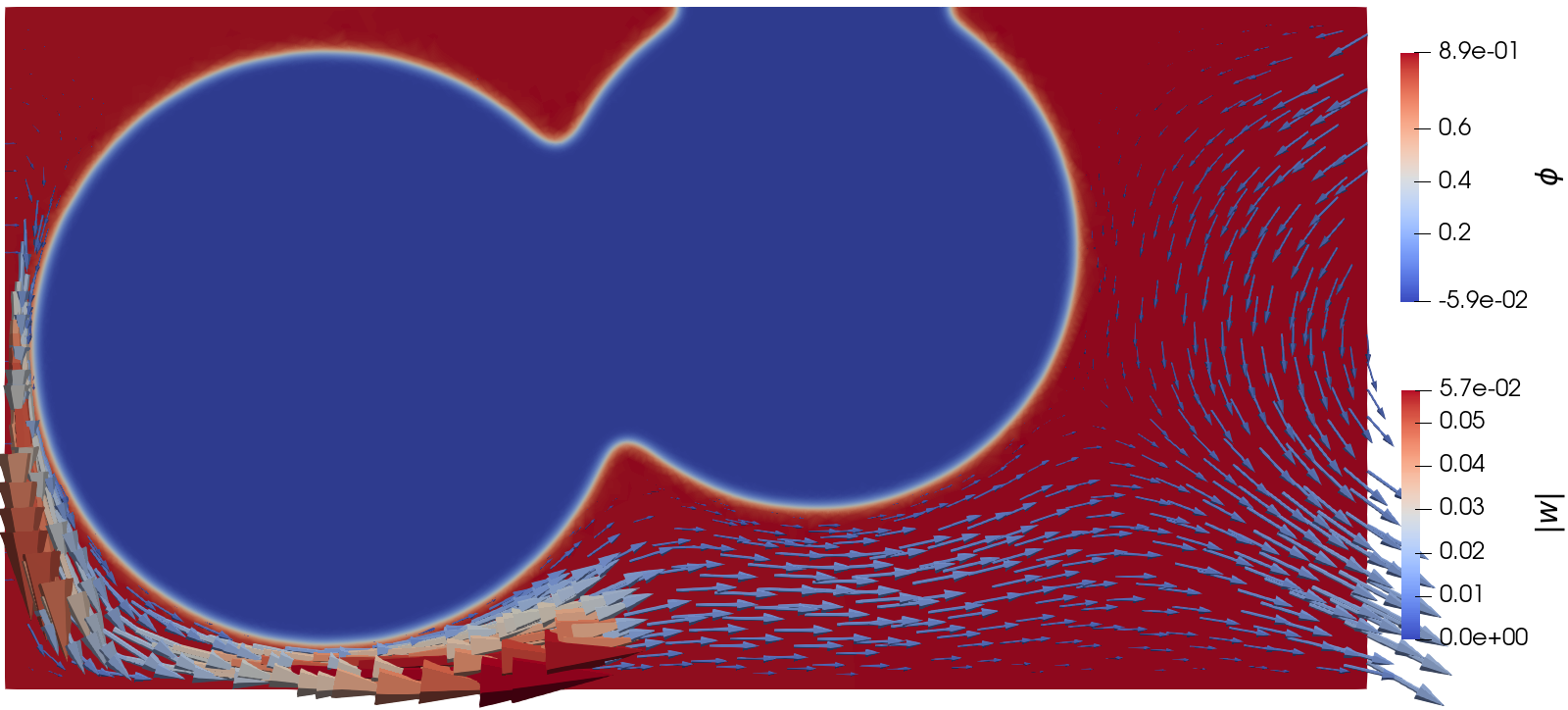}
 \end{subfigure}
 \begin{subfigure}[b]{0.495\textwidth}
 \centering
  \includegraphics[keepaspectratio=true,width=\columnwidth]{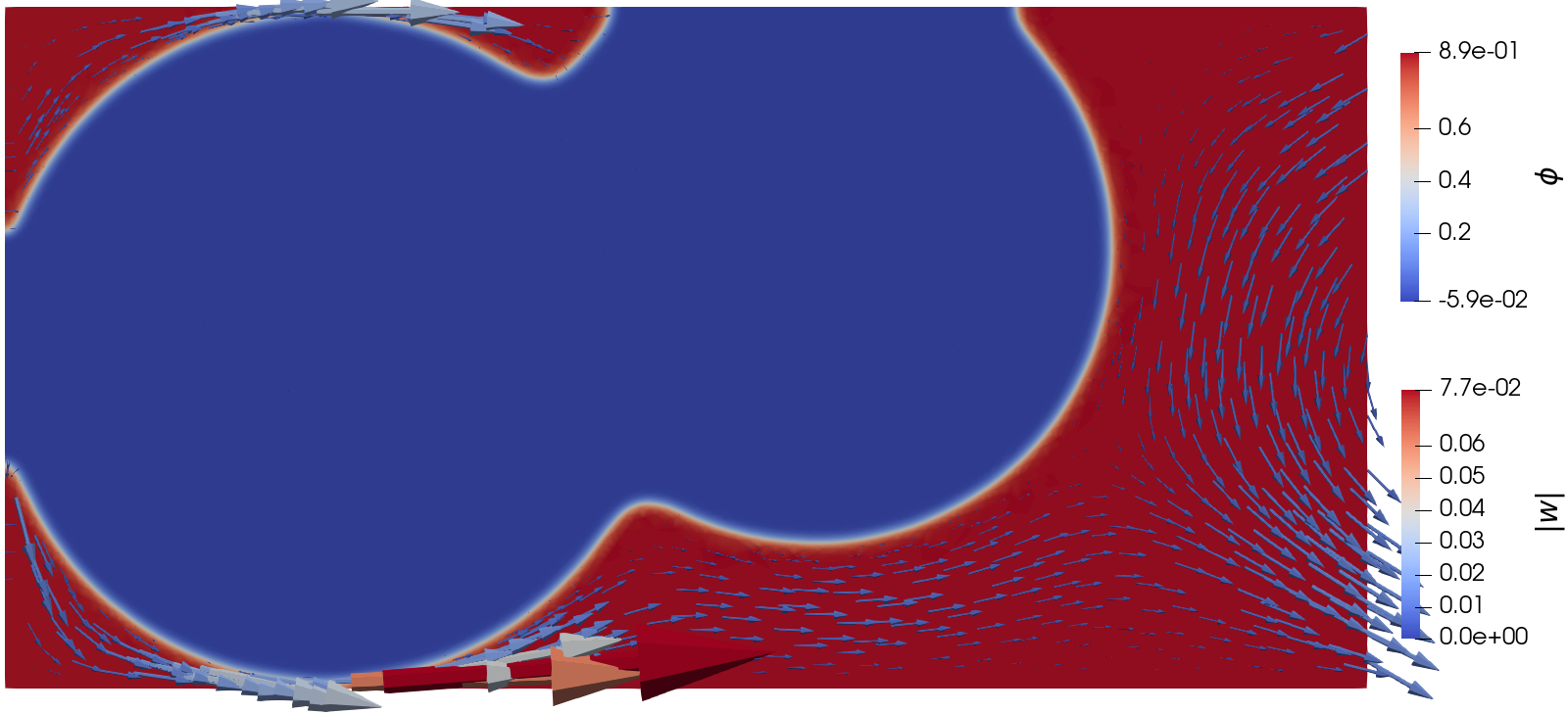}
 \end{subfigure}
 \caption{\revision{Results of the \revisionT{precipitating} reactive channel-flow experiment at time $t \in \{0, 1, 2, 3, 3.58\}$ (from top left to bottom right). The arrows are scaled according to the magnitude of the conserved quantity $\mathbf{w} := \tilde \phi_f \mathbf{v}$. The color scaling is due to the data ranges presented on the right of the respective subfigure.}}
 \label{fig:model_problem_react}
\end{figure}

\revision{
\subsection{Precipitation/dissolution effects} \label{sec:react}
In this section we consider the reactive model \eqref{eq:model} to capture precipitation and dissolution of ions at the interface between the solid and the liquid as well as the motion of the ions in the fluidic phase.

\subsubsection{Thermodynamic consistency}
Firstly, we evaluate the time evolution of the free energy $F[\phi(\cdot,t), \nabla\phi(\cdot,t), \mathbf{v}(\cdot,t), c(\cdot,t)]$ from \eqref{eq:free_energy}. To do so, we employ the lid-driven cavity inspired setup with the initial and boundary conditions from Section~\ref{sec:lid-driven}. We employ two circular solid inclusions in the flow with centers $\mathbf{x}_{1,c} = (0.5,0.5)^T$ and $\mathbf{x}_{2,c} = (1.65, 0.65)^T$ and radius $r_1 = 0.2 = r_2$. The phase-field model parameters are the same as in Table~\ref{tab:parameters_ldc}. Moreover, we employ the precipitation/dissolution model parameters as given in Table~\ref{tab:parameters_react}. We assume a constant initial concentration $c(\cdot,0) = 1.5$, along with zero Neumann boundary conditions for $c$. The resulting time evolution of the free energy is depicted in Figure~\ref{fig:model_problem_two_spheres_energy_react}. We again split the Ginzburg-Landau energy into its two constituents. For better illustration, we scale the kinetic energy by a factor of $10^4$, whose curve reflects the time step where the driving flow is stopped, marked by the dashed line. The total energy is decreasing in accordance with Theorem~\ref{theo:main}.

\begin{table}
\centering
\caption{\revision{Precipitation/dissolution} model parameters. The given dimensions of the parameters are consistent with three spatial dimensions, where $\mathrm{T}$, $\mathrm{L}$, $\mathrm{M}$ and $\mathrm{N}$ denote the dimension symbols of time, length, mass and amount of substance respectively.}
\begin{tabular}{llll}
 \hline
 Parameter & Symbol & Value & Dimension \\ 
 \hline
 Diffusivity & $D$ & $1.0$ & $\mathrm T^{-1} \mathrm L^2$ \\ 
 Solid concentration & $c^\ast$ & $2.0$ & $\mathrm L^{-3} \mathrm N$ \\ 
 Equilibrium concentration & $c_{\text{eq}}$ & $1.0$ & $\mathrm L^{-3} \mathrm N$ \\ 
 Reaction rate constant & $k_c$ & $0.1$ & \revision{$\mathrm T^{-1} \mathrm L$} \\ 
 \revisionT{Geometric scaling constant} & $\tilde \alpha$ & $0.03$ & $\mathrm M^{-1} \mathrm T \mathrm L^2$ \\ 
\hline
\end{tabular}
\label{tab:parameters_react}
\end{table}

\begin{figure}
 \centering
 \includegraphics[width=0.5\linewidth]{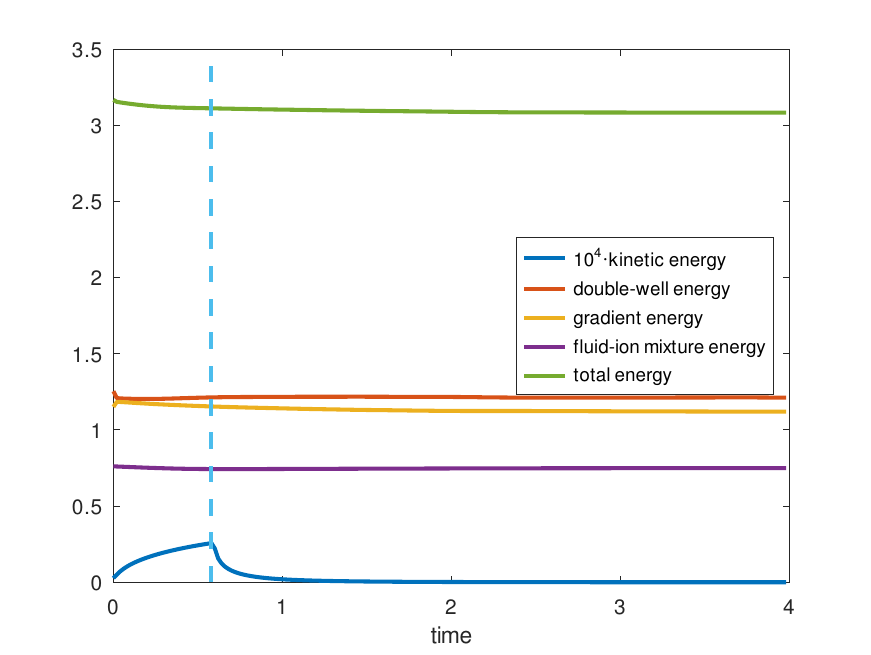}
 \caption{Energy plot for the lid-driven cavity inspired setup with two solid inclusions including precipitation/dissolution effects and stopped driving flow (dashed line).}
 \label{fig:model_problem_two_spheres_energy_react}
\end{figure}

\subsubsection{Clogging in a channel-flow setup}
}
We consider a channel flow with the initial and boundary conditions from Section~\ref{sec:channel_flow}. Moreover, we initialize $c(\cdot, 0) = 1$ in the entire domain, and inject a fixed ion concentration of $c_b = 1.5$ at the inlet. On the other boundaries we set zero Neumann boundary conditions for $c$. We employ two circular solid inclusions in the flow with centers $\mathbf{x}_{1,c} = (0.5,0.5)^T$ and $\mathbf{x}_{2,c} = (1.2, 0.65)^T$ and radius $r_1 = 0.2 = r_2$. The corresponding initial setup is depicted in the top left of Figure~\ref{fig:model_problem_react}. The phase-field model parameters are the same as in Table~\ref{tab:parameters_ldc} except for the regularization parameters $\delta = 0.09$, $\delta_{\text{dw}} = 0.075$ and $\gamma_{\text{dw}} = 0.06$, and the ones for the concentration are the same as in Table~\ref{tab:parameters_react} except for the reaction rate constant $k_c = 5$. We choose two solid inclusions of same radius in order to mitigate (geometric) effects of the Cahn-Hilliard evolution observed and discussed in the previous experiments and hence make the concentration related effects more pronounced. We calculate $200$ time steps with time step size $\tau = 0.02$ resulting in final time $T = 4$.

The results at different times are presented in Figure~\ref{fig:model_problem_react}.
As we inject a higher ion concentration, we expect net precipitation. This can be seen in the top right of Figure~\ref{fig:model_problem_react}, where the solid inclusion closer to the inlet has grown larger. The inclusions coalesce, form a new bulk that grows further, reaches the top boundary (middle left and right), and eventually also the left and bottom boundary (bottom). Hence, we observe clogging and stop the simulation. \revisionU{Note that the computed upper bound of the phase-field variable $\phi$ is considerably smaller than one indicating that some precipitation effects are still present in the fluid.}

\begin{figure}
 \centering
 \begin{subfigure}[b]{0.495\textwidth}
  \centering
  \includegraphics[keepaspectratio=true,width=\columnwidth]{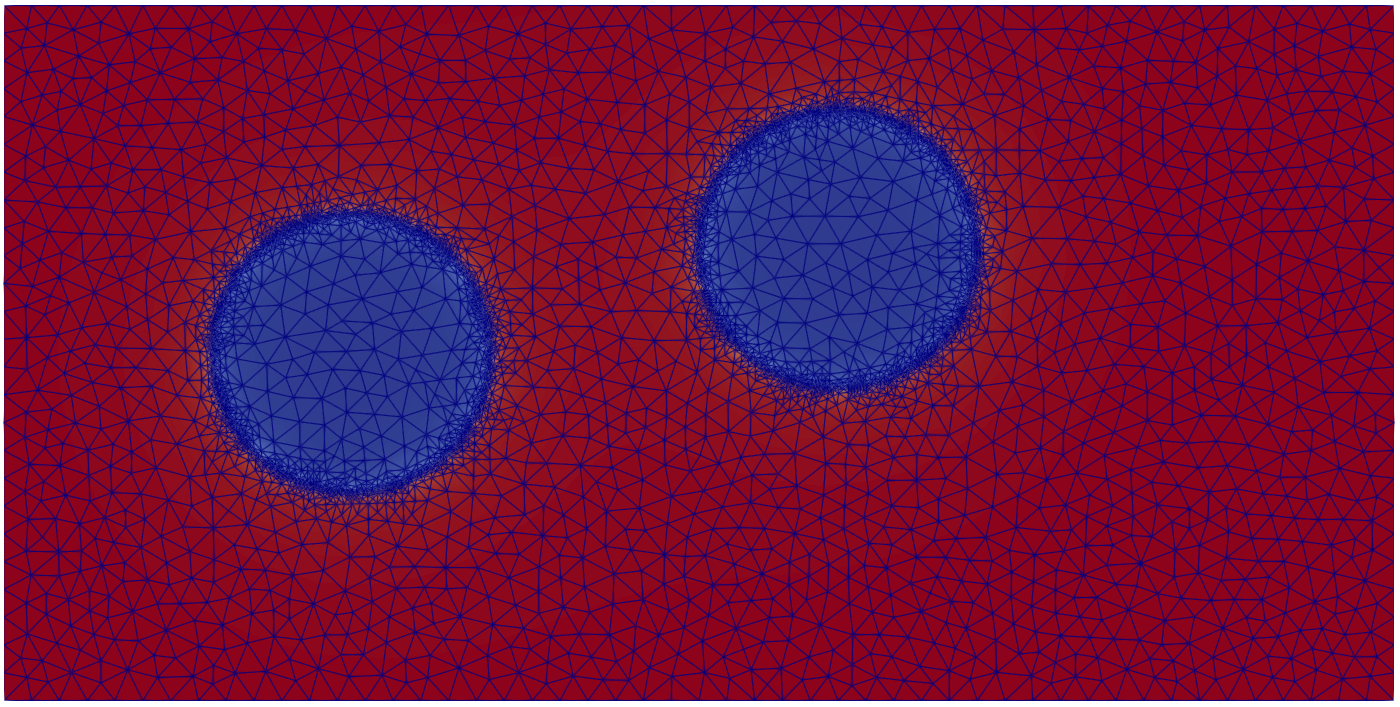}
 \end{subfigure}
 \begin{subfigure}[b]{0.495\textwidth}
  \centering
  \includegraphics[keepaspectratio=true,width=\columnwidth]{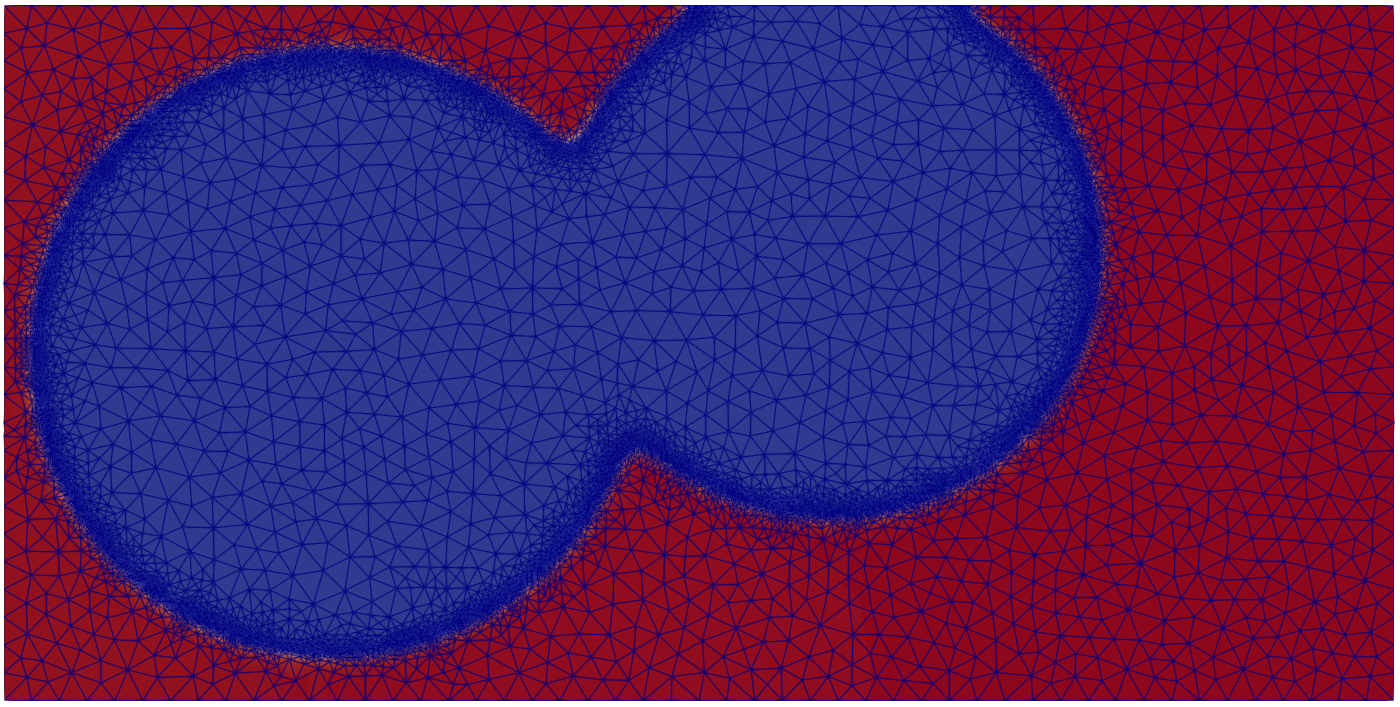}
 \end{subfigure}
 \caption{Refined grid at time $t \in \{0,3\}$.}
 \label{fig:model_problem_react_refinement}
\end{figure}

The distinct temporal evolution of the interfacial layer cannot be resolved without  the use of adaptive grid refinement. As an impression of the dynamical change of the grid we show two grid instances in Figure~\ref{fig:model_problem_react_refinement}.

\revisionT{
\subsubsection{Dissolution in a channel-flow setup}
We consider a channel flow with the initial and boundary conditions from Section~\ref{sec:channel_flow}. In order to assess the model's capability to capture the dissolution of ions, we initialize $c(\cdot, 0) = 1$ in the entire domain, and inject a fixed ion concentration of $c_b = 0.5$ at the inlet. On the other boundaries we set zero Neumann boundary conditions for $c$. We employ one circular solid inclusion in the flow with center $\mathbf{x}_c = (0.5,0.5)^T$ and radius $r = 0.2$. The corresponding initial setup is depicted in the top left of Figure~\ref{fig:model_problem_react_2}. The phase-field model parameters are the same as in Table~\ref{tab:parameters_ldc}, and the ones for the concentration are the same as in Table~\ref{tab:parameters_react} except for the diffusivity $D = 0.1$. For the non-reactive case we have observed (geometric) shrinkage effects of the Cahn-Hilliard evolution in Section~\ref{sec:num_thm}. To distinctly assess the dissolution process, we calculate the experiment for varying values of the reaction rate constant $k_c \in \{0, 0.1, 1\}$. We compute $200$ time steps with time step size $\tau = 0.02$ resulting in final time $T = 4$.

We depict the results for the different reaction rate constants $k_c$ at final time $T = 4$ in Figure~\ref{fig:model_problem_react_2}. Note that the case $k_c = 0$ corresponds to a reaction rate $r(c) = 0$ and hence the resulting reaction term $R = - \frac{q(\phi)}{\varepsilon} \tilde \alpha \sigma \mu$ merely accounts for curvature related effects. As expected, we observe that the higher $k_c$ is chosen the more dissolution occurs.

\begin{figure}
 \centering
 \begin{subfigure}[b]{0.495\textwidth}
  \centering
  \includegraphics[keepaspectratio=true,width=\columnwidth]{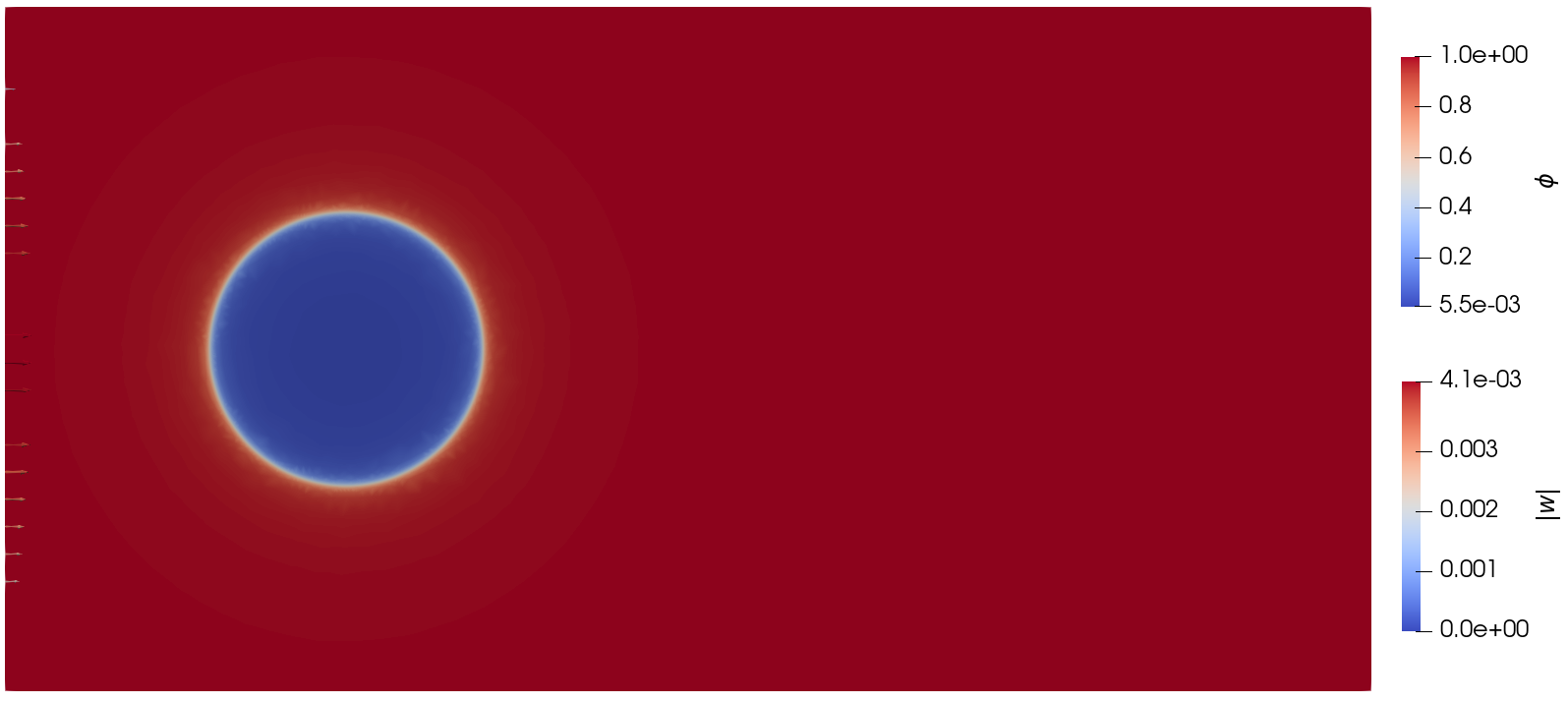}
 \end{subfigure}
 \centering
 \begin{subfigure}[b]{0.495\textwidth}
  \centering
  \includegraphics[keepaspectratio=true,width=\columnwidth]{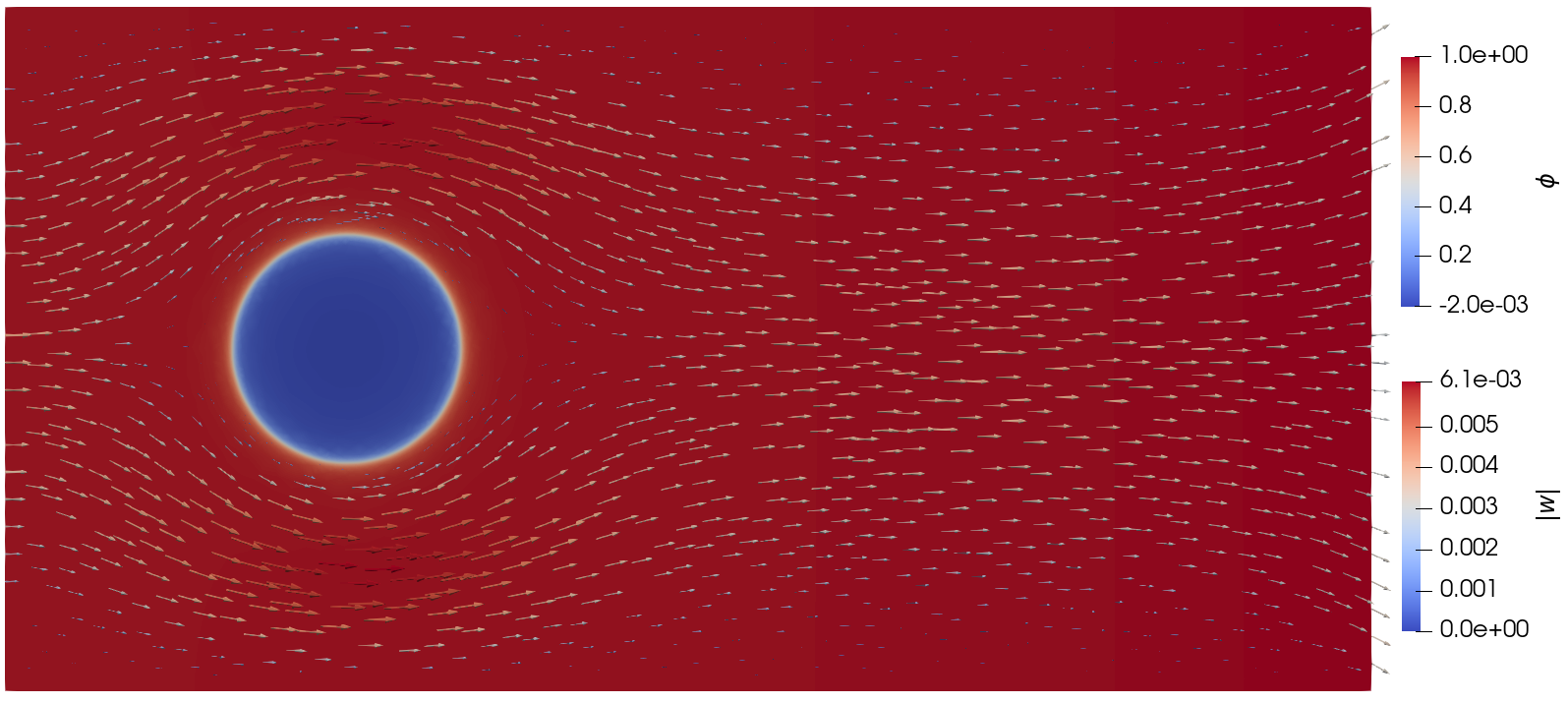}
 \end{subfigure}
 \begin{subfigure}[b]{0.495\textwidth}
  \centering
  \includegraphics[keepaspectratio=true,width=\columnwidth]{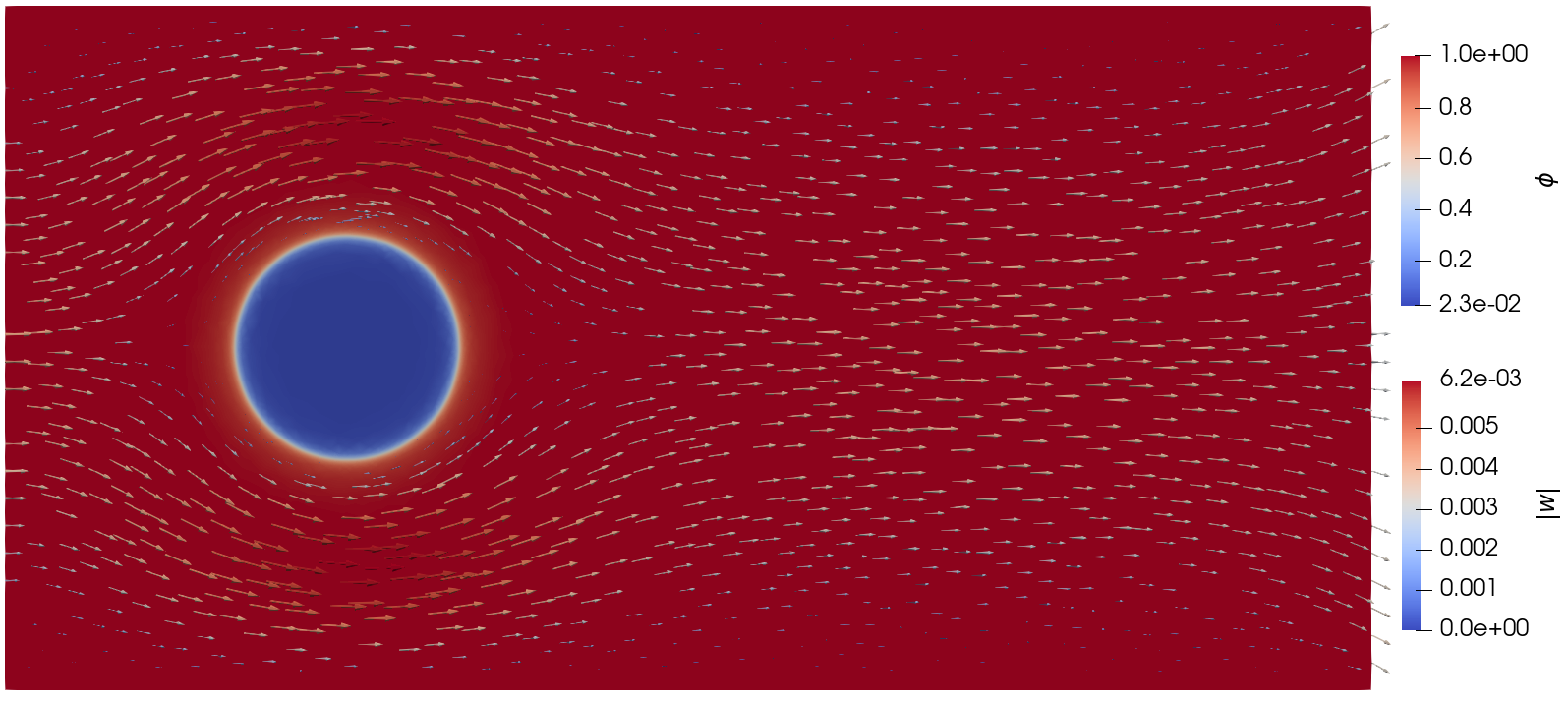}
 \end{subfigure}
 \begin{subfigure}[b]{0.495\textwidth}
 \centering
  \includegraphics[keepaspectratio=true,width=\columnwidth]{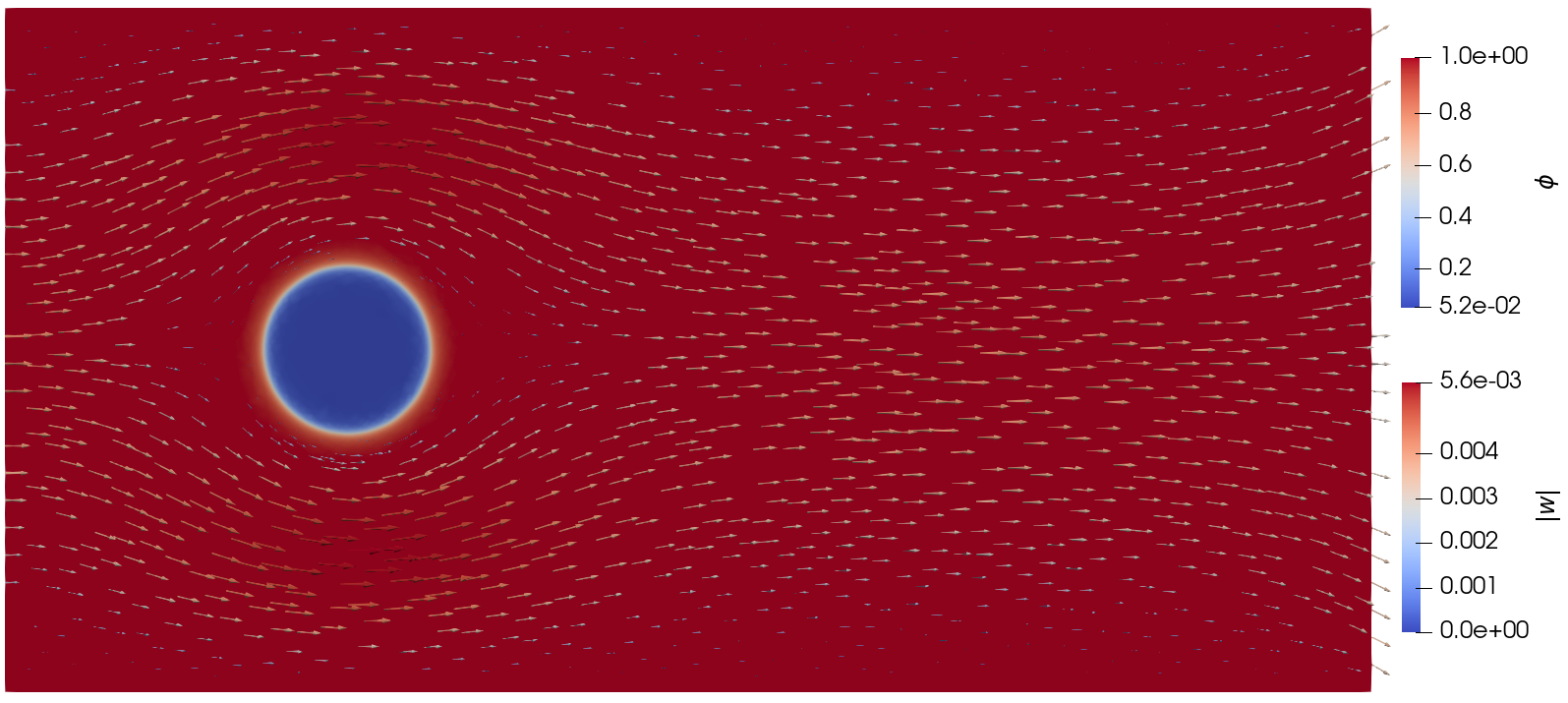}
 \end{subfigure}
 \caption{\revisionT{Results of the dissolving reactive channel-flow experiment at time $T = 4$ with varying reaction rate constant $k_c \in \{0, 0.1, 1\}$ (in the order top right, bottom left, bottom right). The arrows are scaled according to the magnitude of the conserved quantity $\mathbf{w} := \tilde \phi_f \mathbf{v}$. 
 The color scaling is due to the data ranges presented on the right of the respective subfigure.
 }}
 \label{fig:model_problem_react_2}
\end{figure}
}

\section{Conclusions} \label{sec:concl}
We presented a  fully-discrete finite element based   method to solve the initial boundary value problem for 
an incompressible Cahn-Hilliard Navier-Stokes system that
governs the evolution of a fluid-solid system under geometric \revision{and reactive} effects.  We employed a tailor-made first-order time-stepping method \revisionU{that renders the fully-discrete approach to be energy-stable under appropriate assumptions on the discrete solution.}  
For the energy stability we exploited the convex-concave splitting of the non-convex free energy function and a tailored  implicit-explicit discretization of the corresponding fluxes that allow to mimic the proof of the 
energy decay on the continuous level.
We  have shown the applicability of the method for a series of numerical experiments related to the lid-driven cavity problem and a  channel-flow setup \revision{both for non-reactive and reactive flow problems with precipitation/dissolution effects}.  
\\
Further research is devoted to finding higher-order in time methods to match the accuracy of the high-order space discretization. For  a conditionally energy-stable second-order in time scheme which introduces time-step restrictions we refer to \cite{Khanwale2022} whose scope is limited to the case of two fluid phases. We further mention \cite{Brunk2023} where a second order structure-preserving scheme is presented for the case of two fluid phases with constant mass density. The same applies for the use of discontinuous Galerkin methods \revision{\cite{Giesselmann2015} and the extension of the partitioned solution strategy to the reactive model.}

\section*{Acknowledgement}
Funded by Deutsche Forschungsgemeinschaft (DFG, German Research Foundation) under Germany’s Excellence Strategy (Project Number 390740016 – EXC 2075) and the Collaborative Research Center 1313
(Project Number 327154368 – SFB 1313). We acknowledge the support by the Stuttgart Center for Simulation Science (SC SimTech). We thank the anonymous reviewers for their thorough proofreading and valuable comments that helped us to improve the manuscript.

\appendix
\section{Vector calculus}\label{appendix}
\renewcommand{\thesection}{A}
\begin{lemma}\label{auxlemma}
 For (weakly) differentiable vector fields $\mathbf{a}$, $\mathbf{b}$ and $\mathbf{c}$ the identities
 \begin{align}
  &\nabla \cdot (\mathbf{a} \otimes \mathbf{b}) = (\nabla \cdot \mathbf{a}) \mathbf{b} + (\mathbf{a} \cdot \nabla) \mathbf{b}, \qquad \mathbf{b} \cdot \big((\mathbf{a} \cdot \nabla)\mathbf{b}\big) = \frac{1}{2} (\mathbf{a} \cdot \nabla) |\mathbf{b}|^2, \label{eq:aux_1} \\
  &\mathbf{a} \cdot (\mathbf{a} - \mathbf{b}) = \frac{1}{2} |\mathbf{a}|^2 - \frac{1}{2} |\mathbf{b}|^2 + \frac{1}{2} |\mathbf{a} - \mathbf{b}|^2, \label{eq:aux_2} \\
  &((\mathbf{a} \cdot \nabla) \mathbf{b}) \cdot \mathbf{c} = \frac{1}{2} \mathbf{a} \cdot \nabla(\mathbf{b} \cdot \mathbf{c}) + \frac{1}{2} \mathbf{a} \cdot ((\nabla \mathbf{b}) \mathbf{c}) - \frac{1}{2} \mathbf{a} \cdot ((\nabla \mathbf{c}) \mathbf{b}) \label{eq:aux_3}
 \end{align}
 hold.
\end{lemma}

\bibliographystyle{unsrt}
\bibliography{bibliography}

\end{document}

%% file: macro_CR.tex
\usepackage{epstopdf}
\usepackage{stmaryrd} 
\usepackage{fixmath} 
\usepackage{pgf}
\usepackage{wasysym}

\newcommand{\ds}{\displaystyle}

\newcommand{\vecstyle}[1]{{\bf{#1}}}

\newcommand{\vecx}{\vecstyle{x}}

\newcommand{\vecn}{\vecstyle{n}}

\newcommand{\vecv}{\vecstyle{v}}
\newcommand{\setstyle}[1]{{\mathbb #1}}

\def\setR {\setstyle{R}}

